\documentclass{birkmult}
\usepackage{amssymb,amsmath,epsfig,graphics}

\title{GKZ Hypergeometric Structures}
\date{}
\author{Jan Stienstra}
\address{Mathematisch Instituut, Universiteit Utrecht, the Netherlands} \email{stien{@}{math.uu.nl}}
\newtheorem{theorem}{Theorem}
\newtheorem{lemma}{Lemma}
\newtheorem{proposition}{Proposition}
\newtheorem{corollary}{Corollary}
\newtheorem{definition}{Definition}

\newcommand{\CC}{{\mathbb C }}

\newcommand{\LL}{{\mathbb L }}
\newcommand{\MM}{{\mathbb M }}
\newcommand{\NN}{{\mathbb N }}
\newcommand{\PP}{{\mathbb P }}
\newcommand{\QQ}{{\mathbb Q }}
\newcommand{\RR}{{\mathbb R }}
\newcommand{\ZZ}{{\mathbb Z }}
\newcommand{\TT}{{\mathbb T }}

\newcommand{\II}{{\mathbb I }}

\newcommand{\gG}{{\Gamma }}
\newcommand{\gam}{{\gamma }}
\newcommand{\gl}{{\lambda }}
\newcommand{\ga}{{\alpha }}
\newcommand{\gs}{{\sigma }}
\newcommand{\gf}{{\varphi }}
\newcommand{\uga}{{\underline{\alpha }}}
\newcommand{\ugamma}{{\underline{\gamma}}}
\newcommand{\uz}{{\underline{\mf{0}}}}
\newcommand{\eps}{{\varepsilon }}

\newcommand{\eul}{\raisebox{.05ex}{${}_{\Upsilon }$}}

\newcommand{\gA}{{\mathcal A }}
\newcommand{\sA}{{\mathsf A }}

\newcommand{\bu}{\mathbf{u}}

\newcommand{\va}{{\mathbf a }}
\newcommand{\vb}{\mathbf{b}}
\newcommand{\vc}{\mathbf{c}}
\newcommand{\vh}{\mathbf{h}}
\newcommand{\vm}{\mathbf{m}}
\newcommand{\vp}{\mathbf{p}}
\newcommand{\vx}{\mathbf{x}}
\newcommand{\vz}{\mathbf{z}}
\newcommand{\vv}{\mathbf{v}}
\newcommand{\vu}{\mathbf{u}}
\newcommand{\vw}{\mathbf{w}}
\newcommand{\hw}{\widehat{\mathbf{w}}}

\newcommand{\ull}{\underline{\ell}}
\newcommand{\ueps}{\underline{\eps}}
\newcommand{\oeps}{\overline{\eps}}
\newcommand{\ode}{\overline{\delta}}

\newcommand{\ba}{{\mathbf a }}
\newcommand{\bb}{{\mathbf b }}

\newcommand{\ci}{{\mathsf i }}

\newcommand{\cB}{{\mathcal B }}
\newcommand{\cC}{{\mathcal C }}
\newcommand{\cH}{{\mathcal H }}
\newcommand{\cI}{{\mathcal I }}

\newcommand{\cU}{{\mathcal U }}
\newcommand{\cR}{{\mathcal R }}

\newcommand{\cP}{{\mathcal P }}
\newcommand{\cV}{{\mathcal V }}
\newcommand{\cO}{{\mathcal O }}
\newcommand{\cF}{{\mathcal F }}
\newcommand{\oR}{\overline{\mathcal R }}
\newcommand{\gT}{{\mathcal T }}

\newcommand{\hp}{{\widehat{\pi} }}

\newcommand{\rank}{\mathrm{ rank}\,}
\newcommand{\spec}{\mathrm{ Spec}\,}
\newcommand{\Hom}{\mathrm{ Hom}\,}

\newcommand{\mf}[1]{\mathsf{#1}}

\newcommand{\modquot}[2]{\mbox{\raisebox{.2em}{$#1$}\hspace{-2mm} 
{ / }\hspace{-2mm} \raisebox{-.2em}{$#2$}}}

\newcommand{\VE}[1]{\left[\!\!\begin{array}{c} #1\end{array}\!\!\right]}

\begin{document}
\begin{abstract}
This text is based on lectures by the author in the Summer School
\emph{Algebraic Geometry and Hypergeometric Functions} 
in Istanbul in June 2005.
It gives a review of some of the basic aspects of the theory of
hypergeometric structures of Gelfand, Kapranov and Zelevinsky, 
including Differential Equations, Integrals and Series,
with emphasis on the latter. The Secondary Fan is constructed and subsequently
used to describe the `geography' of the domains of convergence
of the $\gG$-series.
A solution to certain Resonance Problems is presented and applied
in the context of Mirror Symmetry. Many examples and some exercises
are given throughout the paper.
\end{abstract}

\maketitle

\section{Introduction}
GKZ stands for \textit{Gelfand, Kapranov} and \textit{Zelevinsky},
who discovered 
fascinating generalizations of the classical hypergeometric structures of Euler, Gauss, Appell, Lauricella, Horn \cite{gkz1,gkz2,gkz4}. 
The main ingredient for these new hypergeometric structures is a finite subset $\gA\subset\ZZ^{k+1}$ which generates $\ZZ^{k+1}$ as an abelian group and for which there exists a group homomorphism 
$h:\ZZ^{k+1}\rightarrow\ZZ$ such that \mbox{$h(\gA)=\{1\}$.}
The latter condition means that $\gA$ lies in a $k$-dimensional affine hyperplane in $\ZZ^{k+1}$.
\begin{figure}
\begin{center}
\begin{picture}(300,130)(30,-10)
\put(35,20){
\begin{picture}(100,80)(0,0)
\put(0,0){\line(1,0){40}}
\put(0,0){\line(0,1){40}}
\put(40,40){\line(-1,0){40}}
\put(40,40){\line(0,-1){40}}
\put(0,0){\circle*{5}}
\put(40,0){\circle*{5}}
\put(0,40){\circle*{5}}
\put(40,40){\circle*{5}}
\put(7,90){\makebox{Gauss}}
\end{picture}}
\put(140,20){
\begin{picture}(100,80)(0,0)
\put(0,0){\line(1,0){40}}
\put(0,0){\line(0,1){40}}
\put(40,40){\line(-1,0){40}}
\put(40,40){\line(0,-1){40}}
\put(60,20){\line(0,1){40}}
\put(40,0){\line(1,1){20}}
\put(40,40){\line(1,1){20}}
\put(0,40){\line(3,1){60}}
\multiput(0,0)(15,5){4}{\line(3,1){10}}
\put(0,0){\circle*{5}}
\put(40,0){\circle*{5}}
\put(0,40){\circle*{5}}
\put(40,40){\circle*{5}}
\put(60,20){\circle*{5}}
\put(60,60){\circle*{5}}
\put(7,90){\makebox{Appell $F_1$}}
\end{picture}}
\put(260,20){
\begin{picture}(100,80)(0,0)
\put(0,30){\line(1,0){45}}
\put(45,30){\line(1,1){22}}
\multiput(0,30)(15,15){2}{\line(1,1){10}}
\multiput(25,55)(15,0){3}{\line(1,0){10}}
\put(0,30){\line(1,2){35}}
\put(0,30){\line(3,-4){35}}
\put(70,55){\line(-1,-2){35}}
\put(69,55){\line(-3,4){35}}
\put(46,30){\line(-1,6){12}}
\put(46,30){\line(-1,-4){11}}
\multiput(25,55)(6,-36){2}{\line(1,-6){3}}
\multiput(25,55)(5,20){2}{\line(1,4){4}}
\put(0,30){\circle*{5}}
\put(46,30){\circle*{5}}
\put(69,55){\circle*{5}}
\put(25,55){\circle*{5}}
\put(35,-16){\circle*{5}}
\put(35,100){\circle*{5}}
\put(-20,90){\makebox{Appell $F_4$}}
\end{picture}}
\end{picture}
\end{center}
\caption{\label{fig:classical examples}}
\end{figure}
Figure \ref{fig:classical examples} shows $\gA$ (the black dots)
sitting in this hyperplane for some classical hypergeometric structures.
In \cite{gkz2,gkz4} these new structures were called \emph{$\gA$-hypergeometric systems}. Nowadays many authors call them \emph{GKZ hypergeometric systems}.
The original name indeed seems somewhat unfortunate,
since \emph{$\gA$-hypergeometric} sounds negative, like   
$'\alpha\gamma\varepsilon o\mu\varepsilon\tau\rho\iota\tau o\varsigma
\;\mu\eta\;$ $'\varepsilon\iota\sigma\iota\tau\omega$ 
(a non-geometer should not enter),
written over the entrance of Plato's academy and in the logo of the 
American Mathematical Society.
Besides the set $\gA$
the construction of GKZ hypergeometric structures requires a vector $\vc\in\CC^{k+1}$.

In these notes we report on the basic theory of GKZ hypergeometric structures
and show how the traditional aspects
\textit{differential equations, integrals, series} are attached to the data $\gA$, $\vc$.
In Section \ref{section GKZ DE} we introduce the GKZ differential 
equations and give examples of GKZ hypergeometric integrals.
In Section \ref{section:gamma series} we discuss GKZ hypergeometric series (so-called $\gG$-series).
We have put details of the GKZ theory for Lauricella's $F_D$ together in Section \ref{section:FD}, so that the reader can 
compare results and view-points on $F_D$ for various
lectures in this School (e.g. \cite{loo}).

The beautiful insight of Gelfand, Kapranov and Zelevinsky was that hypergeometric structures greatly simplify if one introduces extra variables and balances this with an appropriate torus action. More precisely the variables in GKZ theory are the natural 
coordinates on the space $\CC^\gA:=\mathrm{Maps}(\gA,\CC)$ of maps 
from $\gA$ to $\CC$. The torus $\TT^{k+1}:=\mathrm{Hom}(\ZZ^{k+1},\CC^*)$ of 
group homomorphisms from $\ZZ^{k+1}$ to $\CC^*$, acts naturally on 
$\CC^\gA$ and on functions on $\CC^\gA$: 
\emph{for $\sigma\in\TT^{k+1}$ and $\Phi:\CC^\gA\rightarrow\CC$}
\begin{equation}\label{eq:torus action on functions}
(\sigma\cdot\vu)(\va)\,=\, \sigma(\va)\vu(\va)\,,\qquad
(\Phi\cdot\sigma)(\vu)\,=\,\Phi(\sigma\cdot\vu)\,,
\qquad\forall\, \va\in\gA, \;\forall\, \vu\in\CC^\gA\,.
\end{equation}
The GKZ hypergeometric functions associated with $\gA$ and $\vc$
are defined on open domains in $\CC^\gA$, but they 
are not invariant under the  action of $\TT^{k+1}$, unless $\vc=0$.
Rather, for $\vc\in\ZZ^{k+1}$ they transform according to the character
of $\TT^{k+1}$ given by $\vc$. For $\vc\not\in\ZZ^{k+1}$ there is
only an infinitesimal analogue of this transformation behavior,
encoded in one part of the GKZ system of differential equations
(see (\ref{eq:GKZ2})).
On the other hand, the quotient of
any two GKZ hypergeometric functions with a common domain of definition associated with $\gA$ and $\vc$
is always $\TT^{k+1}$-invariant (see (\ref{eq:invariant quotient})).

The important role of the $\TT^{k+1}$-action in GKZ hypergeometric structures motivates a study of the orbit space.
Without going into details, 
this can be described as follows. First take the complement of
the coordinate hyperplanes
$(\CC^*)^\gA:=\mathrm{Maps}(\gA,\CC^*)=
\{\vu\in\CC^\gA\:|\: \vu(\va)\neq 0,\;\forall \va\in\gA\}$.
The above action of $\TT^{k+1}$ preserves this set.
In fact $(\CC^*)^\gA$ is a complex torus and 
$\TT^{k+1}$ can be identified
with a subtorus, acting by left multiplication.
The quotient is the torus
\begin{equation}\label{eq:quotient torus}
\modquot{(\CC^*)^\gA}{\TT^{k+1}}\,=\,\Hom(\LL,\CC^*),
\end{equation}
where $\LL$ is the lattice (= free abelian group) of linear relations in $\gA$. It is often convenient to fix an numbering for the elements
of $\gA$, i.e. $\gA\,=\,\{\va_1,\ldots,\va_N\}$. Then $\LL$ can be described as
\begin{equation}\label{eq:L}
\LL:=\{(\ell_1,\ldots,\ell_N)\in\ZZ^N\;|\;
\ell_1\va_1+\ldots+\ell_N\va_N\,=\,\mf{0}\}\,.
\end{equation}
The rank of $\LL$ and the dimension of the torus in 
(\ref{eq:quotient torus}) are $d:=N-k-1$.
In order to obtain the natural space on which the GKZ hypergeometric
structures live one must compactify the complex torus in (\ref{eq:quotient torus}). For this purpose Gelfand, Kapranov and Zelevinsky 
developed the theory of the \emph{Secondary Fan}.
This is a complete fan of rational polyhedral cones in the real vector 
space $\LL_\RR^\vee:=\Hom(\LL,\RR)$.
Sections \ref{section:secondary fan} and \ref{section:gamma fan}
give full details about the Secondary Fan and the associated toric variety $\cV_\gA$.
Since the Secondary Fan has interesting applications outside the theory
of hypergeometric systems Sections \ref{section:secondary fan} and \ref{section:gamma fan} are written so that they can be read independently of other sections.
The toric variety $\cV_\gA$ provides a very clear picture of the `geography' for the
domains of convergence of the various GKZ hypergeometric series,
since these match exactly with discs about the special points of $\cV_\gA$ coming from the maximal cones in the secondary fan (see Proposition
\ref{prop:convergence+fan}). 
For the examples in Figure \ref{fig:classical examples} the toric varieties and special points associated with the maximal cones in the secondary fan are:
for Gauss the projective line $\PP^1$ with points
$[1,0],\; [0,1]$, for $F_4$ the projective plane $\PP^2$ with points
$[1,0,0],\; [0,1,0],\; [0,0,1]$, for $F_1$ the projective plane blown 
up in the three points $[1,0,0],\; [0,1,0],\; [0,0,1]$ equiped with
the six points of intersection of the exceptional divisors and the proper transforms of the coordinate axes in $\PP^2$.

For most $\gA$ the dimension of local solution spaces for the 
GKZ differential equations equals the volume of the 
$k$-dimensional polytope
$\Delta_\gA\,:=\,\mathrm{convex\: hull\: of\: }\gA$ (see Section
\ref{subsection:dimension solution space}); here the volume is 
normalized  as $k!\times$ the Euclidean volume. Thus for the examples
in Figure \ref{fig:classical examples} the local solution spaces
have dimension $2$, $3$, $4$, respectively.
For generic $\gA$ and $\vc$ the $\gG$-series provide bases of local solutions for the GKZ differential equations. However, in some exceptional,
but very important, cases there are not enough $\gG$-series,
due to a phenomenon called \emph{resonance}.
In Section \ref{section:resonance}
we discuss resonance and demonstrate how one \emph{sometimes} can 
obtain enough solutions by considering infinitesimal deformations of $\gG$-series. \emph{Sometimes} here means under the severe restrictions that $\vc=0$ and that one works in the neighborhood of a point on $\cV_\gA$ which corresponds to a unimodular triangulation of the polytope
$\Delta_\gA$. Very recently Borisov and Horja \cite{BH} found a way to obtain enough solutions for any $\vc\in\ZZ^{k+1}$ and any triangulation.
Their method is close in spirit to the method in Section \ref{section:resonance} and \cite{BH} gives an up-to-date presentation
of this aspect of GKZ hypergeometric structures.

In the 1980's, while Gelfand, Kapranov and Zelevinsky were working
on new hypergeometric structures, physicists discovered
fascinating new structures in string theory: the so-called
\emph{string dualities}. One of these string dualities,
known as \emph{Mirror Symmetry}, soon attracted the attention of
mathematicians, because it claimed very striking consequences for
enumerative geometry. Especially the paper \cite{COGP} of
Candelas, de la Ossa, Green and Parkes with a detailed study
of the quintic in $\PP^4$ played a pivotal role. Batyrev \cite{B1} pointed out that many examples of the Mirror Symmetry phenomenon
dealt with pairs of families of Calabi-Yau hypersurfaces in toric varieties coming from dual polytopes.
In \cite{BB} Batyrev and Borisov extended this kind of Mirror Symmetry to
Calabi-Yau complete intersections in toric varieties.
Batyrev (\cite{B2} thm.14.2) also noticed that the solutions to the 
differential equations which appeared in Mirror Symmetry, were  
solutions to GKZ hypergeometric systems
constructed from the same data as the toric varieties.
The converse is, however, not true: the GKZ system can have solutions
which are not solutions to the system of
differential equations in Mirror Symmetry.
This means that the latter system contains extra differential equations in addition to those of the underlying GKZ system 
(see \cite{HKTY} \S 3.3).
On the other hand, the solutions to the differential equations which 
one encounters in Mirror Symmetry, can all be obtained by a few differentiations from solutions to 
extremely resonant GKZ hypergeometric systems with $\vc=\mf{0}$.
Thus we do not need those extended GKZ systems.
In Section \ref{section:mirror symmetry} we discuss some examples
of this intriguing application of GKZ hypergeometric structures
to String Theory.

The quotient of two solutions of a
GKZ system of differential equations associated with $\gA$ and $\vc$
is $\TT^{k+1}$-invariant. So one can
define (at least locally) a \emph{Schwarz map} 
from the toric variety $\cV_\gA$
to the projectivization of the vector space of (local) solutions.
For Gauss's system, and more generally for Lauricella's $F_D$,
the toric variety $\cV_\gA$ and the projectivized solution space have
the same dimension, equal to the rank of $\LL$. The Schwarz map
for Gauss's system and Lauricella's $F_D$ is discussed extensively in other lectures in this school, e.g. \cite{loo}. Quite in contrast with $F_D$ is the situation for GKZ systems associated with families of Calabi-Yau threefolds. For these the toric variety $\cV_\gA$ has
dimension equal to $\rank \LL$, but the projectivized local
solution space has dimension $1+2\,\rank \LL$. The discussion about the 
\emph{canonical coordinates} and the \emph{pre-potential} in Section
\ref{subsection:Schwarz} can be seen as a description of the image of
the (local) Schwarz map. This is closely related to what in
the (physics) literature is called \emph{Special K\"ahler Geometry}.

All this basically concerns only local aspects of GKZ systems of differential equations.
About singularities, global solutions or global monodromy
of the system not much seems to be known, except for classically studied systems like Gauss's and Lauricella's $F_D$.

Since these notes are intended as an introduction to GKZ hypergeometric structures, we have included throughout the text many
examples and a few exercises. On the other hand we had to omit many 
topics. One of these omissions concerns \emph{$\gA$-discriminants}.
These come up when one identifies $\CC^\gA$ with the space of 
Laurent polynomials in $k+1$ variables with exponents in $\gA$,
$$
\vu\in\CC^\gA\qquad\leftrightarrow\qquad 
\sum_{\va\in\gA}\vu_\va\vx^\va\,, 
$$
and then wonders about the Laurent polynomials with singularities,
i.e. for which there is a point at which all partial derivatives vanish.
For the $\gA$-discriminant and its relation to the secondary fan we refer 
to \cite{gkz3}.
Another omission concerns \emph{symplectic geometry} 
in connection with the secondary fan. We recommend Guillemin's book 
\cite{gui} for further reading on this topic.

\

\noindent
\textbf{Acknowledgments.} This text is an expanded version of the notes for my lectures in the summer school on 
\textit{Arithmetic and Geometry around Hypergeometric Functions}
at Galatasaray University in Istanbul, june 2005.
I want to thank the organizers, in particular Prof. Uludag,
for the hospitality and
for the opportunity to lecture in this summer school.
I am also much indebted to the authors of the many papers 
on hypergeometric structures, triangulations and mirror symmetry
from which I myself learnt this subject.

\tableofcontents

\section{GKZ systems via examples.}
\label{section GKZ DE}
\subsection{Roots of polynomial equations.}
It is clear that in general the zeros of a polynomial
\begin{equation}\label{eq:degree n}
P_\bu(x):=u_0+u_1x+u_2x^2+\ldots+u_nx^n
\end{equation}
are functions of the coefficients $\bu=(u_0,\ldots,u_n)$.
One wonders: \textit{What kind of functions?}
For instance,
it has been known since ancient times that the zeros of a quadratic polynomial
$
ax^2+bx+c
$
are $\frac{1}{2a}(-b\pm\sqrt{b^2-4ac})$.
Similar formulas exist for polynomials of degrees $3$ and $4$, but, according to Galois theory, the zeros of a
general polynomial of degree $\geq 5$
can not be obtained from the polynomial's coefficients by a finite number of algebraic operations.
Changing the point of view
K. Mayr proved that the roots of polynomials are solutions of certain systems of differential equations:

\begin{theorem}\label{thm:Mayr}\textup{(Mayr \cite{mayr})}
If all roots of the equation $P_\bu(\xi)=0$ 
are simple, then a root $\xi$ satisfies the differential equations:
for $i_1+\ldots+i_r=j_1+\ldots+j_r$:
$$
\frac{\partial^r\xi}{\partial u_{i_1}\ldots\partial u_{i_r}}\;=\;
\frac{\partial^r\xi}{\partial u_{j_1}\ldots\partial u_{j_r}}\,.
$$
\end{theorem}
\begin{proof}
By differentiating the equation $P_\bu(\xi)=0$ with respect to $u_i$
we find
$P_\bu'(\xi)\frac{\partial\xi}{\partial u_i}+\xi^i=0$.
This implies
$\frac{\partial\xi}{\partial u_i}=\xi^i\frac{\partial\xi}{\partial u_0}
=\frac{1}{1+i}
\frac{\partial\xi^{1+i}}{\partial u_0}$. 
Induction now gives
$$
\frac{\partial^r\xi}{\partial u_{i_1}\ldots\partial u_{i_r}}\;=\;
\frac{1}{1+i_1+\ldots+i_r}
\frac{\partial^r\xi^{1+i_1+\ldots+i_r}}{\partial u_0^r}\,.
$$
\end{proof}

It obviously suffices to use only those Mayr's differential equations for which $\{i_1,\ldots,i_r\}\cap\{j_1,\ldots,j_r\}=\emptyset$. These
can also be written as

\begin{equation}\label{eq:Mayr DE3}
\prod_{\ell_i<0}\left(\frac{\partial}{\partial u_i}\right)^{-\ell_i}
\xi
\;=\;
\prod_{\ell_i>0}\left(\frac{\partial}{\partial u_i}\right)^{\ell_i}
\xi\qquad\textrm{if}\qquad
\sum_{i=0}^n \ell_i\VE{1\\ i}\,=\,
\left[\begin{array}{c} 0\\ 0\end{array}\right]\,.
\end{equation}

A second system of differential equations, satisfied by the roots of polynomials, follows from the easily checked fact
that for all $s\in\CC^*$:
$$
\xi(su_0,su_1,\ldots,su_n)\,=\,\xi(u_0,\ldots,u_n)\,,\quad
\xi(u_0,su_1,\ldots,s^nu_n)\,=\,s^{-1}\xi(u_0,\ldots,u_n)\,.
$$
When we differentiate this with respect to $s$ and set $s=1$, we find:
$$
\begin{array}{rcl}
\displaystyle{
u_0\frac{\partial\xi}{\partial u_0}+u_1\frac{\partial\xi}{\partial u_1}+
u_2\frac{\partial\xi}{\partial u_2}+
\ldots+u_n\frac{\partial\xi}{\partial u_n}}&=&0\,,
\\[1.5ex]
\displaystyle{
0\,u_0\frac{\partial\xi}{\partial u_0}+
1\,u_1\frac{\partial\xi}{\partial u_1}+
2\,u_2\frac{\partial\xi}{\partial u_2}+
\ldots+n\,u_n\frac{\partial\xi}{\partial u_n}}&=&-\xi\,,
\end{array}
$$
This can be written more transparently as:
\begin{eqnarray}
\label{eq:torus act1}
\xi(tu_0,tsu_1,ts^2u_2\ldots,ts^nu_n)&=&s^{-1}\xi(u_0,\ldots,u_n)
\quad\textrm{for}\quad (t,s)\in (\CC^*)^2\,,
\\[1ex]
\label{eq:torus DE}
\sum_{i=0}^n \VE{1\\ i}
u_i\frac{\partial\xi}{\partial u_i}
&=&
\VE{0\\ -1}\xi\,.
\end{eqnarray}
For more on zeros of 1-variable polynomials and hypergeometric functions see \cite{PT}.

\subsection{Integral with polynomial integrand}
\label{subsection:integrand 1}

Consider the integral
$$
I^{(m)}_\gs\,=\,I^{(m)}_\gs(u_0,\ldots,u_n)
\,:=\,\int_\gs P_\bu(x)^m\,\frac{dx}{x}
$$
with $m\in\ZZ$, $P_\bu(x)$
as in (\ref{eq:degree n}) and 
$\gs$ a circle in $\CC$, with radius $>0$, centred at $0$, independent of $u_0,\ldots,u_n$, not passing through any zero of $P_\bu(x)$.

By differentiating under the integral sign we see
$$
\frac{\partial I^{(m)}_\gs}{\partial u_{i}}\,=\,
m\int_\gs x^i P_\bu(x)^{m-1}\,\frac{dx}{x}
$$
and hence, if $i_1+\ldots+i_r=j_1+\ldots+j_r$, then
$$
\frac{\partial^rI^{(m)}_\gs}{\partial u_{i_1}\ldots\partial u_{i_r}}\;=\;
\frac{\partial^rI^{(m)}_\gs}{\partial u_{j_1}\ldots\partial u_{j_r}}
\,.
$$
As before this can also be written as
\begin{equation}\label{eq:Mayr DE3b}
\prod_{\ell_i<0}\left(\frac{\partial}{\partial u_i}\right)^{-\ell_i}
I^{(m)}_\gs
\;=\;
\prod_{\ell_i>0}\left(\frac{\partial}{\partial u_i}\right)^{\ell_i}
I^{(m)}_\gs\quad\textrm{if}\quad
\sum_{i=0}^n \ell_i\VE{1\\ i}\,=\,
\VE{0\\ 0}\,.
\end{equation}
For $s\in\CC^*$ close to $1$ one checks:
$I^{(m)}_\gs(su_0,\ldots,su_n)=s^mI^{(m)}_\gs(u_0,\ldots,u_n)$
and
$$
I^{(m)}_\gs(u_0,su_1,\ldots,s^nu_n)=
\int_\gs P_\bu(sx)^m\,\frac{dx}{x}=
\int_{s\gs}P_\bu(x)^m\,\frac{dx}{x}=I^{(m)}_\gs(u_0,\ldots,u_n)\,.
$$
More transparently: for $(t,s)\in (\CC^*)^2$ sufficiently close to $(1,1)$
\begin{equation}\label{eq:torus act2}
I^{(m)}_\gs(tu_0,tsu_1,ts^2u_2\ldots,ts^nu_n)\,=\,
t^mI^{(m)}_\gs(u_0,\ldots,u_n)\,.
\end{equation}
By differentating (\ref{eq:torus act2}) with respect to $t$ and $s$ and setting $t=s=1$ we find,
similar to (\ref{eq:torus DE}),
\begin{equation}\label{eq:torus DE2}
\sum_{i=0}^n \VE{1\\ i}
u_i\frac{\partial I^{(m)}_\gs}{\partial u_i}
\,=\,\VE{m\\ 0}I^{(m)}_\gs\,.
\end{equation}
Note the fundamental role of the set 
$
\gA=\left\{\VE{1\\ i}\;|\;
i=0,1,\ldots,n\right\}
$
in (\ref{eq:Mayr DE3})--(\ref{eq:torus DE2}).
Notice also the torus action (\ref{eq:torus action on functions}) 
on the left hand sides of
(\ref{eq:torus act1}) and (\ref{eq:torus act2}).

\subsection{Integral with $k$-variable Laurent polynomial integrand}
\label{subsection:integrand k}

Let us take a Laurent polynomial in $k$ variables
\begin{equation}\label{eq:nvar poly}
P_\bu(x_1,x_2,\ldots,x_k):=\sum_{\mf{a}\in\sA}u_{\mf{a}}\,
x_1^{a_1}x_2^{a_2}\cdots x_k^{a_k}
\end{equation}
where $\mf{a}=(a_1,a_2,\ldots,a_k)$ and $\sA=\{\mf{a}_1,\ldots,\mf{a}_N\}$ is a finite subset of $\ZZ^k$. Consider the integral
\begin{equation}\label{eq:integral k}
I^{(m)}_\gs(\bu)
\,:=\,\int_\gs P_\bu(x_1,\ldots, x_k)^m\,
\frac{dx_1}{x_1}\cdots \frac{dx_k}{x_k}
\end{equation}
with $\bu=(u_{\mf{a}})_{\mf{a}\in\sA}$, $m\in\ZZ$ and with
$\gs=\gs_1\times\ldots\times\gs_k$ a product of $k$ circles $\gs_1,\ldots,\gs_k$ in $\CC$, centred at $0$, independent of $\bu$, so that $P_\bu(x_1,\ldots, x_k)\neq 0$
for all $(x_1,\ldots, x_k)\in\gs_1\times\ldots\times\gs_k$.

By differentiating under the integral sign we see, for 
$\mf{a}=(a_1,\ldots,a_k)$,
$$
\frac{\partial I^{(m)}_\gs(\bu)}{\partial u_{\mf{a}}}\,=\,
m\int_\gs x_1^{a_1}\cdots x_k^{a_k} P_\bu(x_1,\ldots, x_k)^{m-1}\,\frac{dx_1}{x_1}\cdots \frac{dx_k}{x_k}\,.
$$
From this one derives that for every vector 
$(\ell_1,\ldots,\ell_N)\in\ZZ^N$ which satisfies
\begin{equation}\label{eq:LA}
\ell_1+\ldots+\ell_N=0\,,\qquad
\ell_1\mf{a}_1+\ldots+\ell_N\mf{a}_N=0\,,
\end{equation}
the following differential equation holds
\begin{equation}\label{eq:GKZ12}
\prod_{\ell_i<0}\left(\frac{\partial}{\partial u_i}\right)^{-\ell_i}
I^{(m)}_\gs(\bu)
\;=\;
\prod_{\ell_i>0}\left(\frac{\partial}{\partial u_i}\right)^{\ell_i}
I^{(m)}_\gs(\bu)\,;
\end{equation}
for simplicity of notation we write here and henceforth $u_i$ instead of $u_{\mf{a}_i}$.

For $s\in\CC^*$ sufficiently close to $1$ and for $i=1,\ldots, k$ one calculates:
\begin{eqnarray*}
&&
I^{(m)}_\gs(s^{a_{i1}}u_1,s^{a_{i2}}u_2,\ldots,s^{a_{iN}}u_N)\;=\;
\int_\gs P_\bu(x_1,\ldots,sx_i,\ldots, x_k)^m\,
\omega\;=\\
&&=\;
\int_{\gs_1\times\ldots\times s\gs_i\times\ldots\times\gs_k}
\hspace{-5em}P_\bu(x_1,\ldots,x_k)^m\,
\omega\;=\;
\int_{\gs_1\times\ldots\times \gs_i\times\ldots\times\gs_k}
\hspace{-5em}P_\bu(x_1,\ldots,x_k)^m\,
\omega\;=\;
I^{(m)}_\gs(u_1,\ldots,u_N)\,;
\end{eqnarray*}
here $a_{ij}$ denotes the $i$-th coordinate of the vector $\mf{a}_j$
and $\omega=\frac{dx_1}{x_1}\cdots \frac{dx_k}{x_k}$.
This together with $I^{(m)}_\gs(su_1,\ldots,su_N)=s^mI^{(m)}_\gs(u_1,\ldots,u_N)$
can also be written as: 
\begin{equation}\label{eq:torus act3}
I^{(m)}_\gs(ts_1^{a_{11}}s_2^{a_{21}}\cdots s_k^{a_{k1}}u_1,
\ldots,ts_1^{a_{1N}}s_2^{a_{2N}}\cdots s_k^{a_{kN}}u_N)
\,=\,
t^mI^{(m)}_\gs(u_0,\ldots,u_n)
\end{equation}
for $(t,s_1,\ldots,s_k)\in (\CC^*)^{k+1}$ close to $(1,\ldots,1)$.
By differentiating with respect to $t,s_1,\ldots,s_k$ and setting 
$t=s_1=\ldots=s_k=1$ we find
\begin{equation}\label{eq:torus DEk}
\VE{1\\ \mf{a}_1}u_1\frac{\partial I^{(m)}_\gs(\bu)}{\partial u_1}+
\ldots+
\VE{1\\ \mf{a}_N}u_N\frac{\partial I^{(m)}_\gs(\bu)}{\partial u_N}
\,=\,\VE{m\\ \mf{0}}I^{(m)}_\gs(\bu)\,.
\end{equation}
Note the appearance of the set 
$
\gA\,=\,\left\{\VE{ 1\\ \mf{a}}\in
\ZZ^{k+1}\;|\;\mf{a}\in\sA\right\}
$ 
in (\ref{eq:LA}) and (\ref{eq:torus DEk}).
Notice also the torus action (\ref{eq:torus action on functions}) on the 
left hand side of (\ref{eq:torus act3}). 

\

\noindent
\textbf{Remark}
For $m>0$ one can evaluate $I^{(m)}_\gs(\bu)$ using the multinomial and residue theorems. One finds that $I^{(m)}_\gs(\bu)$ is actually a polynomial:
\begin{equation}\label{eq:polynomial solution}
\frac{1}{(2\pi\ci)^k}I^{(m)}_\gs(\bu)\,=\,
\sum_{(m_1,\ldots,m_N)}\frac{m!}{(m_1)!\cdots (m_N)!}
u_1^{m_1}\cdots u_N^{m_N}
\end{equation}
where the sum runs over all $N$-tuples of non-negative integers
$(m_1,\ldots,m_N)$ satisfying 
$m_1+\ldots+m_N=m$ and $m_1\mf{a}_1+\ldots+m_N\mf{a}_N=0$.

In Section \ref{section:mirror symmetry} one can find explicit examples
of these integrals with $m=-1$.

\subsection{Generalized Euler integrals}
\label{subsection:Euler}
In \cite{gkz2,gkz4} Gelfand, Kapranov and Zelevinsky investigate integrals
of the form
\begin{equation}\label{eq:euler}
\int_\gs \prod_i P_i(x_1,\ldots,x_k)^{\ga_i}\,x_1^{\beta_1}\cdots 
x_k^{\beta_k}\,dx_1\cdots dx_k\,,
\end{equation}
which they call \emph{generalized Euler integrals}.
Here the $P_i$ are Laurent polynomials, $\ga_i$ and $\beta_j$ are
complex numbers and $\gs$ is a $k$-cycle.
Since the integrand can be multivalued and can have singularities one must carefully give the precise
meaning of Formula (\ref{eq:euler}) (see \cite{gkz2} \S 2.2).
Having dealt with the technicalities of the precise definition
Gelfand, Kapranov and Zelevinsky view the integrals (\ref{eq:euler})
as functions of the coefficients of the Laurent polynomials $P_i$.
Using the same arguments as we used in Section 
\ref{subsection:integrand k} they then verify that these functions
satisfy a system of differential equations (\ref{eq:GKZ1})-(\ref{eq:GKZ2}) for the appropriate data $\gA$ and $\vc$. Examples can be found in
Sections \ref{subsection:IntFD} and \ref{subsection:cubics intersect}.

\subsection{General GKZ systems of differential equations}
\label{subsection:GKZ}
The systems of differential equations 
(\ref{eq:Mayr DE3})-(\ref{eq:torus DE}),
(\ref{eq:Mayr DE3b})-(\ref{eq:torus DE2}) and
(\ref{eq:GKZ12})-(\ref{eq:torus DEk})
found in the preceding examples
are special cases of systems of differential equations discovered by
Gelfand, Kapranov and Zelevinsky \cite{gkz1,gkz2,gkz4}. The general GKZ system for functions $\Phi$ of $N$ variables $u_1,\ldots,u_N$ is
constructed from a vector $\vc\in\CC^{k+1}$
and an $N$-element subset $\gA=\{\va_1,\ldots,\va_N\}\subset\ZZ^{k+1}$ which generates $\ZZ^{k+1}$ as an abelian group and for which there exists a group homomorphism 
$h:\ZZ^{k+1}\rightarrow\ZZ$ such that $h(\va)=1$ for all $\va\in\gA$.
Let $\LL\subset\ZZ^N$ denote the lattice of relations in $\gA$:
$$
\LL:=\{(\ell_1,\ldots,\ell_N)\in\ZZ^N\;|\;
\ell_1\va_1+\ldots+\ell_N\va_N\,=\,\mf{0}\}\,.
$$
Note that the condition $h(\va)=1$ for all $\va\in\gA$, implies that
$\ell_1+\ldots+\ell_N=0$ for every  
$(\ell_1,\ldots,\ell_N)\in\LL$.
\begin{definition}\label{def:GKZ system}
The GKZ system associated with $\gA$ and $\vc$
consists of 
\begin{itemize}
\item
for every $(\ell_1,\ldots,\ell_N)\in\LL$
one differential equation
\begin{equation}\label{eq:GKZ1}
\prod_{\ell_i<0}\left(\frac{\partial}{\partial u_i}\right)^{-\ell_i}\Phi
\;=\;
\prod_{\ell_i>0}\left(\frac{\partial}{\partial u_i}\right)^{\ell_i}\Phi\,,
\end{equation}
\item
the system of $k+1$ differential equations
\begin{equation}\label{eq:GKZ2}
\va_1\,u_1\frac{\partial  \Phi}{\partial u_1}+
\ldots+\va_N\,u_N\frac{\partial  \Phi}{\partial u_N}\;=\;
\vc\Phi\,.
\end{equation} 
\end{itemize}
\end{definition}

\

\noindent
\textbf{Remark.} It is natural to view $u_1,\ldots,u_N$ as coordinates
on the space $\CC^\gA:=\mathrm{Maps}(\gA,\CC)$.
Then the left hand side of the equation (\ref{eq:GKZ2}) is the 
infinitesimal version of the torus action 
(\ref{eq:torus action on functions}).
If $\Phi_1$ and $\Phi_2$ are two solutions of (\ref{eq:GKZ2}) on some open
set $U\subset\CC^\gA$, their quotient satisfies
\begin{equation}\label{eq:invariant quotient}
\va_1\,u_1\frac{\partial}{\partial u_1}
\left(\frac{\Phi_1}{\Phi_2}\right)+
\ldots+\va_N\,u_N\frac{\partial}{\partial u_N}
\left(\frac{\Phi_1}{\Phi_2}\right)\,=\,0
\end{equation}
and is therefore constant on the intersections of $U$ with the $\TT^{k+1}$-orbits. 
\\
\emph{Thus a basis $\Phi_1,\ldots,\Phi_r$ of the solution space of
(\ref{eq:GKZ1})-(\ref{eq:GKZ2}) induces map from the orbit space
$\modquot{\TT^{k+1}\cdot U}{\TT^{k+1}}$ into the projective space
$\PP^{r-1}$, like the \emph{Schwarz map} for Gauss's hypergeometric
systems.}

\

Another simple, but nevertheless quite useful, 
consequence of the GKZ differential equations is:

\begin{proposition}\label{prop:GKZ derivatives}
If function $\Phi$ satisfies the differential equations
(\ref{eq:GKZ1})-(\ref{eq:GKZ2}) for $\gA$ and $\vc$,
then $\frac{\partial\Phi}{\partial u_j}$ satisfies the differential 
equations (\ref{eq:GKZ1})-(\ref{eq:GKZ2}) for $\gA$ and $\vc-\va_j$.
\end{proposition}
\begin{proof}
The derivation $\frac{\partial}{\partial u_j}$ commutes with all
derivations involved in (\ref{eq:GKZ1}). On the other hand by applying
$\frac{\partial}{\partial u_j}$ to both sides of (\ref{eq:GKZ2}) we get
$$
\va_1\,u_1\frac{\partial}{\partial u_1}
\left(\frac{\partial\Phi}{\partial u_j}\right)+
\ldots+\va_N\,u_N\frac{\partial}{\partial u_N}
\left(\frac{\partial\Phi}{\partial u_j}\right)
+\va_j\frac{\partial\Phi}{\partial u_j}\,=\,
\vc\frac{\partial\Phi}{\partial u_j}\,.
$$
\end{proof}
\subsection{Gauss's hypergeometric differential equation as a GKZ system}
\label{section:Gauss DE}

The most classical hypergeometric differential equation, due to Euler and Gauss, is:
\begin{equation}\label{eq:gauss}
z(z-1)F''+((a+b+1)z-c)F'+abF\;=\;0\,.
\end{equation}
Here $F$ is a function of one variable $z$, ${}'=\frac{d}{dz}$ and $a,b,c$ are additional complex parameters. It is reproduced in the GKZ formalism
by $\vc=(1-c,-a,-b)$ and

\begin{picture}(230,55)(0,20)
\put(50,20){\makebox(100,50){
$
\gA=\left\{
\VE{1\\ 1\\ 1},\VE{-1\\ 0\\ 0},
\VE{0\\ 1\\ 0},\VE{0\\ 0\\ 1}
\right\}\subset\ZZ^3
$
}}
\put(230,28){
\begin{picture}(100,50)(0,0)
\put(0,0){\line(1,0){40}}
\put(0,0){\line(0,1){40}}
\put(40,40){\line(-1,0){40}}
\put(40,40){\line(0,-1){40}}
\put(0,0){\circle*{5}}
\put(40,0){\circle*{5}}
\put(0,40){\circle*{5}}
\put(40,40){\circle*{5}}
\end{picture}}
\end{picture}
\\
and, hence,
$\LL\,=\,\ZZ(1,1,-1,-1)\subset\ZZ^4$.
Indeed, for these data the GKZ system boils down to the following four differential equations for a function $\Phi$ of four variables 
$(u_1,u_2,u_3,u_4)$:
$$
\begin{array}{rcr}
\displaystyle{
\frac{\partial^2\Phi}{\partial u_1\partial u_2}}&=&\displaystyle{
\frac{\partial^2\Phi}{\partial u_3\partial u_4} } \\[2ex]
\displaystyle{
u_1\frac{\partial\Phi}{\partial u_1}\,-\,
u_2\frac{\partial\Phi}{\partial u_2}}&=&\displaystyle{ (1-c)\:\Phi
} \\[2ex]
\displaystyle{
u_1\frac{\partial\Phi}{\partial u_1}\,+\,
u_3\frac{\partial\Phi}{\partial u_3}}&=&\displaystyle{ -a\:\Phi
} \\[2ex]
\displaystyle{
u_1\frac{\partial\Phi}{\partial u_1}\,+\,
u_4\frac{\partial\Phi}{\partial u_4}}&=&\displaystyle{ -b\:\Phi
}
\end{array}
$$
From the second equation we get
$$
\frac{\partial^2\Phi}{\partial u_1\partial u_2}\;=\;
u_2^{-1}\left(
u_1\frac{\partial^2\Phi}{\partial u_1^2}+
c\frac{\partial\Phi}{\partial u_1}
\right)
$$
From the third and fourth equations we get
$$
\frac{\partial^2\Phi}{\partial u_3\partial u_4}\,=\,
u_3^{-1}u_4^{-1}
\left(-u_1\frac{\partial}{\partial u_1}-a\right)
\left(-u_1\frac{\partial}{\partial u_1}-b\right)\Phi$$
Together with the first equation this yields
$$
u_3^{-1}u_4^{-1}
\left(u_1^2\frac{\partial^2\Phi}{\partial u_1^2}+
(1+a+b)u_1\frac{\partial\Phi}{\partial u_1}+ab\Phi\right)
\,=\, 
u_2^{-1}\left(
u_1\frac{\partial^2\Phi}{\partial u_1^2}\,+\,
c\frac{\partial\Phi}{\partial u_1}
\right)\,.
$$
Setting $u_2=u_3=u_4=1$, $u_1=z$ and $F(z)=\Phi(z,1,1,1)$ we find that $F$ satisfies the differential equation (\ref{eq:gauss}).

\subsection{Dimension of the solution space of a GKZ system}
\label{subsection:dimension solution space}
The spaces of (local) solutions of the GKZ differential equations (\ref{eq:GKZ1})-(\ref{eq:GKZ2}) are complex vector spaces.
Theorems 2 and 5 in \cite{gkz1} state that the dimension of the space
of (local) solutions of (\ref{eq:GKZ1})-(\ref{eq:GKZ2}) near a generic point 
is equal to the normalized volume of the $k$-dimensional polytope
$
\Delta_\gA:=\textrm{convex hull}(\gA)\,;
$
here `normalized volume' means $k!$ times the usual Euclidean volume.
In \cite{gkz1a} it is pointed out that the proof in \cite{gkz1} requires
an additional condition on $\gA$. Corollary 8.9 and Proposition 13.15
in \cite{stu} show that this additional condition is satisfied if
the polytope $\Delta_\gA$ admits a unimodular triangulation.
Triangulations of $\Delta_\gA$ and their importance in GKZ hypergeomtric structures are discussed in Section 
\ref{subsection:triangulations}.

\section{$\gG$-series}
\label{section:gamma series}
As before we consider a subset $\gA=\{\va_1,\ldots,\va_N\}\subset\ZZ^{k+1}$ which generates $\ZZ^{k+1}$ as an abelian group and for which there exists a group homomorphism 
$h:\ZZ^{k+1}\rightarrow\ZZ$ such that $h(\va)=1$ for all $\va\in\gA$.
And, still as before, we write:
$$
\LL:=\{(\ell_1,\ldots,\ell_N)\in\ZZ^N\;|\;
\ell_1\va_1+\ldots+\ell_N\va_N\,=\,\mf{0}\}\,.
$$
The condition $h(\va)=1$ for all $\va\in\gA$, implies that
$\ell_1+\ldots+\ell_N=0$ for every  
$(\ell_1,\ldots,\ell_N)\in\LL$.
With $\LL$ and a vector
$\ugamma=( \gamma_1,\ldots, \gamma_N)\in\CC^N$ Gelfand, Kapranov and Zelevinsky \cite{gkz1}
associate what they call a \emph{$\gG$-series}:
\begin{definition}\label{def:gamma series}
The $\gG$-series associated with $\LL$ and  
$\ugamma=( \gamma_1,\ldots, \gamma_N)\in\CC^N$ is
\begin{equation}\label{eq:Gamma series}
\Phi_{\LL,\ugamma}(u_1,\ldots,u_N)\;=\;
\sum_{(\ell_1,\ldots,\ell_N)\in\LL}\prod_{j=1}^N
\frac{u_j^{\gamma_j+\ell_j}}{\Gamma(\gamma_j+\ell_j+1)}\,.
\end{equation}
\end{definition}
\noindent
Here $\gG$ is the $\gG$-function; its definition and main properties are recalled in Section \ref{subsection:gamma function}.
In Section \ref{subsection:gamma examples} we demonstrate how the classical
hypergeometric series of Gauss, Appell and Lauricella appear in the $\gG$-series format.
In Section \ref{subsection:gamma growth} we give estimates for the growth of the coefficients in (\ref{eq:Gamma series}).
Formula (\ref{eq:Gamma series}) requires for $\ugamma\not\in\ZZ^{k+1}$
choices of logarithms for $u_1,\ldots,u_N$. By carefully manoeuvreing
conditions on $\ugamma$ and substitutions setting some $u_j$ equal to $1$,
we can avoid problems and show in Section \ref{subsection:gamma DE}
how a $\gG$-series can be viewed as a power series in $d=N-k-1$ variables with positive radii of convergence.
Nevertheless, a formula avoiding choices of logarithms is desirable.
For that reason we introduce Fourier $\gG$-series in Section
\ref{subsection:gamma fourier}. 
In Section \ref{subsection:gamma DE} we prove that 
$\Phi_{\LL,\ugamma}(u_1,\ldots,u_N)$ can be viewed as a function on
some domain in $(u_1,\ldots,u_N)$-space and satisfies the GKZ differential equations.

\subsection{The $\gG$-function}
\label{subsection:gamma function}
The $\gG$-function is defined for complex numbers $s$ with $\Re s>0$
by the integral
\begin{equation}\label{eq:Gamma def}
\gG (s):=\int_0^\infty t^{s-1}e^{-t}dt\,.
\end{equation}
Using partial integration one immediately checks
$
\gG (s+1)\,=\, s\gG (s) 
$
and, hence, for $n\in\ZZ,\,n>0$
\begin{equation}\label{eq:Gamma ext}
\gG (s+n)\,=\, s(s+1)\ldots(s+n-1)\gG (s)\,.
\end{equation}
Formulas (\ref{eq:Gamma def}) and (\ref{eq:Gamma ext}) imply in
particular
\begin{equation}
\label{eq:factorial}
\gG (1)\,=\,1\,,\qquad
\gG (n+1)\,=\,n! \qquad\textrm{for } n\in\NN\,.
\end{equation}
One can extend the $\gG$-function to a meromorphic function on all of $\CC$ 
by setting 
\begin{equation}\label{eq:Gamma def2}
\gG (s)\,=\,\frac{\gG (s+n)}{s(s+1)\ldots(s+n-1)}\qquad \textrm{with}\quad
n\in\ZZ,\,n>-\Re s\,.
\end{equation} 
The functional equation (\ref{eq:Gamma ext}) shows that this does not depend on the choice of $n$. Formula (\ref{eq:Gamma def})
shows $\gG (s)\,\neq\,0$ if $\Re s>1$ and hence Formula 
(\ref{eq:Gamma def2})
shows that the extended $\gG$-function is holomorphic on 
$\CC\setminus\ZZ_{\leq 0}$ and has at $s=-m\in\ZZ_{\leq 0}$ a first order pole with residue
\begin{equation}\label{eq:Gamma residue}
\mathrm{Res}_{s=-m}\gG (s)\,=\,\frac{(-1)^m}{m!}\,.
\end{equation}
The function $\frac{1}{\gG(s)}$ is holomorphic on the whole complex plane.
Its zero set is $\ZZ_{\leq 0}$ and its Taylor series at 
$-m\in\ZZ_{\leq 0}$ starts like
\begin{equation}\label{eq:Gamma Taylor}
\frac{1}{\gG(s-m)}\,=\,(-1)^mm!\,s\;+\;\ldots
\end{equation}
The coefficients of (classical) hypergeometric series are usually expressed in terms of \emph{Pochhammer symbols} $(s)_n$. These
are defined by $(s)_n\,=\,s(s+1)\cdots(s+n-1)$ and can be
rewritten as quotients of $\gG$-values:
\begin{equation}\label{eq:Pochhammer}
(s)_n\,=\,s(s+1)\cdots(s+n-1)\,=\,
\frac{\gG(s+n)}{\Gamma(s)}\,=\,
(-1)^n\frac{\gG(1-s)}{\Gamma(1-n-s)}\,.
\end{equation}
Note, however, that for integer values of $s$ the Pochhammer symbol 
$(s)_n$ is perfectly well defined, while some of the individual
$\gG$-values in (\ref{eq:Pochhammer}) may become $\infty$.

\subsection{Examples of $\gG$-series}
\label{subsection:gamma examples}
\subsubsection{Gauss's hypergeometric series}
As in the example of Gauss's hypergeometric differential equation 
(Section \ref{section:Gauss DE}) we
take $\LL\,=\,\ZZ(1,1,-1,-1)$ in $\ZZ^4$ and $\ugamma=(0,c-1,-a,-b)\in\CC^4$. 
If $c$ is not an integer $\leq 0$, then, by (\ref{eq:Gamma series}) 
and (\ref{eq:Pochhammer}),
\begin{eqnarray*}
\Phi_{\LL,\ugamma}(u_1,u_2,u_3,u_4)&=&
\sum_{n\in\ZZ}
\frac{u_1^{n}u_2^{c-1+n}u_3^{-a-n}u_4^{-b-n}}{\Gamma(1+n)\Gamma(c+n)
\Gamma(1-n-a)\Gamma(1-n-b)}\\
&=&
\frac{u_2^{c-1}u_3^{-a}u_4^{-b}}{\Gamma(c)
\Gamma(1-a)\Gamma(1-b)}\sum_{n\geq 0}
\frac{(a)_n(b)_n}{n!(c)_n}(u_1u_2u_3^{-1}u_4^{-1})^n
\end{eqnarray*}
and, hence,
$
\Phi_{\LL,\ugamma}(z,1,1,1)=\frac{1}{\Gamma(c)\Gamma(1-a)\Gamma(1-b)}
F(a,b,c|z)
$
with
$$
F(a,b,c|z):=\sum_{n\geq 0}\frac{(a)_n(b)_n}{n!(c)_n}z^n\,.
$$
The power series $F(a,b,c|z)$ is Gauss's hypergeometric series.
Note that if $a$ or $b$ is a positive integer the $\gG$-series is $0$,
but Gauss's hypergeometric series is not $0$.

\subsubsection{The hypergeometric series ${}_pF_{p-1}$}
\label{subsubsection:pFq}
Quite old generalizations of Gauss's hypergeometric series
are the series
$$
{}_pF_{p-1}\left(\left.
\begin{array}{lll}a_1&,\ldots,&a_p\\ c_1&,\ldots,&c_{p-1}\end{array}
\right|z\right)\,:=\,
\sum_{n\geq 0}
\frac{(a_1)_n\cdots (a_p)_n}{n!(c_1)_n\cdots (c_{p-1})_n}
z^n\,.
$$
Like for Gauss's series one easily finds that the series ${}_pF_{p-1}$
match (up to a constant factor) the $\gG$-series for $\LL=\ZZ(1,\ldots,1,-1,\ldots,-1)$ with $p$ $\;1$'s and $p$ $\;(-1)$'s.

\subsubsection{The case ${}_1F_0$.}
The simplest, yet not totally trivial, case of a $\gG$-series arises for
$\LL=\ZZ(1,-1)\subset\ZZ^2$.
The $\gG$-series with $\ugamma=(0,a),\;a\in\CC$ is
$$
\Phi_{\ZZ(1,-1),(0,a)}(u_1,u_2)\;=\;
\sum_{n\in\ZZ}
\frac{u_1^n u_2^{a-n}}{\Gamma(1+n)
\Gamma(1+a-n)}
\;=\;\frac{1}{\gG(1+a)}(u_1+u_2)^a\,;
$$
here we use the generalized binomial theorem and (\ref{eq:Pochhammer}):
$$
\left(\begin{array}{c}a\\ n\end{array}\right)
\,=\,
\frac{a(a-1)\ldots(a-n+1)}{n!}\,=\,
\frac{\gG(1+a)}{\gG(1+n)\gG(1+a-n)}\,.
$$

\

\noindent
\textbf{Remark.}
Note that $\LL=\ZZ(1,-1)$ implies that  
the two elements of $\gA$ are equal.
The GKZ differential equations in this case imply
that the hypergeometric functions are in fact just functions of
the single variable $u_1+u_2$. This illustrates a general fact:
when setting up the theory of GKZ hypergeometric systems one could take
for $\gA$ a list of vectors in $\ZZ^{k+1}$ instead of just a subset (i.e. the elements may occur more than once). But any such apparently more general set up, arises from a case with $\gA$ a genuine set by simply
replacing a variable by a sum of new variables. So by allowing for $\gA$ a list instead of a set one does not get a seriously more general theory.
Therefore we ignore this option in these notes.

\subsubsection{The Appell-Lauricella hypergeometric series}
\label{subsubsection:Appell-Lauricella}
These are generalizations of Gauss's series to $n$ variables  defined by Appell for $n=2$ and Lauricella for general $n$. With the notations 
$\vz^{\vm}\,:=\,z_1^{m_1}\cdots z_n^{m_n}$,
$(\vx)_{\vm}\,:=\,(x_1)_{m_1}\cdots (x_n)_{m_n}$,
$\vm !\,:=\, m_1!\cdots m_n!$,
$|\vm |\,:=\,m_1+\ldots+m_n$
for $n$-tuples of complex numbers
$\vz\,=\,(z_1,\ldots,z_n)$, $\vx\,=\,(x_1,\ldots,x_n)$ and 
of non-negative integers $\vm\,=\,(m_1,\ldots,m_n)$,
the four Lauricella series are
\begin{eqnarray*}
F_A(a,\vb,\vc |\vz)&:=&\sum_{\vm}
\frac{(a)_{|\vm |}(\vb)_{\vm}}{(\vc)_{\vm} \vm !}\,\vz^{\vm}\\
F_B(\va,\vb,c |\vz)&:=&\sum_{\vm}
\frac{(\va)_{\vm}(\vb)_{\vm}}{(c)_{|\vm |} \vm !}\,\vz^{\vm}\\
F_C(a,b,\vc |\vz)&:=&\sum_{\vm}
\frac{(a)_{|\vm |}(b)_{|\vm |}}{(\vc)_{\vm } \vm !}\,\vz^{\vm}\\
F_D(a,\vb,c |\vz)&:=&\sum_{\vm}
\frac{(a)_{|\vm |}(\vb)_{\vm}}{(c)_{|\vm |} \vm !}\,\vz^{\vm}
\end{eqnarray*}
in the summations $\vm$ runs over $\ZZ^n_{\geq 0}$ and the $c$-parameters are not integers $\leq 0$.
In Appell's notation (for $n=2$) these series are called $F_2,F_3,F_4,F_1$ respectively.

One can use (\ref{eq:Pochhammer}) to explore the relations between
the Lauricella series and $\Gamma$-series. 
For the coefficients in $F_D$, for instance, we find
$$
\frac{(a)_{|\vm |}(\vb)_{\vm}}{(c)_{|\vm |} \vm !}\,=\,
\gG(1-a)\Gamma(c)\prod_{j=1}^n\gG(1-b_j)\cdot
\prod_{j=1}^N\frac{1}{\gG(1+\gam_j+\ell_j)}
$$
with $N=2n+2$,
$\quad\ugamma\,=\,
(\gam_1,\ldots,\gam_N)\,=\,(c-1,-b_1,\ldots,-b_n,-a,0,\ldots,0)\,,$
$$
(\ell_1,\ldots,\ell_N)\,=\,(m_1,\ldots,m_n)
\left(
\begin{array}{rrrrrrrrrr}
1&-1&0&\ldots&0&-1&1&0&\ldots&0\\
1&0&-1&\ddots&\vdots&-1&0&1&\ddots&\vdots\\
\vdots&\vdots&\ddots&\ddots&0&\vdots&\vdots&\ddots&
\ddots&0\\
1&0&\ldots&0&-1&-1&0&\ldots&0&1
\end{array}
\right)
$$
So for $\LL$ we take the lattice which is spanned by the rows of the
above $n\times N$-matrix.
Substuting $u_j=1$ for $1\leq j\leq n+2$ and $u_j=z_{j-n-2}$ for 
$n+3\leq j\leq 2n+2$ turns the $\gG$-series into a power series:
$$
\Phi_{\LL,\ugamma}(1,\ldots,1,z_1,\ldots,z_n)=
\left(\prod_{j=1}^{n+2}\gG(1+\gam_j)^{-1}\right)
F_D(a,\vb,c |\vz)\,.
$$

\

\noindent
\textbf{Exercise}
Note that the matrix describing $\LL$ for Lauricella's $F_D$
is 
$ (\mf{1}_n,-\II_n,-\mf{1}_n,\II_n)$ 
where $\mf{1}_n$ is the column vector with $n$ components $1$
and $\II_n$ is the $n\times n$-identity matrix.
Now find the lattice $\LL$ for the Lauricella functions $F_A$, $F_B$ and $F_C$.

\subsection{Growth of coefficients of $\gG$-series}
\label{subsection:gamma growth}

Here is first a simple lemma about the growth behavior of the $\gG$-function.

\begin{lemma} For every $C\in\CC\setminus\ZZ$ there are real constants $P,\,R,\kappa_1,\kappa_2>0$ (depending on $C$) such that for all $M\in \ZZ_{\geq 0}\::$
\begin{equation}\label{eq:schat}
|\Gamma(C+M)|\:\geq\:\kappa_1 R^M M^M\qquad\textrm{and}\qquad
|\Gamma(C-M)|\:\geq\:\kappa_2 P^{-M}M^{-M}\,.
\end{equation}
\end{lemma}
\begin{proof}
From (\ref{eq:Gamma ext}) and the triangle inequalities one derives
\begin{eqnarray*}
|\Gamma(C+M)|&\geq& \prod_{j=0}^{M-1}||C|-j|\cdot|\Gamma(C)|
\:\geq\:Q^M\,M!\,|\Gamma(C)|\:\geq\:\kappa R^M\,M^M\,|\Gamma(C)|\,,
\\
|\Gamma(C-M)|&\geq& \prod_{j=1}^{M}(|C|+j)^{-1}\cdot|\Gamma(C)|
\:\geq\:\frac{|\Gamma(C)|}{(|C|+M)^M}
\:\geq\:P^{-M}M^{-M}|\Gamma(C)|
\end{eqnarray*}
with 
$Q\,:=\,\min_{k\in\NN}\:\frac{|1+|C|-k|}{k}$, $R:=Qe^{-1}$,
$P\,:=\,2\min(1,|C|)$ and some constant $\kappa$ (from Stirling's formula).
\end{proof}

\

Now consider the coefficient
$\prod_{j=1}^N\Gamma(\gamma_j+\ell_j+1)^{-1}$ in the $\gG$-series
(\ref{eq:Gamma series}). Set
$\gamma'_j=\gamma_j$ if $\gamma_j\not\in\ZZ$ and
$\gamma'_j=\gamma_j-\frac{1}{2}$ if $\gamma_j\in\ZZ$.
Note that
$\gG(k-\frac{1}{2})=(k-\frac{3}{2})\cdots\frac{1}{2}
\gG(\frac{1}{2})\,\leq\,(k-1)!\gG(\frac{1}{2})\,=\,\gG(k)\sqrt{\pi}$
for $k\in\ZZ_{\geq 1}$.
Then, using the above lemma, one sees that there are real constants $K,S>0$ such that
$$
\left|\prod_{j=1}^N\frac{1}{\Gamma(\gamma_j+\ell_j+1)}\right|\leq
\left|\prod_{j=1}^N\frac{\sqrt{\pi}}{\Gamma(\gamma'_j+\ell_j+1)}\right|
\leq KS^D\prod_{j=1}^N|\ell_j|^{-\ell_j}
$$
with
$D\;:=\;\frac{1}{2}\sum_{j=1}^N|\ell_j|\;=\;\sum_{\ell_j>0} \ell_j\;=\;-\sum_{\ell_j<0} \ell_j$.
Since
$\prod_{\ell_j<0}|\ell_j|^{-\ell_j}\,\leq\,D^D$
and
$\prod_{\ell_j>0}|\ell_j|^{-\ell_j}\,\leq\,N^DD^{-D}$, our final estimate becomes:

\begin{proposition}\label{prop:growth} 
There are real numbers $K,T>0$, depending on
$\ugamma=(\gam_1,\ldots,\gam_N)$, but independent of
$\ull=(\ell_1,\ldots,\ell_N)$, such that
\begin{equation}\label{eq:estimate}
\left|\prod_{j=1}^N\frac{1}{\Gamma(\gamma_j+\ell_j+1)}\right|\leq
K\,T^{\sum_{j=1}^N|\ell_j|}\,.
\end{equation}
\qed
\end{proposition}

\

\subsection{$\gG$-series and power series}
\label{subsection:gamma power}

Let $J\subset\{1,\ldots,N\}$ be a set with $k+1$ elements, such that the 
vectors $\va_j$ with $j\in J$ are linearly independent. 
Write $J':=\{1,\ldots,N\}\setminus J$.
Let $\ugamma=(\gam_1,\ldots,\gam_N)\in\CC^N$ be such that $\gam_j\in\ZZ$ for $j\in J'$.
Since $\frac{1}{\gG(s)}=0$ if $s\in\ZZ_{\leq 0}$, the $\gG$-series 
(\ref{eq:Gamma series}) constructed with such a $\ugamma$ involves only terms from the set
\begin{equation}\label{eq:L cone} 
\LL_{J,\ugamma}:=\{(\ell_1,\ldots,\ell_N)\in\LL\;|\;
\gam_j+\ell_j\geq 0\quad\textrm{if}\quad j\in J'\}\,.
\end{equation}
The substitution
\begin{equation}\label{eq:power substitution}
u_j=z_j\quad\textrm{if}\quad j\in J'\,,\qquad
u_j=1\quad\textrm{if}\quad j\in J
\end{equation}
therefore turns the $\gG$-series into the power series 
\begin{equation}\label{eq:gamma power}
\sum_{(\ell_1,\ldots,\ell_N)\in\LL_{J,\ugamma}}
\left(\prod_{j=1}^N\frac{1}{\Gamma(\gamma_j+\ell_j+1)}\right)
\prod_{j\in J'} z_j^{\gamma_j+\ell_j}\,.
\end{equation}
The following lemma is needed to convert (\ref{eq:estimate}) into
estimates for the radii of convergence of this power series.

\begin{lemma}\label{lemma:LJ}
Let $J\subset\{1,\ldots,N\}$ be a set with $k+1$ elements, 
such that the 
vectors $\va_j$ with $j\in J$ are linearly independent. 
Write $J':=\{1,\ldots,N\}\setminus J$.
Then there is a positive real constant $\beta$
such that for every $(\ell_1,\ldots,\ell_N)\in\LL$
\begin{equation}\label{eq:LJ}
|\ell_1|+\ldots+|\ell_N|\leq\beta\sum_{j\in J'}|\ell_j|\,.
\end{equation}
\end{lemma}
\begin{proof} 
Take any $d\times N$-matrix $\mf{B}$ whose rows form a $\ZZ$-basis of $\LL$. Let $\bb_1,\ldots,\bb_N$ be its columns. 
Let $\mf{B}_{J'}$ denote the $d\times d$-matrix with columns 
$\bb_j\;(j\in J')$. 
Then the matrix $\mf{B}_{J'}$ is invertible over $\QQ$;
indeed, if it were not, its rows would be linearly dependent and there would be a vector
$(\ell_1,\ldots,\ell_N)\in\LL$ such that
$\ell_j=0$ for $j\in J'$; the relation
$\ell_1\va_1+\ldots+\ell_N\va_N\,=\,\mf{0}$ would contradict the  
linear independence of the vectors $\va_j$ with $j\in J$. 
Now we have the equality of row vectors
for every $(\ell_1,\ldots,\ell_N)\in\LL$
$$
(\ell_1,\ldots,\ell_N)= (\ell)_{J'} (\mf{B}_{J'})^{-1}\mf{B}
$$
where $(\ell)_{J'}$ is the row vector with components 
$\ell_j\;(j\in J')$. So for $\beta$ in (\ref{eq:LJ}) one can take the maximum of the absolute values of the entries of the matrix $(\mf{B}_{J'})^{-1}\mf{B}$.
\end{proof}

\begin{proposition}\label{prop:power convergence}
Let $J\subset\{1,\ldots,N\}$ be a set with $k+1$ elements, such that the 
vectors $\va_j$ with $j\in J$ are linearly independent.
Let $\ugamma=(\gam_1,\ldots,\gam_N)\in\CC^N$ be such that $\gam_j\in\ZZ$ for $j\in J':=\{1,\ldots,N\}\setminus J$.
Then there is an $R\in\RR_{>0}$ such that the power series 
(\ref{eq:gamma power}) converges on the polydisc given by
$|z_j|<R$ for $j=1,\ldots,d$.
\end{proposition}
\begin{proof}
This follows, with $R=T^{-\beta}$, from Proposition \ref{prop:growth} and
Lemma \ref{lemma:LJ}.
\end{proof}

\subsection{Fourier $\gG$-series}
\label{subsection:gamma fourier}
The substitutions in (\ref{eq:power substitution}) depend too rigidly on
the choice of the set $J$ and make it difficult to combine series
constructed with different $J$'s.
In order to get a more flexible framework we make
in the $\gG$-series (\ref{eq:Gamma series}) the substitution of variables
$u_j=e^{2\pi \ci w_j}$ for $j=1,\ldots,N$.
We write $\vw=(w_1,\ldots,w_N)$, 
$\ugamma=( \gamma_1,\ldots, \gamma_N)$ and 
$\ull=(\ell_1,\ldots,\ell_N)$. We also use the dot-product:
$$
\vw\cdot\ull=w_1\ell_1+w_2\ell_2+\ldots+w_N\ell_N\,.
$$
With these new variables and notations the $\gG$-series 
(\ref{eq:Gamma series}) becomes
\begin{equation}\label{eq:Fourier Gamma series}
\Psi_{\LL,\ugamma}(\vw)\;=\;
\sum_{\ull\in\LL}\:
\frac{e^{2\pi \ci\vw\cdot(\ugamma+\ull)}}
{\prod_{j=1}^N\Gamma(\gamma_j+\ell_j+1)}\,.
\end{equation}

As in Section \ref{subsection:gamma power}
we take a set $J\subset\{1,\ldots,N\}$ with $k+1$ elements, such that the 
vectors $\va_j$ with $j\in J$ are linearly independent
and let $\ugamma=(\gam_1,\ldots,\gam_N)\in\CC^N$ be such that $\gam_j\in\ZZ$ for $j\in J':=\{1,\ldots,N\}\setminus J$.
The vector $\sum_{i\in J'}\va_i$ is a
$\ZZ$-linear combination of the vectors $\va_j$ with 
$j\in J$. Such a relation is an element of $\LL$. Thus one sees that
$\LL$ contains an element $\ull=(\ell_1,\ldots,\ell_N)$ with
$\ell_j=1$ for all $j\in J'$. Since $\gG$-series do not change if one adds to $\ugamma$ an element of $\LL$, we can assume without loss of generality that $\ugamma=(\gam_1,\ldots,\gam_N)\in\CC^N$ is such that
$\gam_j\in\ZZ_{\leq 0}$ for $j\in J':=\{1,\ldots,N\}\setminus J$.
Then the series $\Psi_{\LL,\ugamma}(\vw)$ in (\ref{eq:Gamma series}) involves only terms from the set
\begin{equation}\label{eq:L cone2} 
\LL_J:=\{(\ell_1,\ldots,\ell_N)\in\LL\;|\;
\ell_j\geq 0\quad\textrm{if}\quad j\in J'\}\,.
\end{equation}
Using the estimates (\ref{eq:estimate}) we see that the series $\Psi_{\LL,\ugamma}(\vw)$
converges if the imaginary part $\Im\vw$ of $\vw$ satisfies
$\Im\vw\cdot\ull>\frac{\log T}{2\pi}$ for every non-zero $\ull\in\LL_J$.
We return to this issue and put it an appropriate perspective in Section \ref{subsection:gamma fourier secondary fan}.

\subsection{$\gG$-series and GKZ differential equations.}
\label{subsection:gamma DE}
As in the previous section we consider a $k+1$-element subset 
$J\subset\{1,\ldots,N\}$ such that the 
vectors $\va_j$ with $j\in J$ are linearly independent, and a vector 
$\ugamma=(\gam_1,\ldots,\gam_N)\in\CC^N$ be such that $\gam_j\in\ZZ$ for $j\in J':=\{1,\ldots,N\}\setminus J$.
The $\gG$-series 
(\ref{eq:Gamma series}) constructed with such a $\ugamma$ involves only terms from the set $\LL_{J,\ugamma}$ (see (\ref{eq:L cone})).
For $M\in\NN$ we define the $M$-th partial $\gG$-series 
$\Phi_{\LL,\ugamma,M}(u_1,\ldots,u_N)$ to be the 
subseries of (\ref{eq:Gamma series}) consisting of the terms with
$|\ell_1|+\ldots+|\ell_N|\leq M$.
Then it follows, as in Proposition \ref{prop:power convergence}
from Proposition \ref{prop:growth} and
Lemma \ref{lemma:LJ},
that the sequence $\{\Phi_{\LL,\ugamma,M}(u_1,\ldots,u_N)\}_{M\in\NN}$ converges for $M\to\infty$ to $\Phi_{\LL,\ugamma}(u_1,\ldots,u_N)$
if $|u_j|\leq (2T)^{-\beta}$ for $j\in J'$ and $\frac{1}{2}\leq |u_j|\leq 2$ for $j\in J$.
So on this domain the $\gG$-series
$\Phi_{\LL,\ugamma}(u_1,\ldots,u_N)$ becomes a function of $(u_1,\ldots,u_N)$ that can be differentiated term by term.
This shows
\begin{itemize}
\item
for $(\gl_1,\ldots,\gl_N)\in\LL$
\begin{eqnarray*}
&&\prod_{\gl_i<0}\left(\frac{\partial}{\partial u_i}\right)^{-\gl_i}
\Phi_{\LL,\ugamma}
\;=\sum_{(\ell_1,\ldots,\ell_N)\in\LL}\prod_{j=1}^N
\frac{u_j^{\gamma_j+\ell_j+\min(0,\gl_j)}}{\Gamma(\gamma_j+\ell_j+1+\min(0,\gl_j))}\;=
\\
&&=\sum_{(\ell_1,\ldots,\ell_N)\in\LL}\prod_{j=1}^N
\frac{u_j^{\gamma_j+\ell_j+\gl_j-\max(0,\gl_j)}}{\Gamma(\gamma_j+\ell_j
+\gl_j+1-\max(0,\gl_j))}
\;=\;
\prod_{\gl_i>0}\left(\frac{\partial}{\partial u_i}\right)^{\gl_i}
\Phi_{\LL,\ugamma}\,.
\end{eqnarray*}
\item
for $(a_1,\ldots,a_N)\in\ZZ^N$ such that $\sum_{j=1}^N a_j\ell_j=0$ for every $(\ell_1,\ldots,\ell_N)\in\LL$:
$$
\begin{array}{rcl}
\displaystyle{\sum_{j=1}^N
a_j\,u_j\frac{\partial  \Phi_{\LL,\ugamma}}{\partial u_j}}&=&
\displaystyle{\sum_{(\ell_1,\ldots,\ell_N)\in\LL}
\left(\sum_{j=1}^N a_j(\gamma_j+\ell_j)\right)
\prod_{j=1}^N
\frac{u_j^{\gamma_j+\ell_j}}{\Gamma(\gamma_j+\ell_j+1)}}
\\[3ex]
&=&\displaystyle{(\sum_{j=1}^N a_j\gamma_j)\Phi_{\LL,\ugamma}}
\end{array}
$$
\end{itemize}
The latter system of differential equations is equivalent with the system (\ref{eq:GKZ2}) with $\vc=\sum_{j=1}^N \gamma_j\va_j$.
This shows:
\begin{proposition}\label{prop:gamma DE}
As a function on its domain of convergence $\Phi_{\LL,\ugamma}$ satisfies all differential equations of the GKZ system associated with $\gA$ and $\vc=\sum_{j=1}^N \gamma_j\va_j$.
\qed
\end{proposition}

Note that the $\gG$-series $\Phi_{\LL,\ugamma}$ does not change 
if one adds to $\ugamma$ an element of $\LL$ whereas
the differential equations (\ref{eq:GKZ2}) with 
$\vc=\sum_{j=1}^N \gamma_j\va_j$ do not change if one adds to 
$\ugamma$ an element of $\LL\otimes\CC$.

\section{The Secondary Fan}
\label{section:secondary fan}

As before we consider a subset $\gA=\{\va_1,\ldots,\va_N\}\subset\ZZ^{k+1}$ which generates $\ZZ^{k+1}$ as an abelian group and for which there exists a group homomorphism 
$h:\ZZ^{k+1}\rightarrow\ZZ$ such that $h(\va)=1$ for all $\va\in\gA$.
Still as before, we write
$$
\LL:=\{(\ell_1,\ldots,\ell_N)\in\ZZ^N\;|\;
\ell_1\va_1+\ldots+\ell_N\va_N\,=\,\mf{0}\}\,,
$$
and note that
$\ell_1+\ldots+\ell_N=0$ for every  
$(\ell_1,\ldots,\ell_N)\in\LL$.
In order to better keep track of the various spaces involved we write
$\MM$ instead of $\ZZ^{k+1}$.
Thus the input data is a short exact sequence
\begin{equation}\label{eq:lattice exact sequence}
0\longrightarrow\LL\longrightarrow\ZZ^N\longrightarrow\MM\longrightarrow 0
\,.
\end{equation}
The vectors $\va_1,\ldots,\va_N\in\MM$ are the images of the standard basis vectors of $\ZZ^N$.
We set $d:=\rank \LL$ and $k+1:=\rank\MM=N-d$.

Apart from the common input data this section is independent of the sections on GKZ-systems and $\gG$-series. It concentrates on geometric and combinatorial structures associated with $\gA$ (or equivalently $\LL$).

\subsection{Construction of the secondary fan}
\label{subsection:secondary fan construction}

We write $\LL_\RR^\vee:=\mathrm{Hom}(\LL,\RR)$,
$\MM_\RR^\vee:=\mathrm{Hom}(\MM,\RR)$ and
identify $\mathrm{Hom}(\ZZ^N,\RR)$ and $\RR^N$ via the standard bases. The $\RR$-dual of the exact sequence (\ref{eq:lattice exact sequence}) is
\begin{equation}\label{eq:dual exact sequence}
0\longrightarrow\MM_\RR^\vee\longrightarrow\RR^N
\stackrel{\pi}{\longrightarrow}\LL_\RR^\vee\longrightarrow 0
\,.
\end{equation}
Let $\cP:=\{(x_1,\ldots,x_N)\in\RR^N\;|\; x_i\geq 0\,,\;\forall i\}$ be the positive orthant in $\RR^N$ and let
\begin{equation}\label{eq:polytope map}
\hp:\cP\longrightarrow\LL_\RR^\vee
\end{equation}
denote the restriction of $\pi$.
Since the vector $(1,1,\ldots,1)$ lies in
$\ker\pi$ the map $\hp$ is also surjective.

\begin{figure}[b]
\begin{picture}(230,200)(0,0)
\put(0,0){\line(1,0){240}}
\put(0,0){\line(0,1){180}}
\multiput(0,0)(15,5){15}{\line(3,1){10}}
\put(230,100){\makebox{$\hp^{-1}(1)$}}
\put(160,170){\makebox{$\hp^{-1}(-1)$}}
\thicklines
\put(0,50){\line(2,1){250}}
\multiput(0,50)(18,15){9}{\line(6,5){10}}
\put(100,0){\line(2,1){150}}
\multiput(100,33)(18,15){9}{\line(6,5){10}}
\multiput(100,0)(0,13){3}{\line(0,1){8}}
\end{picture}
\caption{\label{fig:polytope family} Fibres $\hp^{-1}(1)$
and $\hp^{-1}(-1)$ for $\LL=\ZZ(-2,1,1)\subset\RR^3$}
\end{figure}

\

\noindent
\textbf{Example.}
Take $\LL=\ZZ(-2,1,1)\subset\RR^3$. Then $\pi$ can be identified with the map
$$
\pi:\RR^3\longrightarrow\RR\,,\qquad \pi(x_1,x_2,x_3)=-2x_1+x_2+x_3\,.
$$
For $t\in\RR$ the polytope $\hp^{-1}(t)$ is the intersection of the positive octant and the plane with equation $-2x_1+x_2+x_3=t$.
Figure \ref{fig:polytope family} illustrates this for $t=1$ and $t=-1$
(with the $x_1$-axis drawn vertically).

\

Let $\bb_1,\ldots,\bb_N\in  \LL_\RR^\vee$ be the images of the standard basis vectors of $\RR^N$ under the map $\pi$. Then,
for $\mf{t}\in  \LL_\RR^\vee$, 
$$
(x_1,\ldots,x_N)\in\hp^{-1}(\mf{t})\qquad\Longleftrightarrow\qquad
\mf{t}=x_1\bb_1+\ldots+x_N\bb_N \quad\textrm{and}\quad 
x_i\geq 0\,,\;\forall i\,.
$$
We see that the fiber $\hp^{-1}(\mf{t})$ is a convex (unbounded) polyhedron.

\begin{lemma}
$\vv=(v_1,\ldots,v_N)\in\cP$ is a vertex of $\hp^{-1}(\mf{t})$ if and only if
$\mf{t}\,=\,\sum_{j=1}^N v_j\bb_j$ and
the vectors $\bb_j$ with $v_j\neq 0$
are linearly independent over $\RR$.
\end{lemma}
\begin{proof}
Suppose $\mf{t}\,=\,\sum_{j=1}^N v_j\bb_j$, 
all $v_j\geq 0$ and
the vectors $\bb_j$ with $v_j\neq 0$
are linearly \emph{dependent} over $\RR$.
Then there is a non-trivial relation $\sum_{j=1}^N x_j\bb_j=0$ with $|x_j|\leq v_j$
for all $j$ and the whole interval 
$\{\vv+s(x_1,\ldots,x_N)\:|\:|s|\leq 1\}$ lies in $\hp^{-1}(\mf{t})$.
Therefore $\vv$, being the midpoint of this interval, can not be a vertex of $\hp^{-1}(\mf{t})$.

Suppose $\vv=(v_1,\ldots,v_N)\in\hp^{-1}(\mf{t})$ is not a vertex of $\hp^{-1}(\mf{t})$. Then there is a non-zero vector 
$\vx=(x_1,\ldots,x_N)\in\RR^N$ such that the interval 
$\{\vv+s\vx\:|\:|s|\leq 1\}$ lies in $\hp^{-1}(\mf{t})=\cP\cap\pi^{-1}(\mf{t})$.
This implies $|x_j|\leq v_j$ for all $j$ and $\sum_{j=1}^N x_j\bb_j=0$.
Consequently, the vectors $\bb_j$ with $v_j\neq 0$
are linearly \emph{dependent} over $\RR$.
\end{proof}

For a vertex $\vv=(v_1,\ldots,v_N)$ of $\hp^{-1}(\mf{t})$ we set 
\begin{equation}\label{eq:Iv}
I_\vv:=\{i\;|\;v_i=0\}\subset\{1,2,\ldots,N\}\,.
\end{equation}
In this way every $\mf{t}\in  \LL_\RR^\vee$ yields a list of subsets
of $\{1,2,\ldots,N\}$:
\begin{equation}\label{eq:vertex list}
T_\mf{t}:=\{I_\vv\;|\; \vv\textrm{ vertex of } \hp^{-1}(\mf{t})\}.
\end{equation}
Since $\pi^{-1}(\mf{t})$ has dimension $N-d$, the cardinality of each $I_\vv$ must be at least $N-d$.

The above lemma provides an alternative description of the list $T_\mf{t}$:
\begin{corollary}\label{cor:list}
A subset $I\subset\{1,\ldots,N\}$ is on
the list $T_\mf{t}$ if and only if
the vectors $\bb_j$ with $j\not\in I$
are linearly independent over $\RR$ and
$\displaystyle{\mf{t}=\sum_{j\not\in I} \tau_j\bb_j}$ with all $\tau_j\in\RR_{>0}$.
\qed
\end{corollary}

\

We now define an equivalence relation on $  \LL_\RR^\vee$ by:
$\;\mf{t}\sim \mf{t}'\quad\Longleftrightarrow\quad T_\mf{t}=T_{\mf{t}'}$.
\\
From Corollary \ref{cor:list} one sees that the equivalence class containing $\mf{t}$ is
\begin{equation}\label{eq:equivclass}
\cC=\bigcap_{I\in T_{\mf{t}}}(\textrm{positive span of }\{\bb_i\}_{i\not\in I}).
\end{equation}
So the equivalence classes are
strongly convex polyhedral cones in $\LL_\RR^\vee$.

\

\begin{definition}\label{def:secondary fan}
 This collection of cones is called the \emph{secondary fan
 of $\gA$ (or $\LL$)}.
\end{definition}

\

For an equivalence class $\cC$ we set
$T_{\cC}:=T_\mf{t}$ for any $\mf{t}\in \cC$. 
It follows from (\ref{eq:equivclass}) that an equivalence class $\cC$ is an open cone of dimension
$d$ if and only if all sets on the list $T_\cC$ have
exactly $N-d$ elements.

\

\noindent
\textbf{Example.}
In the example of $\LL=\ZZ(-2,1,1)\subset\RR^3$ (see Figure \ref{fig:polytope family}) the vertices are given by the
lists 
$$
T_t=\left\{\begin{array}{lll}
\{\{2,3\}\}&\textrm{ if }& t<0\\
\{\{1,2,3\}\}&\textrm{ if }& t=0\\
\{\{1,2\},\{1,3\}\}&\textrm{ if }& t>0
\end{array}\right.
$$

\

\noindent
\textbf{Example.}
For Gauss's hypergeometric structures $\LL=\ZZ(1,1,-1,-1)\subset\RR^4$ and, hence, $\bb_1=\bb_2=1$,
$\bb_3=\bb_4=-1$ in $\RR$. Corollary \ref{cor:list} now yields the lists
$$
T_t=\left\{\begin{array}{lll}
\{\{1,2,3\},\{1,2,4\}\}&\textrm{ if }& t<0\\
\{\{1,2,3,4\}\}&\textrm{ if }& t=0\\
\{\{2,3,4\},\{1,3,4\}\}&\textrm{ if }& t>0
\end{array}\right.
$$

\

\noindent
\textbf{Example.}
For Appell's $F_1$ the lattice $\LL\subset\ZZ^6$ has rank $2$ and is 
generated by the two vectors $(1,1,0,-1,-1,0)$ and $(1,0,1,-1,0,-1)$
which express that the three vertical segments
in Figure \ref{fig:classical examples} are parallel. The vectors $\bb_1,\ldots,\bb_6\in\ZZ^2$ are therefore  
\begin{figure}[t]
\begin{picture}(230,160)(-100,-75)
\put(0,0){\line(1,0){70}}
\put(0,0){\line(-1,0){70}}
\put(0,0){\line(0,1){70}}
\put(0,0){\line(0,-1){70}}
\put(0,0){\line(1,1){70}}
\put(0,0){\line(-1,-1){70}}
\put(75,-3){\makebox(10,10){$\bb_5$}}
\put(-85,-3){\makebox(10,10){$\bb_2$}}
\put(-7,75){\makebox(10,10){$\bb_6$}}
\put(-1,-85){\makebox(10,10){$\bb_3$}}
\put(73,73){\makebox(10,10){$\bb_1$}}
\put(-83,-83){\makebox(10,10){$\bb_4$}}
\put(-50,10){\makebox(40,30){$\begin{array}{l}
\{3,4,5,6\}\\ \{1,3,4,5\}\\ \{1,2,3,5\}\end{array}$}}
\put(10,-40){\makebox(40,30){$\begin{array}{l}
\{2,4,5,6\}\\ \{1,2,3,6\}\\ \{1,2,4,6\}\end{array}$}}
\put(42,10){\makebox(40,30){$\begin{array}{l}
\{2,4,5,6\}\\ \{1,2,3,4\}\\ \{2,3,4,6\}\end{array}$}}
\put(-82,-40){\makebox(40,30){$\begin{array}{l}
\{1,4,5,6\}\\ \{1,2,3,5\}\\ \{1,3,5,6\}\end{array}$}}
\put(-45,-78){\makebox(40,30){$\begin{array}{l}
\{1,4,5,6\}\\ \{1,2,5,6\}\\ \{1,2,3,6\}\end{array}$}}
\put(7,49){\makebox(40,30){$\begin{array}{l}
\{3,4,5,6\}\\ \{1,2,3,4\}\\ \{2,3,4,5\}\end{array}$}}
\end{picture}
\caption{\label{fig:fan F1}Secondary fan for $F_1$}
\end{figure}
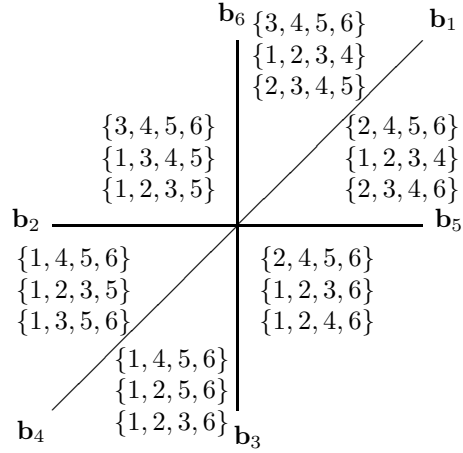
$$
\bb_1=\VE{1\\ 1}\,,\;\bb_2=\VE{-1\\ 0}\,,\;\bb_3=\VE{0\\ -1}\,,\;
\bb_4=\VE{-1\\ -1}\,,\;\bb_5=\VE{1\\ 0}\,,\;\bb_6=\VE{0\\ 1}\,.
$$
Figure \ref{fig:fan F1} shows the secondary fan for $F_1$ and gives for each maximal cone $\cC$ the corresponding list $T_{\cC}$ according
to Corollary \ref{cor:list}.

\

\begin{figure}[t]
\begin{picture}(230,130)(-100,-65)
\put(0,0){\line(-1,0){50}}
\put(0,0){\line(0,-1){50}}
\put(0,0){\line(1,1){50}}
\put(-75,-3){\makebox(10,10){$\bb_5=\bb_6$}}
\put(-3,-65){\makebox(10,10){$\bb_2=\bb_4$}}
\put(53,53){\makebox(10,10){$\bb_1=\bb_3$}}
\put(-30,25){\makebox(40,30){$\begin{array}{l}
\{2,3,4,6\}\\ \{2,3,4,5\}\\ \{1,2,4,6\}\\ \{1,2,4,5\}\end{array}$}}
\put(20,-30){\makebox(40,30){$\begin{array}{l}
\{3,4,5,6\}\\ \{2,3,5,6\}\\ \{1,4,5,6\}\\ \{1,2,5,6\}\end{array}$}}
\put(-45,-42){\makebox(40,30){$\begin{array}{l}
\{1,3,4,6\}\\ \{1,3,4,5\}\\ \{1,2,3,5\}\\ \{1,2,3,6\}\end{array}$}}
\end{picture}
\caption{\label{fig:fan F4} Secondary fan for $F_4$}
\end{figure}
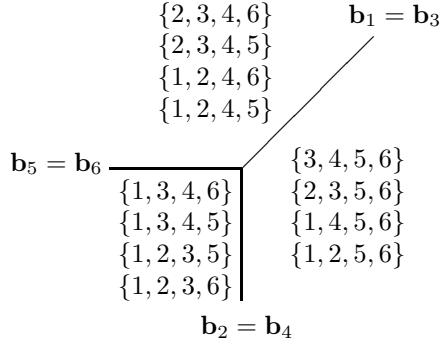

\noindent
\textbf{Example.} 
For Appell's $F_4$ the lattice $\LL\subset\ZZ^6$ has rank $2$ and is 
generated by the two vectors $(1,-1,1,-1,0,0)$ and $(1,0,1,0,-1,-1)$
which express that the three diagonals
in Figure \ref{fig:classical examples} intersect at the centre. The vectors $\bb_1,\ldots,\bb_6\in\ZZ^2$ are   
$$
\bb_1=\bb_3=\VE{1\\ 1}\,,\;\bb_2=\bb_4=\VE{-1\\ 0}\,,\;
\bb_5=\bb_6=\VE{0\\ -1}\,.
$$
Figure \ref{fig:fan F4} shows the secondary fan for $F_4$ and gives for each maximal cone $\cC$ the corresponding list $T_{\cC}$ according
to Corollary \ref{cor:list}.

\

\noindent
\textbf{Example/exercise.} The reader is invited to apply the techniques
demonstrated in the previous examples to the examples in 
Section \ref{subsection:GKZCY} Table \ref{table CY}.

\subsection{Alternative descriptions for secondary fan constructions.}
\label{section:alternatives}
We are going to present geometrically appealing alternative descriptions
for the polyhedra $\hp^{-1}(\mf{t})$ and for the lists $T_{\cC}$ associated with the maximal cones in the secondary fan.
Whereas the constructions
in Section \ref{subsection:secondary fan construction} were completely presented in terms of $\LL$, the alternative descriptions use $\gA$ only.

\subsubsection{Piecewise linear functions associated with $\gA$}
\label{subsection:Polyhedra}
The vectors $\va_1,\ldots,\va_N\in\MM$ are linear functions on the space
$\MM_\RR^\vee:=\mathrm{Hom}(\MM,\RR)$. Let $\MM_\RR:=\MM\otimes\RR$
and denote the pairing between
$\MM_\RR$ and $\MM_\RR^\vee$ by $\langle,\rangle$.
The inclusion $\MM_\RR^\vee\hookrightarrow\RR^N$ is then given by
\begin{equation}\label{eq:M inclusion}
\MM_\RR^\vee\hookrightarrow\RR^N\,,\qquad
\vv\mapsto (\langle \va_1,\vv\rangle,\ldots,\langle \va_N,\vv\rangle)\,.
\end{equation}
For an $N$-tuple $\uga=(\ga_1,\ldots,\ga_N)\in\RR^N$ one has the polyhedron
\begin{equation}\label{eq:A polyhedron}
K_\uga:=\{\vv\in\MM_\RR^\vee\:|\:
\langle \va_j,\vv\rangle\geq -\ga_j,\;\forall j\}\,.
\end{equation}
Recall that  $\bb_1,\ldots,\bb_N\in  \LL_\RR^\vee$ denote the images of the standard basis vectors of $\RR^N$ under the map $\pi$.

\begin{proposition}\label{prop:equivalent polyhedra}
For $\uga=(\ga_1,\ldots,\ga_N)\in\RR^N$ set 
$\mf{t}=\ga_1\bb_1+\ldots+\ga_N\bb_N$.
Then
$$
\hp^{-1}(\mf{t})\:=\:\uga\,+\,K_\uga\,.
$$
\end{proposition}
\begin{proof}
Since $\uga$ is in $\pi^{-1}(\mf{t})$ a point
$\vx$ is in $\pi^{-1}(\mf{t})$ if and only if
$\vx-\uga$ is in $\ker\pi=\MM_\RR^\vee$.
By definition, a point $\vx=(x_1,\ldots,x_N)\in\pi^{-1}(\mf{t})$
lies in $\hp^{-1}(\mf{t})$ if and only if $x_j\geq 0$ for all $j$.
Thus, in view of (\ref{eq:M inclusion}),
$$
\vx\in\hp^{-1}(\mf{t})\qquad\Longleftrightarrow\qquad
\vv:=\vx-\uga\textrm{ satisfies }
\langle \va_j,\vv\rangle+\ga_j\geq 0,\;\forall j\,.
$$
\end{proof}

\

Recall that throughout these notes we assume the existence of a group homomorphism 
$h:\ZZ^{k+1}\rightarrow\ZZ$ such that $h(\va)=1$ for all $\va\in\gA$.
In the present terminology this amounts to the existence of an element
$\vh\in\MM_\RR^\vee$ such that $\langle \va_j,\vh\rangle=1$ for $j=1,\ldots,N$. 
Now fix a direct sum decomposition of real vector spaces
\begin{equation}\label{eq:direct sum}
\MM_\RR^\vee\,=\,\MM_\RR^\circ\,\oplus\,\RR\vh
\end{equation}
and consider the function
\begin{equation}\label{eq:graph}
\mu_\uga:\MM_\RR^\circ\longrightarrow\RR\,,\qquad 
\mu_\uga(\vu)=\min_j\left(\langle \va_j,\vu\rangle+\ga_j\right)\,.
\end{equation}

\begin{proposition}\label{prop:boundary K}
For every $\vu\in\MM_\RR^\circ$ the vector $\vu-\mu_\uga(\vu)\vh$ lies
in the boundary $\partial K_\uga$ of $K_\uga$.
In other words $\partial K_\uga$ is the graph of the function $-\mu_\uga$ on
$\MM_\RR^\circ$.
\end{proposition}
\begin{proof}
Take $\vu\in\MM_\RR^\circ$. Then one checks for every $j$ 
$$
\langle \va_j,\vu-\mu_\uga(\vu)\vh\rangle=
\langle \va_j,\vu\rangle-\mu_\uga(\vu)\geq
\langle \va_j,\vu\rangle-\left(\langle \va_j,\vu\rangle+\ga_j\right)
= -\ga_j.
$$
So $\vu-\mu_\uga(\vu)\vh$ lies in $K_\uga$. If $j$ is such that
$\mu_\uga(\vu)=\langle \va_j,\vu\rangle+\ga_j$, then the above computation
shows $\langle \va_j,\vu-\mu_\uga(\vu)\vh\rangle=-\ga_j$. Therefore
$\vu-\mu_\uga(\vu)\vh$ lies in $\partial K_\uga$.
\end{proof}

If $\mu_\uga(\vu)=\langle \va_j,\vu\rangle+\ga_j$, then
the point $\vu-\mu_\uga(\vu)\vh$ lies in the affine hyperplane
\begin{equation}\label{eq:A hyperplanes}
\cH^\uga_j:= \{\vv\in\MM_\RR^\vee\:|\:
\langle \va_j,\vv\rangle= -\ga_j\}\,,\qquad j=1,\ldots,N\,.
\end{equation}
For generic $\vu\in\MM_\RR^\circ$ (i.e. outside some codimension $1$
closed subset) the minimum in (\ref{eq:graph}) is attained for exactly one $j$. Therefore each codimension $1$ face of the polyhedron $K_\uga$
lies in some unique hyperplane $\cH^\uga_j$.

\

\noindent
\textbf{Remark.} $K_\uga$ can also be described as the closure of that connected component
of $\MM_\RR^\vee\setminus\bigcup_{j=1}^N\cH^\uga_j$ 
that contains the points $t\vh$ for sufficiently large $t$.

\

\begin{figure}
\begin{center}
\setlength\epsfxsize{12cm}
\epsfbox{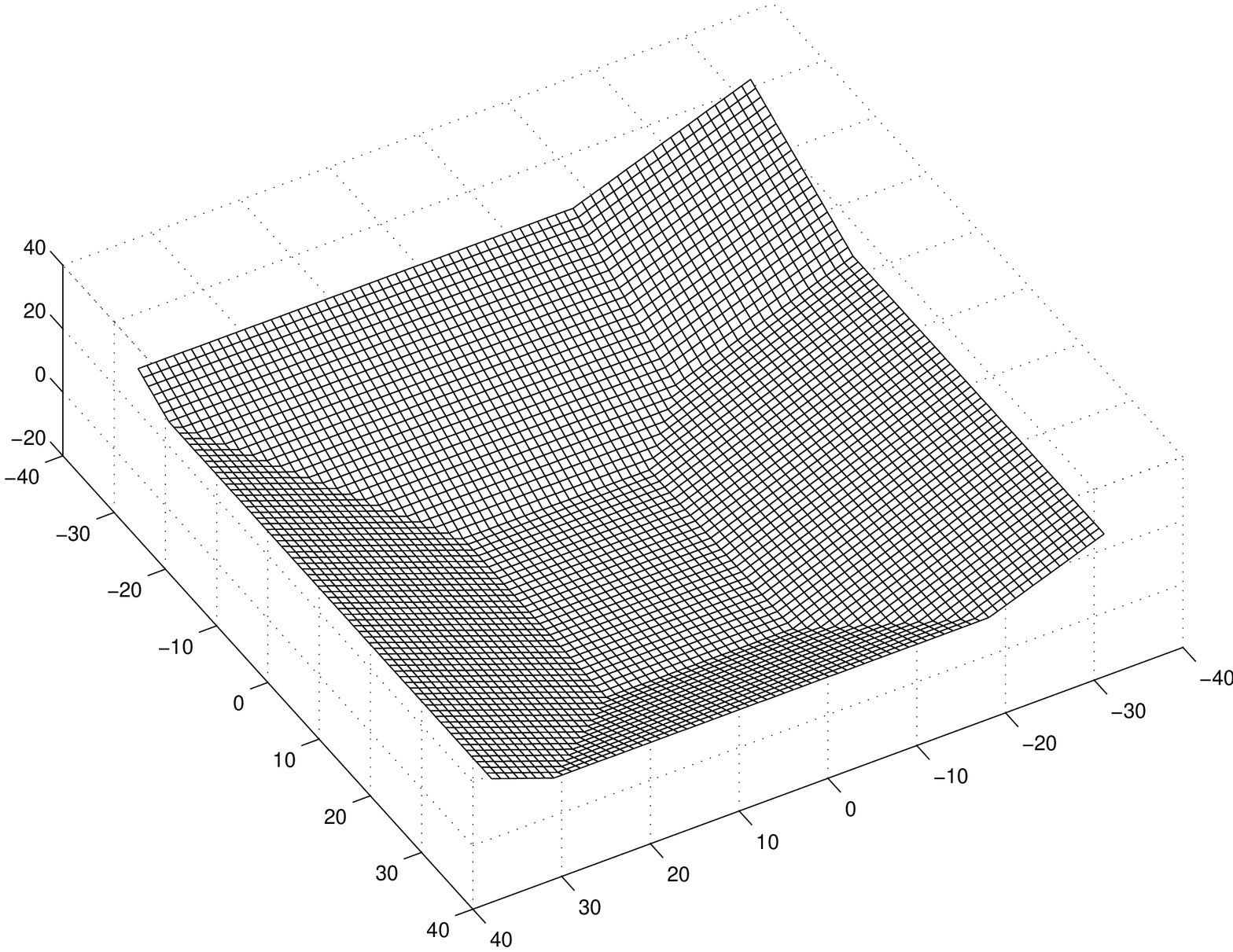}
\end{center}
\caption{\label{fig:polyhedron}
Example of a polyhedron $K_\uga$ for $\gA$ as in (\ref{eq:pentagon A}).}
\end{figure}

\noindent
\textbf{Example.}
Figure \ref{fig:polyhedron} shows (a piece of) 
the polyhedron $K_\uga$ for
\begin{equation}\label{eq:pentagon A}
\gA=\left\{\VE{1\\ 0\\ 1},\,\VE{1\\ 1\\ 1},\,\VE{1\\ -1\\ 0},\,
\VE{1\\ 0\\ 0},\,\VE{1\\ 1\\ 0},\,\VE{1\\ 0\\ -1}\right\}
\end{equation}
and $\uga=(21,35,35,14,21,28)$.
Matching the faces of $K_\uga$ with the vectors in $\gA$
one checks that the list of vertices is
$\{\{1,2,5\},\,\{1,4,5\},\,\{1,3,4\},\,\{4,5,6\},\,\{3,4,6\}\}$.

\subsubsection{Regular triangulations}
\label{subsection:triangulations}
Assume $\uga=(\ga_1,\ldots,\ga_N)\in\RR^N$ with all $\ga_j>0$. Then the dual
of the polyhedron $K_\uga$ in (\ref{eq:A polyhedron}) is, by definition
\begin{equation}\label{eq:dual polyhedron}
K_\uga^\vee:=\{\vw\in\MM_\RR\:|\:\langle\vw,\vv\rangle\,>\,-1\,,\;\forall \vv\in K_\uga\}\,.
\end{equation}

\begin{lemma}\label{lemma:dual polyhedron}
$K_\uga^\vee\,=\,\mathrm{convex\;hull}\;\{\mf{0},\,\frac{1}{\ga_1}\va_1,
\ldots,\frac{1}{\ga_N}\va_N\}$.
\end{lemma}
\begin{proof}
The inclusion $\supset$ follows directly from the definition of $K_\uga$
in (\ref{eq:A polyhedron}).
Now suppose that the two polyhedra are not equal. Then there is a point $\vp$
in
$K_\uga^\vee$ which is separated by an affine hyperplane from
$\mf{0},\,\frac{1}{\ga_1}\va_1,\ldots,\frac{1}{\ga_N}\va_N$.
That means that there is a vector $\vv\in\MM_\RR^\vee$, perpendicular to the 
hyperplane, such that
$\langle\vp,\vv\rangle<-1$ and 
$\langle\frac{1}{\ga_j}\va_1,\vv\rangle>-1$ for $j=1,\ldots,N$.
The last $N$ inequalities imply according to (\ref{eq:A polyhedron})
that $\vv$ is in $K_\uga$, but then the first inequality contradicts
$\vp\in K_\uga^\vee$.
So we conclude that the two polyhedra are equal.
\end{proof}

Next we use the projection from the point $\mf{0}$ to project 
$K_\uga^\vee$ into the hyperplane with equation 
$\langle\vp,\vh\rangle=1$. This maps $K_\uga^\vee$ onto
the polytope
\begin{equation}\label{eq:Delta}
\Delta_\gA:=\textrm{convex hull}\{\va_1,\ldots,\va_N\}.
\end{equation}
The images of the codimension $1$ faces of $K_\uga^\vee$ which
do not contain the vertex $\mf{0}$ induce a subdivision
of $\Delta_\gA$ by the polytopes
\begin{equation}\label{eq:subdivision}
\textrm{convex hull}\{\va_i\}_{i\in I} \qquad \textrm{for}\quad I\in 
T_\mf{t}\,,
\end{equation}
where $\mf{t}=\ga_1\bb_1+\ldots+\ga_N\bb_N$ as in  
Proposition \ref{prop:equivalent polyhedra} and
$T_\mf{t}$ is the corresponding list of vertices of 
$\hp^{-1}(\mf{t})$ as in (\ref{eq:vertex list}).

\setlength{\unitlength}{2pt}
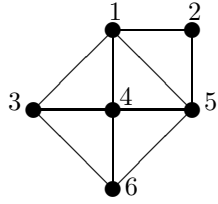
\begin{figure}[t]
\begin{center}
\begin{picture}(80,40)(-10,40)
\put(15,40){\line(1,1){15}}
\put(15,40){\line(-1,1){15}}
\put(15,70){\line(1,0){15}}
\put(15,70){\line(-1,-1){15}}
\put(30,55){\line(0,1){15}}
\put(15,40){\circle*{3}}
\put(0,55){\circle*{3}}
\put(15,55){\circle*{3}}
\put(30,55){\circle*{3}}
\put(15,70){\circle*{3}}
\put(30,70){\circle*{3}}
\put(15,55){\line(0,-1){15}}
\put(15,55){\line(0,1){15}}
\put(0,55){\line(1,0){30}}
\put(30,55){\line(-1,1){15}}
\put(14,72){\makebox(3,3){$1$}}
\put(29,72){\makebox(3,3){$2$}}
\put(-5,55){\makebox(3,3){$3$}}
\put(16,56){\makebox(3,3){$4$}}
\put(32,55){\makebox(3,3){$5$}}
\put(17,39){\makebox(3,3){$6$}}
\end{picture}
\end{center}
\caption{\label{fig:triangulation}
Triangulation corresponding with Figure \ref{fig:polyhedron}}
\end{figure}
\setlength{\unitlength}{1pt}

If the point $\mf{t}$ lies in some maximal cone $\cC$ of the secondary fan,
all members of the list $T_\mf{t}=T_\cC$ have $N-d=k+1$ elements.
The polytopal subdivision of $\Delta_\gA$ is then a \emph{triangulation};
i.e. all polytopes in the subdivision (\ref{eq:subdivision})
are $k$-dimensional simplices.

\begin{definition}\label{def:regular triangulations}
The triangulations of $\Delta_\gA$ obtained in this way are called
\emph{regular triangulations}.
\end{definition}

\begin{definition}\label{def:unimodular triangulations}
One defines the \emph{volume} of a $k$-dimensional simplex
with vertex set $\{\ba_i\}_{i\in I}$ to be
\begin{equation}\label{eq:volume}
\mathrm{volume}(\{\ba_i\}_{i\in I})\;=\;
|\det\left((\ba_i)_{i\in I}\right)|\,.
\end{equation}
A regular triangulation of $\Delta_\gA$ is said to be
\emph{unimodular} if all $k$-dimensional simplices in the
triangulation have volume equal to $1$.
\end{definition}

By abuse of language we will just say ``the triangulation
$T_\cC$'' instead of ``the triangulation corresponding with the maximal cone $\cC$''.
Note that useful information about distances between vertices of
$\hp^{-1}(\mf{t})$ gets lost in the passage to the (purely combinatorial) 
triangulation $T_\cC$.

\

\noindent\textbf{Remark.}
In general there can be
triangulations of $\Delta_\gA$ with vertices
in $\{\ba_1,\ldots,\ba_N\}$, which do not arise from the above construction and are therefore not regular.

\subsection{The Secondary Polytope}
\label{subsection:Secondary Polytope}
The $k$-dimensional polytope $\Delta_\gA$
defined in (\ref{eq:Delta}) is sometimes called the \emph{primary
polytope} associated with $\gA$.
By definition, the regular triangulations of $\Delta_\gA$ correspond bijectively with the maximal cones of the secondary fan.
To a regular triangulation $T_\cC$ we assign the point $q_\cC\in\RR^N$
with 
$$
j^{\mathrm{th}}\textrm{-coordinate of }q_\cC\;=\;
\sum_{I\in T_\cC\:\mathrm{s.t.}\:j\in I}
\mathrm{volume}(\{\ba_i\}_{i\in I})\,,
$$
i.e. the sum of the volumes of the simplices in $T_\cC$ of which 
$\ba_j$ is a vertex.

\begin{definition}\label{def:secondary}
The \emph{secondary polytope} associated with $\gA$ is 
$$
\mathrm{Sec}(\gA)\;=\;\textrm{convex hull }\{q_\cC\:|\:
T_\cC \textrm{ regular triangulation of }\Delta_\gA\;\}
$$
\end{definition}

The map $\RR^N\rightarrow\MM_\RR$ maps the $j$-th standard
basis vector of $\RR^N$ to $\va_j$. Thus the point $q_\cC$ is
mapped to
\begin{eqnarray*}
\sum_{j=1}^N
\sum_{I\in T_\cC\:\mathrm{s.t.}\:j\in I}
\mathrm{volume}(\{\ba_i\}_{i\in I})\,\va_j&=&
\sum_{I\in T_\cC} \mathrm{volume}(\{\ba_i\}_{i\in I})
\left(\sum_{j\in I}\va_j\right)\\
&&\hspace{-3em}=\;(k+1)\times\mathrm{volume}(\Delta_\gA)\times
\mathrm{barycenter}(\Delta_\gA)\,.
\end{eqnarray*}
So the whole secondary polytope is mapped to one point.
Therefore, after some translation in $\RR^N$ we find the secondary polytope in $\LL$:
$$
\mathrm{Sec}(\gA)\subset\LL\otimes\RR\,.
$$

As for the relation between secondary fan and secondary polytope we 
mention the following theorem, which is in a slightly
different formulation proven in \cite{gkz3}.

\begin{theorem}\textup{(\cite{gkz3} p.221 thm.1.7)}
The secondary fan, which lies in $\LL_\RR^\vee$, is in fact the fan of outward pointing vectors perpendicular to the faces of $\mathrm{Sec}(\gA)$.
\qed
\end{theorem}

\

\noindent
\textbf{Example.}
In the example of $\LL=\ZZ(-2,1,1)\subset\ZZ^3$ there are two maximal cones: $\RR_{>0}$ and $\RR_{<0}$.
The corresponding triangulations are:

\begin{picture}(230,30)(0,0)
\put(55,20){\makebox(10,5){$t<0$}}
\put(30,10){\line(1,0){60}}
\put(30,10){\circle*{5}}
\put(90,10){\circle*{5}}
\put(28,0){\makebox(5,5){$2$}}
\put(88,0){\makebox(5,5){$3$}}
\put(205,20){\makebox(10,5){$t>0$}}
\put(180,10){\line(1,0){60}}
\put(180,10){\circle*{5}}
\put(210,10){\circle*{5}}
\put(240,10){\circle*{5}}
\put(178,0){\makebox(5,5){$2$}}
\put(208,0){\makebox(5,5){$1$}}
\put(238,0){\makebox(5,5){$3$}}
\end{picture}

\noindent
The secondary polytope is the line segment between the points
$(0,2,2)$ and $(2,1,1)$ in $\RR^3$.

\

\noindent
\textbf{Example.}
For Gauss's hypergeometric structures $\LL=\ZZ(1,1,-1,-1)\subset\ZZ^4$.
From this one sees that $\LL_\RR^\vee\simeq\RR$ and 
that there are two maximal cones: $\RR_{>0}$ and $\RR_{<0}$.
The corresponding triangulations are:

\begin{picture}(230,65)(-60,0)
\put(10,0){
\begin{picture}(100,80)(0,0)
\put(18,50){\makebox(10,5){$t<0$}}
\put(0,0){\line(1,1){40}}
\put(0,0){\line(1,0){40}}
\put(0,0){\line(0,1){40}}
\put(40,40){\line(-1,0){40}}
\put(40,40){\line(0,-1){40}}
\put(0,0){\circle*{5}}
\put(40,0){\circle*{5}}
\put(0,40){\circle*{5}}
\put(40,40){\circle*{5}}
\put(-10,0){\makebox(5,5){$1$}}
\put(47,0){\makebox(5,5){$3$}}
\put(-10,40){\makebox(5,5){$4$}}
\put(47,40){\makebox(5,5){$2$}}
\end{picture}}
\put(140,0){
\begin{picture}(100,80)(0,0)
\put(18,50){\makebox(10,5){$t>0$}}
\put(40,0){\line(-1,1){40}}
\put(0,0){\line(1,0){40}}
\put(0,0){\line(0,1){40}}
\put(40,40){\line(-1,0){40}}
\put(40,40){\line(0,-1){40}}
\put(0,0){\circle*{5}}
\put(40,0){\circle*{5}}
\put(0,40){\circle*{5}}
\put(40,40){\circle*{5}}
\put(-10,0){\makebox(5,5){$1$}}
\put(47,0){\makebox(5,5){$3$}}
\put(-10,40){\makebox(5,5){$4$}}
\put(47,40){\makebox(5,5){$2$}}
\end{picture}}
\end{picture}

\noindent
The secondary polytope is the line segment between the points
$(2,2,1,1)$ and $(1,1,2,2)$ in $\RR^4$.

\

\noindent
\textbf{Example: }
For $\LL=\ZZ(-3\,1\,1\,1)\subset\ZZ^4$ we have
$\LL_\RR^\vee\simeq\RR$ and $\bb_1=-3$, $\bb_2=\bb_3=\bb_4=1$.
Corollary \ref{cor:list} shows for $t\in\LL_\RR^\vee\simeq\RR$:
$$
T_t=\left\{\begin{array}{lll}
\{\;\{1,3,4\}\,,\:\{1,2,4\}\,,\:\{1,2,3\}\;\}&\textrm{ if }& t>0\\
\{\;\{1,2,3,4\}\;\}&\textrm{ if }& t=0\\
\{\;\{2,3,4\}\;\}&\textrm{ if }& t<0
\end{array}\right.
$$
So there are two maximal cones: $\RR_{>0}$ and $\RR_{<0}$.
In terms of triangulations:

\setlength{\unitlength}{1.4pt}

\begin{picture}(230,50)(25,20)
\put(30,15){\makebox(100,50){
$
\gA=\left\{\VE{1\\ 0\\ 0},\,\VE{1\\ 1\\ 1},\,\VE{1\\ -1\\ 0},\,
\VE{1\\ 0\\ -1}\right\}
$
}}
\put(90,-2){\makebox(100,50){
\put(35,65){\makebox(100,5){$t>0$}}
\put(95,40){\line(-1,1){15}}
\put(95,40){\circle*{3}}
\put(80,55){\circle*{3}}
\put(95,55){\circle*{3}}
\put(110,70){\circle*{3}}
\put(95,55){\line(-1,0){15}}
\put(95,55){\line(0,-1){15}}
\put(95,55){\line(1,1){15}}
\put(110,70){\line(-2,-1){30}}
\put(110,70){\line(-1,-2){15}}
}}
\put(40,-2){\makebox(100,50){
\put(35,65){\makebox(100,5){$t<0$}}
\put(95,40){\line(-1,1){15}}
\put(95,40){\circle*{3}}
\put(80,55){\circle*{3}}
\put(95,55){\circle*{3}}
\put(110,70){\circle*{3}}
\put(110,70){\line(-2,-1){30}}
\put(110,70){\line(-1,-2){15}}
}}
\end{picture}

\setlength{\unitlength}{1pt}

\noindent
The secondary polytope is the line segment between the points
$(0,3,3,3)$ and $(3,2,2,2)$ in $\RR^4$.

\begin{figure}[t]
\begin{center}
\begin{picture}(250,320)(13,-10)
\put(-15,-18){\makebox(60,20){${}_{(2,3,2,4,1,3)}$}}
\put(185,-18){\makebox(60,20){${}_{(2,2,2,5,2,2)}$}}
\put(265,22){\makebox(60,20){${}_{(1,3,3,4,2,2)}$}}
\put(-15,227){\makebox(60,20){${}_{(4,3,2,0,1,5)}$}}
\put(75,267){\makebox(60,20){${}_{(1,5,4,0,1,4)}$}}
\put(265,267){\makebox(60,20){${}_{(1,3,5,0,4,2)}$}}
\put(105,65){\makebox(60,20){${}_{(1,4,3,3,1,3)}$}}
\put(95,182){\makebox(60,20){${}_{(5,1,2,0,3,4)}$}}
\put(230,130){\makebox(60,20){${}_{(3,1,2,4,3,2)}$}}
\put(188,228){\makebox(60,20){${}_{(3,1,4,0,5,2)}$}}
\thicklines
\put(15,32){\line(0,1){162}}
\put(33,15){\line(1,0){162}}
\multiput(32,21)(20,10){3}{\line(2,1){10}}
\put(34,215){\line(1,0){82}}
\put(32,223){\line(2,1){51}}
\put(215,32){\line(0,1){85}}
\put(232,23){\line(2,1){51}}
\put(272,241){\line(2,1){13}}
\multiput(95,75)(0,24){7}{\line(0,1){12}}
\multiput(117,55)(25,0){4}{\line(1,0){12}}
\multiput(218,55)(25,0){2}{\line(1,0){12}}
\put(266,55){\line(1,0){10}}
\put(277,255){\line(-1,0){162}}
\put(295,237){\line(0,-1){162}}
\put(154,215){\line(6,1){88}}
\put(218,153){\line(2,5){29}}
\put(205,145){\line(-1,1){60}}

\thinlines
\put(15,0){\line(1,1){15}}
\put(15,0){\line(-1,1){15}}
\put(15,30){\line(1,0){15}}
\put(15,30){\line(-1,-1){15}}
\put(30,15){\line(0,1){15}}
\put(15,0){\circle*{3}}
\put(0,15){\circle*{3}}
\put(15,15){\circle*{3}}
\put(30,15){\circle*{3}}
\put(15,30){\circle*{3}}
\put(30,30){\circle*{3}}
\put(15,15){\line(-1,0){15}}
\put(15,15){\line(1,1){15}}
\put(15,0){\line(0,1){30}}
\put(15,0){\line(1,2){15}}

\put(15,200){\line(1,1){15}}
\put(15,200){\line(-1,1){15}}
\put(15,230){\line(1,0){15}}
\put(15,230){\line(-1,-1){15}}
\put(30,215){\line(0,1){15}}
\put(15,200){\circle*{3}}
\put(0,215){\circle*{3}}
\put(15,215){\circle*{3}}
\put(30,215){\circle*{3}}
\put(15,230){\circle*{3}}
\put(30,230){\circle*{3}}
\put(15,200){\line(0,1){30}}
\put(15,200){\line(1,2){15}}

\put(215,0){\line(1,1){15}}
\put(215,0){\line(-1,1){15}}
\put(215,30){\line(1,0){15}}
\put(215,30){\line(-1,-1){15}}
\put(230,15){\line(0,1){15}}
\put(215,0){\circle*{3}}
\put(200,15){\circle*{3}}
\put(215,15){\circle*{3}}
\put(230,15){\circle*{3}}
\put(215,30){\circle*{3}}
\put(230,30){\circle*{3}}
\put(200,15){\line(1,0){30}}
\put(215,0){\line(0,1){30}}
\put(215,15){\line(1,1){15}}

\put(95,40){\line(1,1){15}}
\put(95,40){\line(-1,1){15}}
\put(95,70){\line(1,0){15}}
\put(95,70){\line(-1,-1){15}}
\put(110,55){\line(0,1){15}}
\put(95,40){\circle*{3}}
\put(80,55){\circle*{3}}
\put(95,55){\circle*{3}}
\put(110,55){\circle*{3}}
\put(95,70){\circle*{3}}
\put(110,70){\circle*{3}}
\put(95,55){\line(-1,0){15}}
\put(95,55){\line(0,-1){15}}
\put(95,55){\line(1,1){15}}
\put(110,70){\line(-2,-1){30}}
\put(110,70){\line(-1,-2){15}}

\put(295,40){\line(1,1){15}}
\put(295,40){\line(-1,1){15}}
\put(295,70){\line(1,0){15}}
\put(295,70){\line(-1,-1){15}}
\put(310,55){\line(0,1){15}}
\put(295,40){\circle*{3}}
\put(280,55){\circle*{3}}
\put(295,55){\circle*{3}}
\put(310,55){\circle*{3}}
\put(295,70){\circle*{3}}
\put(310,70){\circle*{3}}
\put(295,55){\line(0,-1){15}}
\put(295,55){\line(1,1){15}}
\put(280,55){\line(1,0){30}}
\put(280,55){\line(2,1){30}}

\put(95,240){\line(1,1){15}}
\put(95,240){\line(-1,1){15}}
\put(95,270){\line(1,0){15}}
\put(95,270){\line(-1,-1){15}}
\put(110,255){\line(0,1){15}}
\put(95,240){\circle*{3}}
\put(80,255){\circle*{3}}
\put(95,255){\circle*{3}}
\put(110,255){\circle*{3}}
\put(95,270){\circle*{3}}
\put(110,270){\circle*{3}}
\put(110,270){\line(-1,-2){15}}
\put(110,270){\line(-2,-1){30}}

\put(295,240){\line(1,1){15}}
\put(295,240){\line(-1,1){15}}
\put(295,270){\line(1,0){15}}
\put(295,270){\line(-1,-1){15}}
\put(310,255){\line(0,1){15}}
\put(295,240){\circle*{3}}
\put(280,255){\circle*{3}}
\put(295,255){\circle*{3}}
\put(310,255){\circle*{3}}
\put(295,270){\circle*{3}}
\put(310,270){\circle*{3}}
\put(280,255){\line(1,0){30}}
\put(280,255){\line(2,1){30}}

\put(135,200){\line(1,1){15}}
\put(135,200){\line(-1,1){15}}
\put(135,230){\line(1,0){15}}
\put(135,230){\line(-1,-1){15}}
\put(150,215){\line(0,1){15}}
\put(135,200){\circle*{3}}
\put(120,215){\circle*{3}}
\put(135,215){\circle*{3}}
\put(150,215){\circle*{3}}
\put(135,230){\circle*{3}}
\put(150,230){\circle*{3}}
\put(135,230){\line(0,-1){30}}
\put(135,230){\line(1,-1){15}}

\put(215,120){\line(1,1){15}}
\put(215,120){\line(-1,1){15}}
\put(215,150){\line(1,0){15}}
\put(215,150){\line(-1,-1){15}}
\put(230,135){\line(0,1){15}}
\put(215,120){\circle*{3}}
\put(200,135){\circle*{3}}
\put(215,135){\circle*{3}}
\put(230,135){\circle*{3}}
\put(215,150){\circle*{3}}
\put(230,150){\circle*{3}}
\put(200,135){\line(1,0){30}}
\put(215,120){\line(0,1){30}}
\put(230,135){\line(-1,1){15}}

\put(255,220){\line(1,1){15}}
\put(255,220){\line(-1,1){15}}
\put(255,250){\line(1,0){15}}
\put(255,250){\line(-1,-1){15}}
\put(270,235){\line(0,1){15}}
\put(255,220){\circle*{3}}
\put(240,235){\circle*{3}}
\put(255,235){\circle*{3}}
\put(270,235){\circle*{3}}
\put(255,250){\circle*{3}}
\put(270,250){\circle*{3}}
\put(270,235){\line(-1,0){30}}
\put(270,235){\line(-1,1){15}}

\end{picture}
\end{center}
\caption{\label{fig:pentagon} The secondary polytope and all regular triangulations
for $\gA$ as in (\ref{eq:pentagon A}).}
\end{figure}
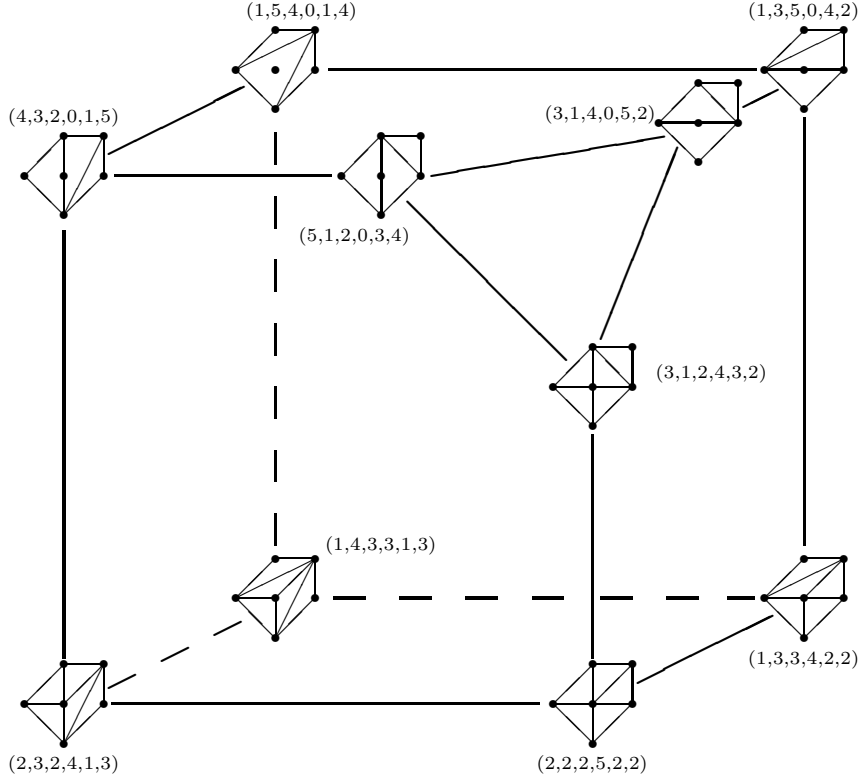

\

\noindent
\textbf{Example: } Figure \ref{fig:pentagon} shows the secondary polytope with at each vertex the corresponding regular triangulation of
$\Delta_\gA$
for $\gA$ as in (\ref{eq:pentagon A}).

\

\section{The toric variety associated with the Secondary Fan}
\label{section:gamma fan}

\subsection{Construction of the toric variety for the secondary fan}
\label{subsection:secondary toric}

The secondary fan is a complete fan of strongly convex polyhedral cones
in $\LL_\RR^\vee:=\Hom(\LL,\RR)$
which are generated by vectors from the lattice
$\LL_\ZZ^\vee:=\Hom(\LL,\ZZ)$.
By the general theory of toric varieties \cite{F} this lattice-fan  pair
gives rise to a toric variety. 
We are going to describe the construction of the toric variety
for the case of $\LL_\ZZ ^\vee$ and the secondary fan.
Before starting we must point out that \cite{F} works with a fan of \emph{closed cones}, while the cones in our Definition 
\ref{def:secondary fan} of the secondary fan are \emph{not closed};
see also Formula (\ref{eq:equivclass}). This difference, however, only affects a few minor
subtleties in the formulation at intermediate stages. The monoids
(\ref{eq:LC}) and therefore also the resulting toric varieties
are the same as in \cite{F}.

We denote the pairing between
$\LL_\RR:=\LL\otimes_\ZZ\RR$ and $\LL_\RR^\vee$ by $\langle,\rangle$.
For each cone $\cC$ in the secondary fan (see (\ref{eq:equivclass}))
one considers 
the affine scheme $\cU_\cC\::=\:{\spec\:}\ZZ[\LL_\cC]$ associated with
the monoid ring $\ZZ[\LL_\cC]$ of the monoid
\footnote{Alternative terminology: monoid = semi-group}
\begin{equation}\label{eq:LC}
\LL_\cC\::=\{\:\ull\in\LL\:|\: \langle\omega,\ull\rangle\geq 0\quad
\textrm{for all}\quad\omega\in\cC\:\}\,.
\end{equation}
In down-to-earth terms, a complex point of $\cU_\cC$ is just a homomorphism from the
additive monoid $\LL_\cC$ to the multiplicative monoid $\CC$,
sending $\mf{0}\in\LL_\cC$ to $1\in\CC$.

For cones $\cC$ and $\cC'$ in the secondary fan such that
$\cC'$ is contained in the closure $\overline{\cC}$ of $\cC$, there are inclusions
$$
\cC'\subset\overline{\cC}\,,\qquad
\LL_{\cC'}\supset\LL_{\cC}\,,\qquad \ZZ[\LL_{\cC'}]\supset\ZZ[\LL_{\cC}]
\,,\qquad \cU_{\cC'}\subset\cU_{\cC}\,;
$$
more precisely, the following lemma shows that the inclusion $\cU_{\cC'}\hookrightarrow\cU_{\cC}$ is an open immersion associated with
the inversion of an element in the ring $\ZZ[\LL_\cC]$.

\begin{lemma}\label{lemma:open subset}
In the above situation $\LL_{\cC'}=\LL_{\cC}\,+\,\ZZ\gl$
for some $\gl\in\LL_{\cC}$.
\end{lemma}
\begin{proof} If $\cC'=\cC$, the result is trivial. So assume
$\cC'\neq\cC$.
Since $\cC$ is a rational polyhedral cone it is spanned by
finitely many $\omega_1,\ldots,\omega_p\in\LL_\ZZ^\vee$, i.e.
every point in $\cC$ is a linear combination with non-negative real coefficients of $\omega_1,\ldots,\omega_p$. Moreover since 
$\cC'\subset\overline{\cC}$ there is a
$\gl\in\LL_{\cC}$ such that $\langle\omega',\gl\rangle= 0$ for all
$\omega'\in\cC'$ and $\langle\omega,\gl\rangle> 0$ for all
$\omega\in\cC$. Take $\mu\in\RR_{>0}$ such that for $j=1,\ldots,p$
we have $\langle\omega_j,\gl\rangle>\mu$
if $\omega_j\not\in\overline{\cC'}$.
For every $\ull\in\LL_{\cC'}$ and every non-negative integer 
$r>-\frac{1}{\mu}\min_{j}\langle\omega_j,\ull\rangle$, 
one now easily checks
that $\langle\omega,\ull+r\gl\rangle\geq 0$ for every $\omega\in\cC$, and hence $\ull+r\gl\in\LL_{\cC}$.
\end{proof}

\begin{definition}\label{def:toric variety}
The \emph{toric variety associated with the secondary fan} is the scheme
that results from glueing the affine schemes
$\cU_{\cC}$, where $\cC$ ranges over all cones in the secondary fan,
using the open immersions $\cU_{\cC'}\hookrightarrow\cU_{\cC}$
for $\cC'\subset\overline{\cC}$.
We denote this toric variety by $\cV_\gA$.
\end{definition}

For every cone $\cC$ of the secondary fan the monoid $\LL_\cC$
splits as a disjoint union $\LL_\cC\,=\,\LL_\cC^0\coprod\LL_\cC^+$
where $\LL_\cC^0$ (resp. $\LL_\cC^+$) is the set of elements  
which do (resp. do not) have an inverse in the
additive monoid $\LL_\cC$. One easily checks that
$$
\LL_\cC^0\::=\{\:\ull\in\LL\:|\: \langle\omega,\ull\rangle= 0\quad
\forall\omega\in\cC\:\}\,, \qquad
\LL_\cC^+\::=\{\:\ull\in\LL\:|\: \langle\omega,\ull\rangle> 0\quad
\forall\omega\in\cC\:\}\,.
$$ 
If $\cC=\{0\}$, then $\LL_\cC^0=\LL_\cC$ and $\LL_\cC^+=\emptyset$.
If $\cC\neq\{0\}$, the elements of $\LL_\cC^+$
generate a proper ideal $I_\cC$ in the ring $\ZZ[\LL_{\cC}]$ and
one has in $\cU_{\cC}$ the closed subscheme 
$$
\cB_\cC:={\spec}\:\modquot{\ZZ[\LL_\cC]}{I_\cC}\,.
$$

Let us see what this amounts to for the complex points of $\cU_{\cC}$,
viewed as homomorphisms from the additive monoid $\LL_\cC$ into the 
multiplicative monoid $\CC$. Each such homomorphism has to send
invertible elements to invertible elements, i.e. $\LL_\cC^0$ into
$\CC^*$. If $\cC=\{0\}$, the set of complex points of $\cU_{\cC}$
can therefore be identified with the 
set (in fact, $d$-dimensional torus) of group homomorphisms from $\LL$ into $\CC^*$:
\begin{equation}\label{eq:big cell}
\cU_{\{0\}}(\CC)\,=\,\Hom(\LL,\,\CC^*)\,=\,
\LL_\ZZ ^\vee\otimes_\ZZ \CC^*\,.
\end{equation}
If $\cC\neq\{0\}$, the complex points of $\cB_{\cC}$ are those monoid
homomorphisms that send all elements of $\LL_\cC^+$ to $\{0\}$.
So, the set of complex points of $\cB_{\cC}$ can be identified with the 
set of group homomorphisms from $\LL_\cC^0$ into
$\CC^*$:
$$
\cB_{\cC}(\CC)\,=\,\Hom(\LL_\cC^0,\,\CC^*)\,.
$$
This is a torus of dimension $d-\dim\cC$.

If $\cC$ is a maximal cone, then $\LL_\cC^0=\{0\}$ and $\cB_{\cC}(\CC)$ is only one point, which we denote as $\vp_\cC$.
For every positive real number $r<1$ the homomorphisms
$\LL_\cC\rightarrow\CC$ mapping $\LL_\cC^+$ into
the disc of radius $r$ centred at $0$ in $\CC$ form an open neighborhood of $\vp_\cC$, which we will also call the
disc of radius $r$ about $\vp_\cC$ in $\cU_\cC(\CC)$.

\

\noindent
\textbf{Example.}
For Gauss's hypergeometric structures $\LL=\ZZ(1,1,-1,-1)\subset\RR^4$.
So, $\LL_\RR^\vee\simeq\RR$ and 
the secondary fan has two maximal cones: $\RR_{>0}$ and $\RR_{<0}$.
One can easily see that the associated toric variety is the projective line $\PP^1$
(see \cite{F} p.6).

\

\noindent
\textbf{Example.}
For Appell's $F_4$ the secondary fan is shown in Figure \ref{fig:fan F4}.
One can easily see that the associated toric variety is the projective plane $\PP^2$ (see \cite{F} p.6-7).

\

\noindent
\textbf{Example.}
For Appell's $F_1$ the secondary fan is shown in Figure \ref{fig:fan F1}.
One can easily see that the associated toric variety is the projective plane $\PP^2$ with three points blown up (see \cite{F}).

\subsection{Convergence of Fourier $\gG$-series and the secondary fan}
\label{subsection:gamma fourier secondary fan}
We want to use the toric variety associated with the secondary fan 
to put the domains of convergence of the Fourier $\gG$-series 
(\ref{eq:Fourier Gamma series})
in the proper perspective.
Let us write $\hw\in\LL_\CC^\vee$ for the image of $\vw\in\CC^N$ under the natural projection $\CC^N\longrightarrow\LL_\CC^\vee:=\Hom(\LL,\CC)$ (linear dual; cf. (\ref{eq:dual exact sequence})). Then 
$\vw\cdot\ull=\langle\hw,\ull\rangle$ for all $\ull\in\LL$ and
(\ref{eq:Fourier Gamma series}) can be rewritten as
\begin{equation}\label{eq:FG}
\Psi_{\LL,\ugamma}(\vw)\;=\;e^{2\pi \ci\vw\cdot\ugamma}
\sum_{\ull\in\LL}\:
\frac{e^{2\pi \ci\langle\hw,\ull\rangle}}
{\prod_{j=1}^N\Gamma(\gamma_j+\ell_j+1)}\,.
\end{equation}
A vector $\hw\in\LL_\CC^\vee$ defines a homomorphism
$$
\LL\rightarrow\CC^\ast\,,\qquad\ull\mapsto 
e^{2\pi \ci \langle\hw,\ull\rangle}
$$
and, hence, a complex point of the toric variety $\cV_\gA$.
This point lies in the
disc of radius $r<1$ about the special point $\vp_\cC$ corresponding to a
maximal cone $\cC$ of the secondary fan if and only if 
$\langle\Im\hw,\ull\rangle>-\frac{\log r}{2\pi}$
for every non-zero $\ull\in\LL_\cC$; 
this means that $\Im\hw$ should lie `sufficiently far' inside the cone $\cC$. 

Recall that a maximal cone $\cC$ of the secondary fan corresponds to a
regular triangulation of the polytope $\Delta_\gA$. The index sets 
of the vertices of the maximal simplices in this triangulation constitute
a list $T_\cC$ of subsets of $\{1,\ldots,N\}$ with $N-d$ elements,
and according to (\ref{eq:equivclass}) 
\begin{equation}\label{eq:equivclass2}
\cC=\bigcap_{J\in T_\cC}
(\textrm{positive span of }\{\bb_j\}_{j\not\in J}).
\end{equation}
Now note that for $\ull=(\ell_1,\ldots,\ell_N)\in\LL$ almost tautologically $\ell_j=\langle\ull,\bb_j\rangle$. This shows
that for $\LL_J$ as defined in (\ref{eq:L cone2})
\begin{equation}\label{eq:LJC} 
\LL_J\subset\LL_\cC\qquad\textrm{for every}\quad J\in T_\cC\,.
\end{equation}

\

The above arguments together with those in Section
\ref{subsection:gamma fourier} show:

\begin{proposition}\label{prop:convergence+fan}
Let $\cC$ be a maximal cone of the secondary fan.
Let $J\in T_\cC$.
Let $\ugamma=(\gam_1,\ldots,\gam_N)\in\CC^N$ be such that
$\gam_j\in\ZZ_{\leq 0}$ for $j\not\in J$.
Then there is a positive real constant $r<1$ (depending on $\ugamma$)
such that the Fourier $\gG$-series $\Psi_{\LL,\ugamma}(\vw)$ in
(\ref{eq:FG}) converges for every $\vw\in\CC^N$ for which 
$\hw$ defines a point in the disc of radius $r$ about the special point $\vp_\cC$
in the toric variety $\cV_\gA$.
\qed
\end{proposition}

\subsection{Solutions of GKZ differential equations and the secondary fan}
\label{subsection:differential equations secondary fan}
Let us look for local solutions to the GKZ 
differential equations (\ref{eq:GKZ1})-(\ref{eq:GKZ2}) associated
with an $N$-element subset $\gA=\{\va_1,\ldots,\va_N\}\subset\ZZ^{k+1}$
and a vector $\vc\in\CC^{k+1}$. 
Let $\cC$ be a maximal cone in the secondary fan of $\gA$.
According to Proposition \ref{prop:convergence+fan},
every vector $\ugamma=(\gam_1,\ldots,\gam_N)\in\CC^N$
which satisfies
\begin{eqnarray}
\label{eq:c2gamma}
\gam_1\va_1+\ldots+\gam_N\va_N&=&\vc\,,
\\
\label{eq:C2gamma}
\exists\;J\in T_\cC\quad\textrm{such that}&&\gam_j\in\ZZ_{\leq 0}
\quad\textrm{for}\quad j\not\in J\,,
\end{eqnarray}
yields a Fourier $\gG$-series $\Psi_{\LL,\ugamma}(\vw)$ 
converging for every $\vw\in\CC^N$ for which 
$\hw$ defines a point in a sufficiently small disc about the point $\vp_\cC$ in $\cV_\gA$. According to Section \ref{subsection:gamma DE}
the corresponding $\gG$-series satisfies the GKZ 
differential equations (\ref{eq:GKZ1})-(\ref{eq:GKZ2}) for $\gA$
and $\vc$.  
If $\ugamma\equiv\ugamma'\bmod\LL$, then the two
Fourier $\gG$-series are equal. 
Lemma \ref{lemma:number Gamma series} will imply that the number of $\LL$-congruence classes of solutions to (\ref{eq:c2gamma})-(\ref{eq:C2gamma}) is finite  and, hence,
\emph{the Fourier $\gG$-series we obtain in this way have a common domain of convergence.}

\

\noindent
\textbf{Remark.}
Because of the factor $e^{2\pi \ci\vw\cdot\ugamma}$ in (\ref{eq:FG})
the Fourier $\gG$-series $\Psi_{\LL,\ugamma}(\vw)$ will in general
not descend to a function on some disc about $\vp_\cC$ in $\cV_\gA$.
On the other hand, if $\ugamma$ and $\ugamma'$ both satisfy
(\ref{eq:c2gamma})-(\ref{eq:C2gamma}), then $\ugamma-\ugamma'\in\LL_\CC$ and
$\vw\cdot(\ugamma-\ugamma')\,=\,\langle\hw,\ugamma-\ugamma'\rangle$
for every $\vw\in\CC^N$.
This means that the quotient 
$\Psi_{\LL,\ugamma}(\vw)\Psi_{\LL,\ugamma'}^{-1}(\vw)$
does descend to a function on some disc about $\vp_\cC$ in $\cV_\gA$.

\begin{lemma}\label{lemma:number Gamma series}
Fix $\vc\in\CC^{k+1}$ and a $k+1$-element set 
$J\subset\{1,\ldots,N\}$ such that the 
vectors $\va_j$ with $j\in J$ are linearly independent.
Then the number of classes modulo $\LL$ of vectors $\ugamma=(\gam_1,\ldots,\gam_N)\in\CC^N$ which satisfy
Equation (\ref{eq:c2gamma}) and
$\gam_j\in\ZZ$ for $j\in J':=\{1,\ldots,N\}\setminus J$,
is equal to $|\det\left((\va_j)_{j\in J}\right)|$.
\end{lemma}
\begin{proof} Since the 
vectors $\va_j$ with $j\in J$ are linearly independent,
the equation $\sum_{j=1}^N \gamma_j\va_j\,=\,\vc$ can be solved in
parametric form with the components $\gam_j$ for 
$j\in J'$ as free parameters. Every solution is the sum of one particular solution of the inhomogeneous system (e.g. the solution with
$\gam_j=0$ for $j\in J'$) and a solution of the homogeneous system.
So it suffices to determine the number of $\LL$-equivalence classes of
solutions of the equation $\sum_{j=1}^N \gamma_j\va_j\,=\,\mf{0}$
with $\gam_j\in\ZZ$ for $j\in J'$. The solutions themselves lie in
$\LL\otimes\QQ$.

Take any $d\times N$-matrix $\mf{B}$ whose rows form a $\ZZ$-basis of $\LL$. This amounts to choosing an isomorphism $\LL\simeq\ZZ^d$.
Let $\mf{B}_{J'}$ (resp. $\mf{B}_J$) 
denote the submatrix of $\mf{B}$ formed by the
columns with index in $J'$ (resp. in $J$). As in the proof of Lemma \ref{lemma:LJ}
one sees that the matrix $\mf{B}_{J'}$ is invertible over $\QQ$
and that the set of
solutions of $\sum_{j=1}^N \gamma_j\va_j\,=\,\mf{0}$
with $\gam_j\in\ZZ$ for $j\in J'$ is $\ZZ^d(\mf{B}_{J'})^{-1}\subset\QQ^d\simeq\LL\otimes\QQ$;
the notation $\ZZ^d(\mf{B}_{J'})^{-1}$ refers to the fact that here
$\ZZ^d$ consists of row vectors. 
The number of classes modulo $\LL$ of such solutions is therefore
$$
\sharp\left(\modquot{\ZZ^d(\mf{B}_{J'})^{-1}}{\ZZ^d}\right)\,=\,
\sharp\left(\modquot{\ZZ^d}{\ZZ^d\mf{B}_{J'}}\right)\,=\,
|\det(\mf{B}_{J'})|\,.
$$
Thus we must prove:
\begin{equation}\label{eq:detA=detB}
|\det(\mf{B}_{J'})|\,=\,|\det\left((\va_j)_{j\in J}\right)|\,.
\end{equation}
\textit{Proof of} (\ref{eq:detA=detB}):
Let $\mf{A}$ (resp. $\mf{A}_J$ resp. $\mf{A}_{J'}$)
denote the matrix with columns $\va_j$ with $j\in\{1,\ldots,N\}$ 
(resp. $j\in J$ resp. $j\in J'$). Then $\mf{A}\,\mf{B}^t=0$ and hence
\begin{equation}\label{eq:AB}
\mf{A}_J^{-1}\mf{A}_{J'}\,=\,-(\mf{B}_{J'}^{-1}\mf{B}_J)^t\,.
\end{equation}
As in Cramer's rule one sees that the matrix entries on the left hand side of (\ref{eq:AB}) are all of the form
$\pm(\det\mf{A}_J)^{-1}(\det\mf{A}_I)$ with $I\subset\{1,\ldots,N\}$
such that $\sharp I=k+1$ and $\sharp (I\cap J)=k$.
The corresponding matrix entries on the right hand side of (\ref{eq:AB}) are 
$\pm(\det\mf{B}_{J'})^{-1}(\det\mf{B}_{I'})$ with
$I':=\{1,\ldots,N\}\setminus I$. Thus we see
$$
|\det\mf{A}_J|^{-1}|\det\mf{A}_I|\,=\,
|\det\mf{B}_{J'}|^{-1}|\det\mf{B}_{I'}|\,,
$$
first for every $I\subset\{1,\ldots,N\}$
such that $\sharp I=k+1$ and $\sharp (I\cap J)=k$ and then, by induction,
for every $k+1$-element subset $I\subset\{1,\ldots,N\}$.
Consequently there are coprime positive integers $a,\,b$ such that
\begin{equation}\label{eq:AB2}
b|\det\mf{A}_I|\,=\,a|\det\mf{B}_{I'}|
\end{equation}
for every $k+1$-element subset $I\subset\{1,\ldots,N\}$.
Now recall that the columns $\va_1,\ldots,\va_N$ of $\mf{A}$
generate $\ZZ^{k+1}$. This implies that the greatest common divisor
of the numbers $\det\mf{A}_I$ is $1$. So $a=1$ in (\ref{eq:AB2}).
On the other hand, the rows of $\mf{B}$ form a $\ZZ$-basis of $\LL$.
Therefore for every prime number $p$
the rows the matrix $\mf{B}\bmod p\ZZ$ are linearly independent over the 
field $\modquot{\ZZ}{p\ZZ}$ and at least one of the numbers
$\det\mf{B}_{I'}$ must be not divisible by $p$. This shows
$b=1$ in (\ref{eq:AB2}) and finishes the proof of Formula (\ref{eq:detA=detB}).
\end{proof}

\begin{lemma}\label{lemma:resonance}
Let $\cC$ be a maximal cone of the secondary fan.
Let $\ugamma^{1},\ldots,\ugamma^{p}$ be solutions to the equations
(\ref{eq:c2gamma})-(\ref{eq:C2gamma}) such that
$\ugamma^{i}\not\equiv\ugamma^{j}\bmod \LL$ for $i\neq j$.
Then the Fourier $\gG$-series 
$\Psi_{\LL,\ugamma^{1}}(\vw),\ldots,\Psi_{\LL,\ugamma^{p}}(\vw)$
are linearly independent over $\CC$.
\end{lemma}
\begin{proof}
Fix a positive real constant $r<1$ such that all the given
Fourier $\gG$-series converge for every $\vw\in\CC^N$ for which 
$\hw$ defines a point in the disc of radius $r$ about the special point 
$\vp_\cC$ in the toric variety $\cV_\gA$ 
(cf. Proposition \ref{prop:convergence+fan}).
Next choose $\vw\in\CC^N$ such that $\hw$ defines a point in that disc
and such that no two of the numbers 
$\Im \vw\cdot (\ugamma^j+\ull)$ with $1\leq j\leq p$ and $\ull\in\LL_\cC$
such that the $\ull$-th term in the Fourier $\gG$-series 
$\Psi_{\LL,\ugamma^j}(\vw)$ is not $0$,
are equal. For this $\vw$ the set of those numbers 
$\Im \vw\cdot (\ugamma^j+\ull)$ assumes its minimum for a unique 
pair, say $(\ugamma^m,\,\ull^m)$.
This implies
$$
\lim_{t\to\infty}e^{-2\pi \ci t\vw\cdot(\ugamma^m+\ull^m)}
\Psi_{\LL,\ugamma^j}(t\vw)\,=\,0\quad\textrm{if}\quad j\neq m\,,
\qquad\textrm{resp.}\quad\neq\, 0\quad\textrm{if}\quad j=m\,.
$$
The linear independence claimed in the lemma now follows immediately.
\end{proof}

It follows from Lemma \ref{lemma:number Gamma series} that the number of $\LL$-congruence classes of solutions to the Equations (\ref{eq:c2gamma})-(\ref{eq:C2gamma}) is less than or equal to
$$
\sum_{J\in T_\cC}|\det\left((\va_j)_{j\in J}\right)|\,=\,
\mathrm{volume}\: \Delta_\gA\,.
$$

\begin{definition}\label{def:resonance}
Let $\cC$ be a maximal cone of the secondary fan of 
$\gA=\{\va_1,\ldots,\va_N\}$
and let $\vc\in\CC^{k+1}$. One says that $\vc$ is \emph{$\cC$-resonant} if
the number of $\LL$-congruence classes of solutions to the equations (\ref{eq:c2gamma})-(\ref{eq:C2gamma}) is less than
$\mathrm{volume}\: \Delta_\gA$.
\end{definition}

This means that $\vc$ is $\cC$-resonant if and only if there
is a $\ugamma=(\gam_1,\ldots,\gam_N)\in\CC^N$
which satisfies
$\gam_1\va_1+\ldots+\gam_N\va_N\,=\,\vc$,
and for which there are two different sets $J_1$ and $J_2$ on the
list $T_\cC$ such that
$\gam_j\in\ZZ$ for $j\in \{1,\ldots,N\}\setminus (J_1\cap J_2)$.

\begin{corollary}\label{cor:nonresonant}
If $\vc$ is not $\cC$-resonant, the Fourier $\gG$-series
$\Psi_{\LL,\ugamma}(\vw)$ associated with solutions $\ugamma$
of the equations (\ref{eq:c2gamma})-(\ref{eq:C2gamma}) are linearly independent and span
a space of local solutions of the GKZ 
differential equations (\ref{eq:GKZ1})-(\ref{eq:GKZ2}) of
dimension equal to $\mathrm{volume}\: \Delta_\gA$.
According to the discussion in Section 
\ref{subsection:dimension solution space} this is then the full space of local solutions if (for instance)
the polytope $\Delta_\gA$ admits a unimodular triangulation.
\qed
\end{corollary}

\section{Extreme resonance in GKZ systems}
\label{section:resonance}
In this section $\cC$ is a maximal cone of the secondary fan of 
$\gA=\{\va_1,\ldots,\va_N\}$ for which the corresponding regular triangulation of $\Delta_\gA$
is unimodular, i.e. 
$$
|\det\left((\ba_j)_{j\in J}\right)|=1\qquad\textrm{for every}\quad
J\in T_\cC.
$$
This means that for every $J\in T_\cC$ the set 
$\{\ba_j\}_{j\in J}$ is a $\ZZ$-basis of $\ZZ^{k+1}$.
Consequently, for every $\vc\in\ZZ^{k+1}$ all solutions 
$\ugamma=(\gam_1,\ldots,\gam_N)$ of the equations
(\ref{eq:c2gamma})-(\ref{eq:C2gamma}) lie in $\ZZ^N$ and are
therefore congruent modulo $\LL$.
So a vector $\vc\in\ZZ^{k+1}$ is $\cC$-resonant, in an extreme way:
all Fourier $\gG$-series coming from solutions of 
(\ref{eq:c2gamma})-(\ref{eq:C2gamma}) are equal!
\emph{In this section we will demonstrate how one can obtain, locally
near the point $\vp_\cC$ on $\cV_\gA$, more
solutions of the GKZ differential equations 
(\ref{eq:GKZ1})-(\ref{eq:GKZ2}) from an
`infinitesimal deformation' of this Fourier $\gG$-series.}

\begin{definition}\label{def:ring R}
For $\gA$ and $\cC$ as above we define the ring 
\begin{equation}\label{eq:ring R}
\cR_{\gA,\cC}:=\modquot{\ZZ[E_1,\ldots,E_N]}{(\cI_\gA\,+\,\cI_\cC)}
\end{equation}
where $\ZZ[E_1,\ldots,E_N]$ is just the polynomial ring over $\ZZ$
in $N$ variables,\\ 
$\cI_\gA$ is the ideal generated by the
linear forms which are the components of the vector
\begin{equation}\label{eq:ideal A}
E_1\va_1+\ldots +E_N\va_N\,,
\end{equation}
and $\cI_\cC$ is the ideal generated by the monomials
\begin{equation}\label{eq:ideal C}
E_{i_1}\cdots E_{i_s}\qquad\textrm{with} \quad
\{i_1,\ldots,i_s\}\not\subset J\quad\textrm{for all}\quad J\in T_\cC\,.
\end{equation}
We write $\eps_j$ for the image of $E_j$ in $\cR_{\gA,\cC}$.
\end{definition}

\

\noindent
So in $\cR_{\gA,\cC}$ we have the relations
\begin{eqnarray}
\label{eq:linear relations}
\eps_1\va_1+\ldots \eps_N\va_N&=&0\,,
\\
\label{eq:monomial relations}
\eps_{i_1}\cdots \eps_{i_s}&=&0\qquad\textrm{if} \quad
\{i_1,\ldots,i_s\}\not\subset J\quad\textrm{for all}\quad J\in T_\cC\,.
\end{eqnarray}
Relation (\ref{eq:linear relations}) means that the vector $\ueps=(\eps_1,\ldots,\eps_N)\in\cR_{\gA,\cC}^N$ lies in
$\LL\otimes_\ZZ \cR_{\gA,\cC}$.

\

\noindent
\textbf{Remark.} The ideal $\cI_\cC$ is well-known in combinatorial 
algebra \cite{sta}, where it is called the \emph{Stanley-Reisner ideal}
of the triangulation $T_\cC$.
The ring $\modquot{\ZZ[E_1,\ldots,E_N]}{\cI_\cC}$
is called the \emph{Stanley-Reisner ring}.

\

The following facts about the ring $\cR_{\gA,\cC}$ are proven in
\cite{sti} \S 2.

\begin{proposition}\label{prop:ring R}
\begin{enumerate}
\item
$\cR_{\gA,\cC}$ is a free $\ZZ$-module of rank equal to 
$\mathrm{volume}\: \Delta_\gA$.
\item
$\cR_{\gA,\cC}$ is a graded ring and each $\eps_j$ has degree $1$.
\item
Denoting the homogeneous part of degree $i$ in $\cR_{\gA,\cC}$ 
by $\cR_{\gA,\cC}^{(i)}$ one has isomorphisms (see also
\S\ref{subsection:secondary fan construction})
\begin{equation}\label{eq:R01}
\cR_{\gA,\cC}^{(0)}=\ZZ\,,\qquad
\cR_{\gA,\cC}^{(1)}\simeq \LL_\ZZ^\vee\,,\quad \eps_j\mapsto \vb_j\,.
\end{equation}
\item
The Poincar\'e series of the graded ring $\cR_{\gA,\cC}$ is
\begin{equation}\label{eq:poincare}
\sum_{i\geq 0} \left(\rank\,\cR_{\gA,\cC}^{(i)}\right) T^i\,=\,
\sum_{m=0}^{k+1} S_{\cC,m}\,T^m(1-T)^{k+1-m}\,,
\end{equation}
where $S_{\cC,0}=1$ and $S_{\cC,m}$, for $m\geq 1$, is the number of simplices with $m$ vertices in the triangulation of $\Delta_\gA$
corresponding with $\cC$.
In particular $\cR_{\gA,\cC}^{(i)}=0$ for $i\geq k+1$
and the elements $\eps_1,\ldots,\eps_N$ are nilpotent.
\end{enumerate}
\qed
\end{proposition}

\

\noindent
For the examples at the end of Sections 
\ref{subsection:Secondary Polytope} and 
\ref{subsection:secondary fan construction}
(see also Figures \ref{fig:classical examples}, 
\ref{fig:fan F1}, \ref{fig:fan F4})
we find:

\

\noindent
\textbf{Example.}
For $\LL=\ZZ(-2,1,1)\subset\ZZ^3$ there is only one
unimodular triangulation, namely $T_\cC=\{\{1,2\},\{1,3\}\}$.
One easily checks that in this case
$$
\cR_{\gA,\cC}\,=\,\ZZ\oplus\ZZ\eps\,,\qquad \eps^2=0\,,\qquad
(\eps_1,\eps_2,\eps_3)=(-2\eps,\eps,\eps)\,.
$$

\

\noindent
\textbf{Example.}
For Gauss $\LL=\ZZ(1,1,-1,-1)\subset\ZZ^4$ and there are two
unimodular triangulations, which both lead to
$$
\cR_{\gA,\cC}\,=\,\ZZ\oplus\ZZ\eps\,,\qquad \eps^2=0\,.
$$
For one triangulation
$(\eps_1,\eps_2,\eps_3,\eps_4)$ is $(-\eps,-\eps,\eps,\eps)$,
for the other $(\eps,\eps,-\eps,-\eps)$.

\

\noindent
\textbf{Example: }
For $\LL=\ZZ(-3\,1\,1\,1)\subset\ZZ^4$  there is only one
unimodular triangulation, namely 
$T_\cC=\{\;\{1,3,4\}\,,\:\{1,2,4\}\,,\:\{1,2,3\}\;\}$.
One easily checks that in this case
$$
\cR_{\gA,\cC}\,=\,\ZZ\oplus\ZZ\eps\oplus\ZZ\eps^2
\,,\qquad \eps^3=0\,,\qquad
(\eps_1,\eps_2,\eps_3,\eps_4)=(-3\eps,\eps,\eps,\eps)\,.
$$

\

\noindent
\textbf{Example.}
For Appell's $F_1$ $\LL=\ZZ(1,-1,0,-1,1,0)\oplus\ZZ(1,0,-1,-1,0,1)$.
There are six unimodular triangulations 
(see Figure \ref{fig:fan F1}).
One can check that for the triangulation 
$T_\cC=\{\{3,4,5,6\},\, \{1,2,3,4\},\, \{2,3,4,5\}\}$
the relations (\ref{eq:ideal A})-(\ref{eq:ideal C}) yield
\begin{eqnarray*}
&&
(\eps_1,\eps_2,\eps_3,\eps_4,\eps_5,\eps_6)=
\eps(1,-1,0,-1,1,0)+\delta(1,0,-1,-1,0,1)\,,
\\
&&
\eps_1\eps_5\,=\,\eps_1\eps_6\,=\,\eps_2\eps_6\,=\,0\,,
\end{eqnarray*}
and hence:
$$
\cR_{\gA,\cC}\,=\,\ZZ\oplus\ZZ\eps\oplus\ZZ\delta
\,,\qquad \eps^2\,=\,\delta^2\,=\,\eps\delta\,=\,0\,.
$$

\

\noindent
\textbf{Example.}
For Appell's $F_4$ $\LL=\ZZ(1,-1,1,-1,0,0)\oplus\ZZ(1,0,1,0,0,-1,-1)$.
There are three unimodular triangulations
(see Figure \ref{fig:fan F4}).
One can check that for the triangulation
$T_\cC=\{\{1,3,4,6\},\, \{1,3,4,5\},\, \{1,2,3,5\},\, \{1,2,3,6\}\}$
the relations (\ref{eq:ideal A})-(\ref{eq:ideal C}) yield
\begin{eqnarray*}
&&
(\eps_1,\eps_2,\eps_3,\eps_4,\eps_5,\eps_6)=
\eps(1,-1,1,-1,0,0)+\delta(1,0,1,0,0,-1,-1)\,,
\\
&&
\eps_2\eps_4\,=\,\eps_5\eps_6\,=\,0\,,
\end{eqnarray*}
and hence:
$$
\cR_{\gA,\cC}\,=\,\ZZ\oplus\ZZ\eps\oplus\ZZ\delta\oplus\ZZ\eps\delta
\,,\qquad \eps^2\,=\delta^2\,=\,\,0\,.
$$

\

\

For $z\in\CC$ and nilpotent $\eps$ one can define
$\frac{1}{\gG(z+\eps)}$ as an element of $\CC[\eps]$ by using the 
Taylor expansion of the function $\frac{1}{\gG}$ at $z$:
$$
\frac{1}{\gG(z+\eps)}\,:=\,\frac{1}{\gG(z)}\,+\,\eps
\left(\frac{1}{\gG}\right)'(z)\,+\,\frac{\eps^2}{2}
\left(\frac{1}{\gG}\right)''(z)\,+\,\frac{\eps^3}{3!}
\left(\frac{1}{\gG}\right)'''(z)\,+\,\ldots\,.
$$
One defines similarly $\gG(1+\eps)$. Thus for $z\in\CC$ and 
nilpotent $\eps$ also $\frac{\gG(1+\eps)}{\gG(z+1+\eps)}$ has been
defined. From (\ref{eq:Pochhammer}) one sees that for $m\in\ZZ$:
\begin{equation}\label{eq:gammadeform}
\frac{\gG(1+\eps)}{\gG(m+1+\eps)}\:=\,
\left\{\begin{array}{lll}
\displaystyle{\frac{1}{(1+\eps)(2+\eps)\cdots(m+\eps)}}
&\textrm{ if }&m>0\\[2ex]
1&\textrm{ if }&m=0\\
\eps(\eps-1)(\eps-2)\cdots(\eps+m+1)&\textrm{ if }&m<0
\end{array}\right.
\end{equation}
Finally, for $z\in\CC$, $u\in\CC^*$ (with a choice of a branch of 
$\log u$) and nilpotent $\eps$ one has naturally
$$
e^{\eps\,z}\,:=\,\sum_{m\geq 0}\frac{1}{m!}\eps^m\,z^m\,,
\qquad
u^\eps\,:=\,e^{\eps\,\log u}\,.
$$
We are ready to present our deformation of the (Fourier) $\gG$-series:

\begin{definition}\label{def:fourier gamma deform}
For $\ugamma=(\gam_1,\ldots,\gam_N)\in\ZZ^N$
and $\ueps=(\eps_1,\ldots,\eps_N)\in\cR_{\gA,\cC}^N$ we define
\begin{eqnarray}\label{eq:Fourier Gamma series deform}
\Psi_{\LL,\ugamma,\ueps}(\vw)&:=&
\sum_{\ull\in\LL}\:\prod_{j=1}^N
\frac{\Gamma(1+\eps_j)}{\Gamma(\gamma_j+\ell_j+1+\eps_j)}\;
e^{2\pi \ci\vw\cdot(\ugamma+\ull+\ueps)}\,,
\\
\label{eq:Gamma series deform}
\Phi_{\LL,\ugamma,\ueps}(\vu)&:=&
\sum_{\ull\in\LL}\:\prod_{j=1}^N
\frac{\Gamma(1+\eps_j)\:u_j^{\gamma_j+\ell_j+\eps_j}}
{\Gamma(\gamma_j+\ell_j+1+\eps_j)}\,.
\end{eqnarray}
\end{definition}

\

\noindent
\textbf{Remark.} From the point of view of deforming $\ugamma$
it seems more natural to consider
\begin{eqnarray}\label{eq:naive deform 1}
\Psi_{\LL,\ugamma+\ueps}(\vw)&:=&
\sum_{\ull\in\LL}\:\prod_{j=1}^N
\frac{1}{\Gamma(\gamma_j+\ell_j+1+\eps_j)}\;
e^{2\pi \ci\vw\cdot(\ugamma+\ull+\ueps)}\,,
\\
\label{eq:naive deform 2}
\Phi_{\LL,\ugamma+\ueps}(\vu)&:=&
\sum_{\ull\in\LL}\:\prod_{j=1}^N
\frac{u_j^{\gamma_j+\ell_j+\eps_j}}
{\Gamma(\gamma_j+\ell_j+1+\eps_j)}\,;
\end{eqnarray}
i.e.
$$
\Psi_{\LL,\ugamma+\ueps}(\vw)=
\frac{\Psi_{\LL,\ugamma,\ueps}(\vw)}{\prod_{j=1}^N\Gamma(1+\eps_j)}
\,,\qquad
\Phi_{\LL,\ugamma+\ueps}(\vu)=
\frac{\Phi_{\LL,\ugamma,\ueps}(\vu)}{\prod_{j=1}^N\Gamma(1+\eps_j)}\,.
$$
Indeed, expanding these functions in coordinates
with respect to a basis of $\cR_{\gA,\cC}$ 
is for (\ref{eq:naive deform 1}) and (\ref{eq:naive deform 2})
essentially just Taylor expansion, if one views the expressions as
functions of $\ugamma$, while the interpretation 
as (multi-valued) local solutions of GKZ differential equations
with values in $\cR_{\gA,\cC}\otimes\CC$ (see below)
are equally true for 
(\ref{eq:naive deform 1})-(\ref{eq:naive deform 2}) in place of 
(\ref{eq:Fourier Gamma series deform})-(\ref{eq:Gamma series deform}).
We prefer, however,
the latter because
their coordinates are series with rational coefficients,
whereas the coefficients of the coordinate series of the former
involve interesting, but mysterious non-rational numbers
like the Euler-Masceroni constant and values of Riemann's zeta-function.
We can be slightly more informative about the coefficients in
(\ref{eq:naive deform 1})-(\ref{eq:naive deform 2}):
there is the well-known formula for the $\gG$-function due to Gauss
$$
\gG(s)\,=\,\lim_{n\to\infty}\left[
\frac{n!\,n^s}{s(s+1)\cdots (s+n)}\right]\,,
$$
from which one easily derives the expansion
$$
\log\gG(1+s)\,=\,-\eul s\,+\,
\sum_{m=2}^\infty (-1)^m\zeta (m)\frac{s^m}{m}
$$
where $\eul$ denotes the Euler-Masceroni constant
and $\zeta$ is Riemann's zeta-function. By exponentiating and re-expanding
one finds the Taylor expansion for $\gG(1+s)$
and then eventually the expansion of 
$\left[\prod_{j=1}^N\Gamma(1+\eps_j)\right]^{-1}$.

\

\begin{lemma}\label{lemma:fourier gamma deform}
There are finitely many $\ull^{(1)},\ldots,\ull^{(r)}\in\LL_\cC$ (with $\LL_\cC$ as in (\ref{eq:LC}))
such that the series $\Psi_{\LL,\ugamma,\ueps}(\vw)$ and
$\Phi_{\LL,\ugamma,\ueps}(\vu)$ involve only terms with
$$
\ull\in\bigcup_{i=1}^r\;(-\ull^{(i)}\,+\,\LL_\cC)\,.
$$ 
In particular for $\ugamma=\mf{0}$ the series involve only terms
with $\ull\in\LL_\cC$.
\end{lemma}
\begin{proof}
It follows immediately from (\ref{eq:gammadeform}) and
(\ref{eq:ideal C}) that for the terms which appear with non-zero coefficient, the set $\{j\:|\:\gamma_j+\ell_j\,<0\}$ is contained in 
some $J$ on the list $T_\cC$. Suppose 
$\{j\:|\:\gamma_j+\ell_j\,<0\}\subset J\in T_\cC$. Then $\gamma_j+\ell_j\geq 0$
for every $j\in J':=\{1,\ldots,N\}\setminus J$.
The vector $\sum_{j\in J'}\max(0,\gamma_j)\,\va_j$ is a
$\ZZ$-linear combination of the vectors $\va_i$ with 
$i\in J$, because the triangulation is unimodular. Such a relation is an element of $\LL$. Thus one sees that
$\LL$ contains an element 
$\ull^J=(\ell_1^J,\ldots,\ell_N^J)$ with
$\ell^J_j=\max(0,\gamma_j)$ for all $j\in J'$.
So $\ell^J_j+\ell_j\geq 0$ for every $j\in J'$. In the notation
introduced in (\ref{eq:L cone2}) this can be written as $\ull^J\in\LL_J$
and
$\ull^J+\ull\in\LL_J$. The lemma now follows from (\ref{eq:LJC}).
\end{proof}

\

Partial sums (with finitely many  terms) of the series 
(\ref{eq:Fourier Gamma series deform}) 
resp. (\ref{eq:Gamma series deform})
can be evaluated as elements in the ring $\cR_{\gA,\cC}\otimes\CC$
and be written in coordinates with respect to a $\ZZ$-basis of the finite rank  $\ZZ$-module $\cR_{\gA,\cC}$.
These coordinates are again partial sums of series.
In \cite{sti} \S3 one finds estimates on the growth of the coefficients of these series and on a common domain of convergence.
Thus $\Psi_{\LL,\ugamma,\ueps}(\vw)$ and $\Phi_{\LL,\ugamma,\ueps}(\vu)$
are functions with values in $\cR_{\gA,\cC}\otimes\CC$.
The function $\Psi_{\LL,\ugamma,\ueps}(\vw)$ is defined for
$\vw\in\CC^N$ with $\Im\hw$ `sufficiently far' inside the cone
$\cC$ (cf. \S\ref{subsection:gamma fourier secondary fan}).
Because of the appearance of logarithms $\Phi_{\LL,\ugamma,\ueps}(\vu)$ 
is actually a multi-valued function, defined on some open disc
about $\mf{0}$ in $\CC^\gA$ with the divisor $u_1\cdots u_N=0$ removed. The multi-valuedness is easily described using the relation 
$u_j=e^{2\pi \ci w_j}$ which matches 
$w_j$ with a choice of $\log u_j$. A different choice adds an integer to
$w_j$. Now note that for $\mf{m}\in\ZZ^N$
$$
\Psi_{\LL,\ugamma,\ueps}(\vw+\mf{m})\,=\,
e^{2\pi \ci\mf{m}\cdot\ueps}\:\Psi_{\LL,\ugamma,\ueps}(\vw)\,.
$$
This formula can also be read as a precise expression for local monodromy.
Since
$\{\mf{m}\cdot\ueps\;|\;\mf{m}\in\ZZ^N\}\,=\,\cR_{\gA,\cC}^{(1)}$,
we can summarize our analysis of the multi-valuedness of
$\Phi_{\LL,\ugamma,\ueps}(\vu)$ as follows:

\begin{proposition}\label{prop:Phi deform}
$\Phi_{\LL,\ugamma,\ueps}(\vu)$ is a multi-valued function
with values in $\cR_{\gA,\cC}\otimes\CC$. Different branches
of this function are related by multiplication with an element
$e^{2\pi \ci\omega}$ with $\omega\in\cR_{\gA,\cC}^{(1)}$.
\qed
\end{proposition}

\

The same arguments as those used in Section \ref{subsection:gamma DE}
show immediately

\begin{proposition}\label{prop:Phi deform DE}
The $\cR_{\gA,\cC}\otimes\CC$-valued function $\Phi_{\LL,\ugamma,\ueps}(\vu)$ 
satisfies the GKZ system of differential equations (\ref{eq:GKZ1})-(\ref{eq:GKZ2})
for $\gA$ and $\vc=\sum_{j=1}^N \gamma_j\va_j$.

The $\CC$-valued functions which arise as coordinates of
$\Phi_{\LL,\ugamma,\ueps}(\vu)$ with respect to a basis of
$\cR_{\gA,\cC}$ satisfy the same GKZ system of differential equations.
\qed
\end{proposition}

\

\noindent
\textbf{Example.}
For $\LL=\ZZ(-3,1,1,1)$ and $\cC$ the cone corresponding to the unimodular triangulation
$T_\cC=\{\{1,2,3\},\{1,2,4\},\{1,3,4\}\}$ 
$$
\cR_{\gA,\cC}\,=\,\ZZ\oplus\ZZ\eps\oplus\ZZ\eps^2
\,,\qquad \eps^3=0\,,\qquad
(\eps_1,\eps_2,\eps_3,\eps_4)=(-3\eps,\eps,\eps,\eps)\,.
$$
For $\uz=(0,0,0,0)$ one then finds, using (\ref{eq:gammadeform}) and setting $z=u_1^{-3}u_2u_3u_4$,
\begin{eqnarray*}
\Phi_{\LL,\uz,\ueps}(\vu)&=& 
\sum_{m\in\ZZ}\:
\frac{\Gamma(1-3\eps)}
{\Gamma(1-3m-3\eps)}
\left(\frac{\Gamma(1+\eps)}
{\Gamma(1+m+\eps)}\right)^3
u_1^{-3m-3\eps}
u_2^{m+\eps}
u_3^{m+\eps}
u_4^{m+\eps}
\\
&=&z^{\eps}\left(1\,+\,\sum_{m\geq 1}
\frac{(-3\eps)(-3\eps-1)\cdots(-3\eps-3m+1)}
{((1+\eps)\cdots (m+\eps))^3}\:z^m\right)
\\
&=&\left(\left.\left.1\,+\,\eps\log z\,+\,\frac{\eps^2}{2}\log^2 z
\right)\right(1\,+\,\eps G_1(z)\,+\,\eps^2 G_2(z)\right)
\\
&=& 1\,+\,(\log z\,+\,G_1(z))\eps
\,+\,
(\textstyle{\frac{1}{2}}\log^2 z\,+\,G_1(z)\log z\,+\,G_2(z))\eps^2
\end{eqnarray*}
with
\begin{eqnarray*}
G_1(z)&=&3\sum_{m\geq 1}(-1)^m\frac{(3m-1)!}{(m!)^3}\:z^m
\\
G_2(z)&=&9\sum_{m\geq 1}(-1)^m\frac{(3m-1)!}{(m!)^3}
\left(\sum_{j=m+1}^{3m-1}\frac{1}{j}\right)\:z^m\,.
\end{eqnarray*}
Similarly, for $\ugamma=(-1,0,0,0)$ we obtain
\begin{eqnarray*}
\Phi_{\LL,\ugamma,\ueps}(\vu)&=& 
\sum_{m\in\ZZ}\:
\frac{\Gamma(1-3\eps)}
{\Gamma(-3m-3\eps)}
\left(\frac{\Gamma(1+\eps)}
{\Gamma(m+1+\eps)}\right)^3
u_1^{-1-3m-3\eps}
u_2^{m+\eps}
u_3^{m+\eps}
u_4^{m+\eps}
\\
&=&
u_1^{-1}\sum_{m\geq 0}\:
\frac{(-3\eps)(-3\eps-1)\cdots(-3\eps-3m)}
{((1+\eps)\cdots (m+\eps))^3}\:
(u_1^{-3}u_2u_3u_4)^{m+\eps}
\\
&=&
u_1^{-1}\left(\left.\left.1\,+\,\eps\log z\,+\,\frac{\eps^2}{2}\log^2 z
\right)\right(\eps F_1(z)\,+\,\eps^2 F_2(z)\right)
\\
&=& u_1^{-1}F_1(z)\eps\,+\,u_1^{-1}(F_1(z)\log z\,+\,F_2(z))\eps^2
\end{eqnarray*}
with
\begin{eqnarray*}
F_1(z)&=&-3\sum_{m\geq 0}\:(-1)^m\frac{(3m)^!}{(m!)^3} z^m\,,
\\
F_2(z)&=&-9\sum_{m\geq 1}\:(-1)^m\frac{(3m)^!}{(m!)^3}
\left(\sum_{j=m+1}^{3m}\frac{1}{j}\right) z^m\,.
\end{eqnarray*}
Note that in agreement with Proposition \ref{prop:GKZ derivatives}
$$
\Phi_{\LL,\ugamma,\ueps}(\vu)\,=\,
\frac{\partial}{\partial u_1}\Phi_{\LL,\uz,\ueps}(\vu)\,=\,
-3u_1^{-1}z\frac{\partial}{\partial z}\Phi_{\LL,\uz,\ueps}(\vu)\,.
$$
The components of $\Phi_{\LL,\uz,\ueps}(\vu)$
are three linearly independent solutions of the GKZ system of differential equations with $\vc=\mf{0}$, whereas the components of $\Phi_{\LL,\ugamma,\ueps}(\vu)$ yield only two linearly independent solutions of the GKZ system for $\vc=-\va_1$.
Since in this case $\mathrm{volume}\: \Delta_\gA\,=\,3$ we find 
enough solutions for $\vc=\mf{0}$, but
not enough for $\vc=-\va_1$ (see Section
\ref{subsection:dimension solution space}).

\

The phenomenon observed at the end of the previous example --
namely that our method yields enough solutions if $\vc=\mf{0}$,
but misses solutions if $\vc\neq\mf{0}$ -- occurs quite generally.
Below, in Theorem \ref{thm:solution basis}, we quote \cite{sti} Theorem 5
and also recall some conclusions (e.g. Proposition \ref{prop:Phi deform})
found earlier in the present notes:

\begin{theorem}\label{thm:solution basis}
Let $\cC$ be a maximal cone of the secondary fan of 
$\gA$ for which the corresponding regular triangulation of $\Delta_\gA$
is unimodular. Let $\uz=(0,\ldots,0)$.
Then the coordinates of the
$\cR_{\gA,\cC}\otimes\CC$-valued function $\Phi_{\LL,\uz,\ueps}(\vu)$
with respect to a basis of the free $\ZZ$-module $\cR_{\gA,\cC}$
constitute a basis for the local solution space of 
the GKZ system of differential equations (\ref{eq:GKZ1})-(\ref{eq:GKZ2})
for $\gA$ and $\vc=\mf{0}$.
These multi-valued functions are invariant under the action 
(\ref{eq:torus action on functions}) of the torus $\TT^{k+1}$
and descend therefore to multi-valued functions
on a disc minus a divisor centered at the point $\vp_\cC$ in the toric
variety $\cV_\gA$.
The multi-valuedness of these functions
is given by multiplying $\Phi_{\LL,\uz,\ueps}(\vu)$
with elements in the group 
$\{e^{2\pi \ci\omega}\:|\:\omega\in\cR_{\gA,\cC}^{(1)}\}$.
\qed
\end{theorem}

\

\noindent
\textbf{Remark.}
In \cite{boo} Anne de Boo carefully re-examined the preceding method
and improved it by also taking $\ugamma$ into account.
In this way he obtained full local solution spaces for 
GKZ systems of differential equations 
for many more instances of the triangulation of $\Delta_\gA$
and of the parameter $\vc$.

Very recently Borisov and Horja \cite{BH} found a way to obtain enough solutions for any $\vc\in\ZZ^{k+1}$ and any triangulation.
Their method is close in spirit to the method in this Section \ref{section:resonance}. We recommend \cite{BH} for further reading on this aspect of GKZ hypergeometric structures.

\section{GKZ for Lauricella's $F_D$}
\label{section:FD}
Since Lauricella's $F_D$ also plays an important role in other lectures in this School, we put details of the GKZ theory for Lauricella's $F_D$ together in this section.
\subsection{Series, $\LL$, $\gA$ and the primary polytope $\Delta_\gA$}
\label{subsection:FDLA}

Recall that in Section \ref{subsubsection:Appell-Lauricella}
we found, starting from the power series expansion of  
Lauricella's $F_D$ in $k-1$ variables
$$
F_D(a,\vb,c |\vz):=\sum_{\vm}
\frac{(a)_{|\vm |}(\vb)_{\vm}}{(c)_{|\vm |} \vm !}\,\vz^{\vm}\,,
$$ 
that the lattice $\LL$ is generated by the rows of the following
$(k-1)\times (2k)$-matrix
\begin{equation}\label{eq:FD L matrix}
\left(
\begin{array}{rrrrrrrrrr}
1&-1&0&\ldots&0&-1&1&0&\ldots&0\\
1&0&-1&\ddots&\vdots&-1&0&1&\ddots&\vdots\\
\vdots&\vdots&\ddots&\ddots&0&\vdots&\vdots&\ddots&
\ddots&0\\
1&0&\ldots&0&-1&-1&0&\ldots&0&1
\end{array}
\right)\,.
\end{equation}
So for $\gA$ we can take the set of columns of the
$(k+1)\times (2k)$-matrix
\begin{equation}\label{eq:FD A matrix}
\left(
\begin{array}{rrrrrrrrrr}
0&0&0&\ldots&0&1&1&1&\ldots&1\\
1&0&0&\ldots&0&1&0&0&\ldots&0\\
0&1&0&\ldots&0&0&1&0&\ldots&0\\
0&0&1&\ddots&\vdots&0&0&1&\ddots&\vdots\\
\vdots&\vdots&\ddots&\ddots&0&\vdots&\vdots&\ddots&
\ddots&0\\
0&0&\ldots&0&1&0&0&\ldots&0&1
\end{array}
\right)\,.
\end{equation}
This notation is consistent with the main part of this text:
$\gA$ is a subset of $\ZZ^{k+1}$; moreover
$N\,=\,2k$ and $d=\rank\LL=k-1$.

The primary polytope $\Delta_\gA$ is the direct product of a $(k-1)$-simplex and a $1$-simplex and, for $k=3$, looks like 
the prism in Figure \ref{fig:classical examples}. 
The vectors $\va_1,\ldots,\va_k$
are in the bottom face of the prism; $\va_{k+1},\ldots,\va_{2k}$
are in the top face. The numbering is such that the difference vectors
$\va_{k+j}-\va_j$, for $j=1,\ldots,k$ are all equal.

\subsection{Integrals and differential equations for $F_D$}
\label{subsection:IntFD}
In \cite{loo} Lauricella's $F_D$ in variables
$z_0,\ldots,z_n$ is introduced via the integrals
\begin{equation}\label{eq:integral FD}
F_\ga(z_0,\ldots,z_n)\,:=\,
\int_{\ga}(z_0-\zeta)^{-\mu_0}\cdots
(z_n-\zeta)^{-\mu_n}\,d\zeta
\end{equation}
over suitable intervals $\ga$, with endpoints in $\{z_0,\ldots,z_n,\infty\}$. Note that because of the translation invariance 
property 
\begin{equation}\label{eq:translation FD}
F_\ga(z_0+a,\ldots,z_n+a)=F_\ga(z_0,\ldots,z_n)
\end{equation}
the integral (\ref{eq:integral FD}) is in fact a function
of just $n$ variables: $z_1-z_0,\ldots,z_n-z_0$.

GKZ theory can deal efficiently with (multiplicative) torus actions
on the variables, but it can not accommodate for translation invariance like (\ref{eq:translation FD}). So we eliminate
the translation invariance during the passage to GKZ and consider 
the integrals (with the same $\mu_0,\ldots,\mu_n$)
\begin{equation}\label{eq:GKZ integral FD}
I_\gs(u_1,\ldots,u_{2n})
\,=\,\int_{\gs}(u_1+u_{n+1}\xi)^{-\mu_1}\cdots
(u_n+u_{2n}\xi)^{-\mu_n}\xi^{-\mu_0}\,d\xi\,,
\end{equation}
which are of the type considered in Section \ref{subsection:Euler}.

The GKZ differential equations satisfied by these integrals can be found
with the methods used in Section \ref{subsection:integrand k}.
For instance, for $j=1,\ldots,n$
\begin{eqnarray*}
\frac{\partial I_\gs}{\partial u_j}&=&-\mu_j
\int_{\gs}(u_1+u_{n+1}\zeta)^{-\mu_1}\cdots
(u_n+u_{2n}\zeta)^{-\mu_n}\xi^{-\mu_0}
\,\frac{d\zeta}{u_j+u_{j+n}\zeta}
 \\
\frac{\partial I_\gs}{\partial u_{j+n}}&=&
-\mu_j
\int_{\gs}(u_1+u_{n+1}\zeta)^{-\mu_1}\cdots
(u_n+u_{2n}\zeta)^{-\mu_n}\xi^{-\mu_0}
\,\frac{\zeta d\zeta}{u_j+u_{j+n}\zeta}
\end{eqnarray*}
and, hence, for $i,\,j=1,\ldots,n$ 
$$
\frac{\partial^2 I_\gs}{\partial u_i\partial u_{j+n}}\,=\,
\frac{\partial^2 I_\gs}{\partial u_j\partial u_{i+n}}\,,
$$
i.e.
$I_\gs$ satisfies the differential equations (\ref{eq:GKZ1}) with $\LL$ as in (\ref{eq:FD L matrix}) and $k=n$.

Similarly, for $s\in\CC$ close to $1$, we have
\begin{eqnarray*}
&&I_\gs(u_1,\ldots,u_n,su_{n+1},\ldots,su_{2n})\,=\,
s^{\mu_0-1}
I_\gs(u_1,\ldots,u_{2n})\,,\\
&&I_\gs(u_1,\ldots,u_{j-1},su_j,u_{j+1},\ldots,
u_{j+n-1},su_{j+n},u_{j+n+1},\ldots,u_{2n})\,= \\
&&\qquad=\,s^{-\mu_j}
I_\gs(u_1,\ldots,u_{2n})\,.
\end{eqnarray*}
This leads to the differential equations (\ref{eq:GKZ2}) 
with $k=n$, $\gA$ as in (\ref{eq:FD A matrix}) and 
$\vc=(\mu_0-1,-\mu_1,\ldots,-\mu_n)^t$.

As we have seen in Section \ref{subsubsection:Appell-Lauricella}
the power series $F_D(a,\vb,c |\vz)$ is, up to a constant factor,
the $\gG$-series associated with the above $\LL$ and with
$\quad\ugamma\,=\,
(\gam_1,\ldots,\gam_N)\,=\,(c-1,-b_1,\ldots,-b_{k-1},-a,0,\ldots,0)$.
The parameter $\vc$ in the GKZ differential equations (\ref{eq:GKZ2})
is therefore
$$
\vc\,=\,\sum_{j=1}^{2k}\gam_j\va_j\,=\,(-a,c-a-1,-b_1,\ldots,-b_{k-1})^t
\,=:\,(c_0,c_1,c_2,\ldots,c_k)^t.
$$
The system of differential equations (\ref{eq:GKZ2}) can now be written as
\begin{eqnarray}
\label{eq:FDDE1}
&& \frac{\partial  \Phi}{\partial u_{j+k}}\,=\,-u_{j+k}^{-1}\left(
u_j\frac{\partial  \Phi}{\partial u_j}-c_j\Phi\right)
\qquad\textrm{for}\quad
j=1,\ldots,k,
\\
\label{eq:FDDE2}
&& u_1\frac{\partial  \Phi}{\partial u_1}+
\ldots+\,u_k\frac{\partial  \Phi}{\partial u_k}\,=\,
(-c_0+c_1+\ldots+c_k)\Phi\,.
\end{eqnarray}
The system (\ref{eq:GKZ1}) is equivalent with the following  $\frac{1}{2}k(k-1)$ differential equations 
\begin{equation}\label{eq:FDDE3}
\frac{\partial^2 \Phi  }{\partial u_i\partial u_{j+k}}\,=\,
\frac{\partial^2  \Phi}{\partial u_j\partial u_{i+k}}
\qquad\textrm{for}\quad 1\leq i< j\leq k\,.
\end{equation}
Next we substitute (\ref{eq:FDDE1}) into (\ref{eq:FDDE3}) and set
$$
u_j=z_j\quad\textrm{if}\quad 1\leq j\leq k,\qquad
u_j=1\quad\textrm{if}\quad k+1\leq j\leq 2k\,.
$$
The result is the system of $\frac{1}{2}k(k-1)$ differential equations
\begin{equation}\label{eq:FDDE4}
(z_i-z_j)\frac{\partial^2 \Phi  }{\partial z_i\partial z_j}\,=\,
c_i\frac{\partial  \Phi}{\partial z_j}-
c_j\frac{\partial  \Phi}{\partial z_i}
\qquad\textrm{for}\quad 1\leq i< j\leq k\,.
\end{equation}
The above substitution turns (\ref{eq:FDDE2}) into
\begin{equation}\label{eq:FDDE5}
z_1\frac{\partial  \Phi}{\partial z_1}+
\ldots+\,z_k\frac{\partial  \Phi}{\partial z_k}\,=\,
(-c_0+c_1+\ldots+c_k)\Phi\,.
\end{equation}
The system of differential equations
(\ref{eq:FDDE4})-(\ref{eq:FDDE5})
is then equivalent with the GKZ system (\ref{eq:GKZ1})-(\ref{eq:GKZ2}) 
for Lauricella's $F_D$. The Equations (\ref{eq:FDDE4}) appear in this form also in \cite{loo}, and (\ref{eq:FDDE5}) appears in loc.cit. in an `integrated' form:
$$
\Phi(e^tz_1,\ldots,e^tz_k)=e^{(-c_0+c_1+\ldots+c_k)t}
\Phi(z_1,\ldots,z_k)\,.
$$
The match with \cite{loo} becomes exact, if one eliminates in
loc. cit. the translation invariance by setting $z_0=0$
( like we did in passing from (\ref{eq:integral FD}) to 
(\ref{eq:GKZ integral FD})).

\subsection{Triangulations of $\Delta_\gA$, secondary polytope and fan
for $F_D$}
\label{subsection:secondary polytope FD}
Consider a triangulation $\gT$ of the prism $\Delta_\gA$ by $k$-dimensional simplices with vertices in the set $\gA$. Then the bottom $(k-1)$-simplex $[\va_1,\ldots,\va_k]$ must be a face of exactly one $k$-simplex in the triangulation, say $\sigma_1$. Let 
$\va_{k+s_1}$ be the vertex of $\sigma_1$ opposite to the face
$[\va_1,\ldots,\va_k]$. So $1\leq s_1\leq k$. The face of $\sigma_1$ opposite to the vertex $s_1$ has vertices $\va_{k+s_1}$ and $\va_i$ with $1\leq i\leq k\,,\;i\neq s_1$. This must be a face of exactly one other $k$-dimensional simplex in the triangulation, say $\sigma_2$. 
Let $\va_{k+s_2}$ be the remaining vertex of $\sigma_2$. So 
$1\leq s_2\leq k$ and $s_2\neq s_1$. The face of $\sigma_2$ opposite to the vertex $\va_{s_2}$ has vertices $\va_{k+s_1},\,\va_{k+s_2}$
and $\va_i$ with $1\leq i\leq k\,,\;i\neq s_1\,,\,s_2$. This must be a face of exactly one other $k$-dimensional simplex, say $\sigma_3\,.$ Let $\va_{k+s_3}$ be the remaining vertex of $\sigma_3$. So 
$1\leq s_3\leq k$ and $s_3\neq s_1\,,\,s_2$. And so on. Thus  the triangulation $\gT$ of $\Delta_\gA$ determines a permutation $\tau$ of $\{1,\ldots,k\}$ with $\tau(i)\,=\,s_i\,.$

There is an obvious converse to this procedure associating to a permutation $\tau$ of $\{1,2,3,\ldots,k\}$ the triangulation with maximal
simplices $\sigma_1^{(\tau)},\ldots,\sigma_k^{(\tau)}$ where
\begin{equation}\label{eq:permtriang}
\sigma_j^{(\tau)}:=\textrm{convex hull }(\{\va_{\tau (i)} 
\:|\:j\leq i\leq k\}
\cup\{\va_{k+\tau (i)} \:|\:1\leq i\leq j\})\,.
\end{equation}

These triangulations are unimodular; i.e.
all $k$-simplices have \mbox{volume $1$.}
So when constructing the secondary polytope one only has to count for every triangulation how many simplices come together in the points
$\va_1,\ldots,\va_N$. With the above formula for the simplex
$\sigma_j^{(\tau)}$ one easily finds that the vector associated 
with the permutation $\tau$ is 
$(\tau^{-1}(1),\ldots,\tau^{-1}(k),
k+1-\tau^{-1}(1),\ldots,k+1-\tau^{-1}(k))$. 

The secondary polytope is the convex hull of these points as $\tau$ runs through all permutations of $\{1,2,3,\ldots,k\}$.
By translating over the vector corresponding to the identity permutation
the secondary polytope moves to the convex hull of the points
$(\tau^{-1}(1)-1,\ldots,\tau^{-1}(k)-k,
1-\tau^{-1}(1),\ldots,k-\tau^{-1}(k))$
in the space $\LL_\RR:=\LL\otimes\RR$.
\begin{figure}[t]
\begin{picture}(350,150)(100,50)
\put(0,0){\makebox(400,250){
\setlength\epsfxsize{7cm} 
\epsfbox{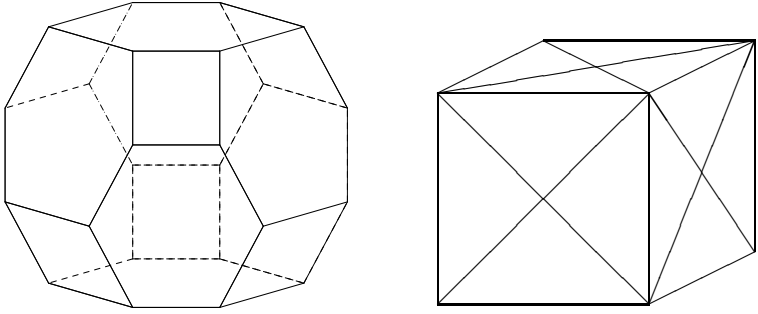}
}}
\put(300,50){
\setlength{\unitlength}{.5pt}
\begin{picture}(250,250)(13,-15)
\put(15,15){\line(0,1){160}}
\put(175,15){\line(0,1){160}}
\put(255,55){\line(0,1){160}}
\put(15,15){\line(1,0){160}}
\put(15,175){\line(1,0){160}}
\put(95,215){\line(1,0){160}}
\put(15,175){\line(2,1){80}}
\put(175,15){\line(2,1){80}}
\put(175,175){\line(2,1){80}}
\put(15,15){\line(1,1){160}}
\put(175,15){\line(-1,1){160}}
\put(15,175){\line(6,1){240}}
\put(95,215){\line(2,-1){80}}
\put(175,15){\line(2,5){80}}
\put(255,55){\line(-2,3){80}}
\end{picture}
\setlength{\unitlength}{1pt}
}
\end{picture}
\caption{\label{fig:FD Secondary}
Secondary polytope (left) and Secondary fan (right) for $F_D$ with $k=4$. All cones in the fan have their apex at the centre of the cube.
Shown are the intersections of the cones with some faces of the cube.
The reader is invited to label the maximal cones with the permutations of $\{1,2,3,4\}$.}
\end{figure}

\

\noindent
\textbf{Example/Exercise.}
The reader is invited to determine with the above algorithm 
the permutations corresponding to the maximal cones of
the secondary fan of Appell's $F_1$ (= Lauricella's $F_D$ with $k=3$) shown in Figure \ref{fig:fan F1}.

\

Recall from Section \ref{subsection:secondary fan construction}
that the secondary fan is a partition of the real vector space
$\LL_\RR^\vee=\mathrm{Hom}(\LL,\RR)$ into rational cones, all with
their apex in $\mf{0}$. Corollary \ref{cor:list} and Formula (\ref{eq:equivclass}) describe these cones. The vectors $\bb_1,\ldots,\bb_N$ are the images of the standard basis vectors of $\RR^N$ under the natural surjection
$\RR^N\longrightarrow \LL_\RR^\vee$.
In the present situation we choose the rows of the matrix
(\ref{eq:FD L matrix}) as a basis for $\LL_\RR=\LL\otimes\RR$.
On $\LL_\RR^\vee$ we use coordinates with respect to the dual basis.
The columns of (\ref{eq:FD L matrix}) then represent the vectors $\bb_1,\ldots,\bb_N$ in these coordinates.

Now consider a vector $\mf{t}=(t_2,\ldots,t_k)$ in $\LL_\RR^\vee$.
Put $t_1=0$. Then $\mf{t}$ defines a partial ordering 
$<_\mf{t}$ on the set $\{1,2,\ldots,k\}$ by
$$
i <_\mf{t} j\qquad\Leftrightarrow\qquad t_i<t_j\,.
$$ 
The indexing and ordering is such that for $h=1,\ldots,k$
\begin{equation}\label{eq:t lin}
\mf{t}\,=\,\sum_{i >_\mf{t} h}(t_i-t_h)\bb_{i+k}\,+\,
\sum_{i <_\mf{t} h}(t_h-t_i)\bb_i\,.
\end{equation}
One also easily checks that these are the only expresssions for
$\mf{t}$ as positive linear combination of a linearly independent
subset of $\{\bb_1,\ldots,\bb_N\}$.
Corollary \ref{cor:list} and Formula (\ref{eq:equivclass}) now tell exactly in which cone of the secondary fan $\mf{t}$ lies.
In particular, $\mf{t}$ lies in the interior of a maximal cone
if and only if no two of the numbers $t_1,t_2,\ldots,t_k$ are equal.
In that case, the ordering $<_\mf{t}$ is a total ordering, or what amounts to the same a permutation  of $\{1,2,\ldots,k\}$. More precisely,
if we associate with $\mf{t}$ the permutation $\tau$ defined by
$$
\tau(1)<_\mf{t}\tau(2)<_\mf{t}\ldots <_\mf{t}
\tau(k-1)<_\mf{t}\tau(k)\,,
$$
then the index set which effectively appears in (\ref{eq:t lin})
is the complement of the index set in (\ref{eq:permtriang})
with $h=\tau(j)$.

Thus we have shown:

\begin{corollary}\label{cor:FD fan}
The maximal cones in the secondary fan for $F_D$ are the connected components of the complement in $\RR^{k-1}$ of the union of hyperplanes
$$
\bigcup_{1\leq i<j\leq k} H_{ij}\,,\qquad\textrm{with equation for}\quad
H_{ij}:\quad t_i=t_j\,,
$$
where on the right $t_2,t_3,\ldots,t_k$ are
coordinates on $\RR^{k-1}$ and $t_1=0$.
\qed
\end{corollary}

\

\noindent
\textbf{Remarks.} Most GKZ systems do not have a secondary fan which is
cut out by a hyperplane arrangement. $F_D$ is something special.

\

\noindent
\textbf{Exercise.}
In Looijenga's lectures \cite{loo} the natural domain of definition
for the \emph{Schwarz map} is $\PP(V_n^\circ)$. The notations are:
$\PP(V_n^\circ)$ is $(\CC^{n+1})^\circ$ modulo the natural
$\CC^*$-action with weights $(1,1,\ldots,1)$ and modulo translations over $\CC(1,1,\ldots,1)$,
$$
(\CC^{n+1})^\circ:=\CC^{n+1}\setminus\bigcup_{i<j}\{
\textrm{hyperplane with equation }\;z_i=z_j\}\,.
$$
\emph{How does $\PP(V_n^\circ)$ relate to the toric variety $\cV_\gA$?}

\section{A glimpse of Mirror Symmetry}
\label{section:mirror symmetry}
\subsection{GKZ data from Calabi-Yau varieties}
\label{subsection:GKZCY}

One of the manifestations of the Mirror Symmetry phenomenon is a relation between two families of $3$-dimensional Calabi-Yau varieties,
matching complex geometry on one family with symplectic geometry on the other.

In the language of complex geometry a smooth Calabi-Yau variety is a compact
smooth K\"ahler manifold X with trivial canonical bundle, i.e.
$\Omega_X^{\dim X}\simeq\cO_X $,
which also satisfies
$H^0(X,\Omega_X^i)=0$ for $0<i<\dim X$.
Note: not all definitions in the literature require this second condition. Moreover, there are definitions which allow certain types of singularities.

A Calabi-Yau variety of dimension $1$ is an \textit{elliptic curve}.
A Calabi-Yau variety of dimension $2$ is a \textit{K3 surface}.
A Calabi-Yau variety of dimension $3$ is usually called a
\textit{Calabi-Yau threefold}.
Standard examples of Calabi-Yau varieties, all given as complete intersections in a product of projective spaces, are shown in the 
second column of Table \ref{table CY}.
From the homogeneous degrees of the defining equations and
the coordinates of the ambient projective space one builds a lattice
$\LL$ for use in GKZ context. This is shown in the third column 
of Table \ref{table CY}.
The lattice $\LL$ comes naturally with an embedding into some $\ZZ^N$
and the quotient $\MM:=\modquot{\ZZ^N}{\LL}$ is torsion free,
isomorphic to $\ZZ^{k+1}$, $k+1=N-d$. As in Section 
\ref{section:secondary fan} we let $\va_1,\ldots,\va_N\in\MM$ 
denote the images of the standard basis vectors of $\ZZ^N$
and $\gA=\{\va_1,\ldots,\va_N\}$. 
\begin{table}
\begin{center}
\begin{tabular}{|c|c|l|}\hline 
dim. & Calabi-Yau variety&\hspace{4em}$\mf{B}$
\\
\hline
1&cubic curve in $\PP^2$&$(-3,1,1,1)$
\\
1&$\bigcap$ two quadrics in $\PP^3$&$(-2,-2,1,1,1,1)$
\\
1&curve of degree $(2,2)$ in $\PP^1\times\PP^1$&
$\left(\begin{array}{rrrrr}
-2&1&0&1&0\\
-2&0&1&0&1
\end{array}\right)$
\\
1&$\bigcap$ two surf. deg. $(1,1,1)$ in $(\PP^1)^3$&
$\left(\begin{array}{rrrrrrrr}
-1&-1&1&0&0&1&0&0\\
-1&-1&0&1&0&0&1&0\\
-1&-1&0&0&1&0&0&1
\end{array}\right)$
\\
\hline
2&quartic surface in $\PP^3$&$(-4,1,1,1,1)$
\\
2&$\bigcap$ quadric and cubic in $\PP^4$&$(-2,-3,1,1,1,1,1)$
\\
2&$\bigcap$ three quadrics in $\PP^5$&$(-2,-2,-2,1,1,1,1,1,1)$
\\
2&surface of deg. $(2,2,2)$ in $(\PP^1)^3$&
$\left(\begin{array}{rrrrrrr}
-2&1&0&0&1&0&0\\
-2&0&1&0&0&1&0\\
-2&0&0&1&0&0&1
\end{array}\right)$
\\
\hline
3&quintic hypersurface in $\PP^4$&$(-5,1,1,1,1,1)$
\\
3&$\bigcap$ two cubics in $\PP^5$&
$(-3,-3,1,1,1,1,1,1)$
\\
3&$3$-fold of deg. $(3,3)$ in $\PP^2\times\PP^2$&
$\left(\begin{array}{rrrrrrr}
-3&1&0&1&0&1&0\\
-3&0&1&0&1&0&1
\end{array}\right)$
\\
3&$\bigcap$ four quadrics in $\PP^7$&
$(-2,-2,-2,-2,1,1,1,1,1,1,1,1)$
\\
3&$3$-fold of deg. $(2,2,2,2)$ in $(\PP^1)^4$&
$\left(\begin{array}{rrrrrrrrr}
-2&1&0&0&0&1&0&0&0\\
-2&0&1&0&0&0&1&0&0\\
-2&0&0&1&0&0&0&1&0\\
-2&0&0&0&1&0&0&0&1
\end{array}\right)$
\\
\hline
&$\bigcap$ means `intersection of'.&$\LL$ = $\ZZ$-span of rows of $\mf{B}$
\\
\hline
\end{tabular}
\caption{Standard examples of Calabi-Yau varieties.}\label{table CY}
\end{center}
\end{table}
Using Corollary \ref{cor:list} and Formula (\ref{eq:equivclass}) one
checks that in these examples the positive span of the last $d$ columns of matrix $\mf{B}$ is a maximal cone $\cC$ in the secondary fan of $\LL$.
Using (\ref{eq:detA=detB}) one checks that the triangulation of $\Delta_\gA$ corresponding to $\cC$ is unimodular.
Next one computes the ring $\cR_{\gA,\cC}$ in
Definition \ref{eq:ring R} and one finds that it
is (isomorphic to) the cohomology ring of the 
ambient space in the second column of Table \ref{table CY}:
$$
\cR_{\gA,\cC}\,=\,\left\{
\begin{array}{lll}
\modquot{\ZZ[\eps]}{(\eps^{r+1})}&\qquad&
\textrm{if the ambient space is }\quad \PP^r \\[1ex]
\modquot{\ZZ[\delta_1,\ldots,\delta_r]}{(\delta_1^2,\ldots,\delta_r^2)}
&\qquad&
\textrm{if the ambient space is }\quad (\PP^1)^r\\[1ex]
\modquot{\ZZ[\delta_1,\delta_2]}{(\delta_1^3,\delta_2^3)}
&\qquad&
\textrm{if the ambient space is }\quad (\PP^2)^2
\end{array}\right.
$$
Choosing a $\ZZ$-basis for $\MM$ we write 
$\va_j=(a_{1j},\ldots,a_{(k+1)\,j})$.
With $\gA$ one now associates the \emph{Laurent polynomial}
in the variables $x_1,\ldots,x_{k+1}$ with undetermined coefficients $u_j=\vu_{\va_j}$ (cf.(\ref{eq:nvar poly})): 
\begin{equation}\label{eq:mirror polynomial}
P_\gA(\vx)\,=\,P_\gA(x_1,\ldots,x_{k+1})\,=\,\sum_{\va\in\gA} \vu_\va\vx^\va\,=\,
\sum_{j=1}^N u_j\prod_{i=1}^{k+1}x_i^{a_{ij}}\,.
\end{equation}
Since for each $\va_j$ the coordinates sum to $1$,
this Laurent polynomial is homogeneous: $P_\gA(t\vx)\,=\,tP_\gA(\vx)$
for every $t\in\CC^*$.
As the coefficients $\vu$ vary the zero loci of $P_\gA(\vx)$
sweep out a family of hypersurfaces in 
$\modquot{(\CC^*)^{k+1}}{\CC^*}\,=\,(\CC^*)^k$.
Both $(\CC^*)^k$ and the hypersurfaces can be suitably compactified.
This family of compactified hypersurfaces is then the \emph{mirror
in the sense of \cite{BB} of the family of Calabi-Yau varieties}
in the second column in Table \ref{table CY}.
The members of this mirror family are
Calabi-Yau varieties if the original Calabi-Yau varieties have codimension $1$ in the ambient space. In case of codimension $>1$ the mirror
family consists of 
\emph{generalized Calabi-Yau varieties} in the sense of \cite{BB}.

We will now discuss details for the first three examples of Calabi-Yau threefolds. Other examples can be treated in the same way.

\subsection{The quintic in $\PP^4$}
\label{subsection:quintic}
This is the original example with which Mirror Symmetry entered
the mathematical arena; see \cite{COGP}.
The matrix $\mf{B}$ in Table \ref{table CY} is of the form
$\mf{B}\,=\,(\widetilde{\mf{B}}\;\II_d)\,,$
where $\II_d$ is the $d\times d$-identity matrix. The matrix
$\mf{A}=(\II_{N-d}\;-\widetilde{\mf{B}}^t)$ then satisfies
$\mf{B}\mf{A}^t=0$ and its columns generate $\ZZ^{k+1}$.
We apply row operations (i.e. a basis transformation in $\ZZ^{k+1}$)
so that the Laurent polynomial $P_\gA(\vx)$ in 
(\ref{eq:mirror polynomial}) assumes a pleasant form:
$$
\mf{A}=\left(\begin{array}{rrrrrr}
1&0&0&0&0&5\\
0&1&0&0&0&-1\\
0&0&1&0&0&-1\\
0&0&0&1&0&-1\\
0&0&0&0&1&-1
\end{array}\right)
\quad\leadsto\quad
\left(\begin{array}{rrrrrr}
1&1&1&1&1&1\\
0&1&0&0&0&-1\\
0&0&1&0&0&-1\\
0&0&0&1&0&-1\\
0&0&0&0&1&-1
\end{array}\right)
$$
We let $\gA\,=\,\{\va_1,\ldots,\va_N\}$ denote the columns of the
right-hand matrix. Then
$$
P_\gA(\vx)\,=\, x_1\left(u_1+u_2x_2+u_3x_3+u_4x_4+u_5x_5
+u_6(x_2x_3x_4x_5)^{-1}\right)\,.
$$

\

\noindent
\textbf{Remark.} The Laurent polynomial $P_\gA(\vx)$ can be
dehomogenized by setting $x_1=1$ and subsequently be homogenized
to the degree $5$ polynomial in $5$ variables
$$
\begin{array}{lll}
\widetilde{P}_\gA(\mf{X})&=&
u_1X_1X_2X_3X_4X_5+u_2X_2^2X_3X_4X_5+u_3X_2X_3^2X_4X_5+\\
&& +u_4X_2X_3X_4^2X_5+u_5X_2X_3X_4X_5^2
+u_6X_1^5\,.
\end{array}
$$
The polynomial $\widetilde{P}_\gA(\mf{X})$ defines a family
of special quintic hypersurfaces in $\PP^4$, which is the mirror of 
the family of general quintic hypersurfaces in $\PP^4$.
Traditionally the mirror family is presented as a quotient
of the hypersurface
$$
Z_1^5+Z_2^5+Z_3^5+Z_4^5+Z_5^5-5\psi Z_1Z_2Z_3Z_4Z_5=0
$$
by a specific action of the group $\left(\modquot{\ZZ}{5\ZZ}\right)^3$
(see \cite{COGP} \S 2).
The hypergeometric integrals and series constructed from the
periods of this `Fermat-like quintic' are however the same as those
coming from $\widetilde{P}_\gA(\mf{X})$.

\

The periods of the mirror Calabi-Yau hypersurfaces are given by integrals
\begin{equation}\label{eq:quintic period}
I^-_\gs(\bu)
\,:=\,\frac{1}{(2\pi \ci)^4}\int_\gs P_\gA(1,x_2,x_3,x_4,x_5)^{-1}\,
\frac{dx_2}{x_2}\frac{dx_3}{x_3}\frac{dx_4}{x_4} \frac{dx_5}{x_5}\,.
\end{equation}
As shown in Section \ref{subsection:integrand k}, these integrals
viewed as functions of $u_1,\ldots,u_6$, satisfy the GKZ system of differential equations (\ref{eq:GKZ1})-(\ref{eq:GKZ2}) with
$\gA$ as above and $\vc=-\va_1$.

If the numbers $\left|u_ju_1^{-1}\right|$ for $2\leq j\leq 6$ are 
sufficiently small, the domain of integration for one of the above period integrals $I^-_\gs(\bu)$ can be taken to be 
$\gs=\{|x_2|=|x_3|=|x_4|=|x_5|=1\}$. Using geometric series,
the binomial and residue theorems, one obtains for this period integral
the series expansion:
\begin{equation}\label{eq:quintic series}
I^-_\gs(\bu)\,=\,u_1^{-1}\sum_{n\geq 0} (-1)^n
\frac{(5n)!}{(n!)^5} z^n\qquad\textrm{with}\quad
z=u_1^{-5}u_2u_3u_4u_5u_6\,.
\end{equation}
Now look at the series $\Phi_{\LL,\ugamma,\ueps}(\vu)$
for $\LL=\ZZ(-5,1,1,1,1,1)$ and $\ugamma=(-1,0,0,0,0,0)$
defined in (\ref{eq:Gamma series deform}). In this case
$\ueps=(\eps_1,\eps_2,\eps_3,\eps_4,\eps_5,\eps_6)=
(-5\eps,\eps,\eps,\eps,\eps,\eps)$. 
The triangulation is given by the list $T_\cC$ consisting of the 
five sets one gets
by deleting from $\{1,2,3,4,5,6\}$ one number $>1$.
The minimal set not contained in a set on the list $T_\cC$ is
$\{2,3,4,5,6\}$. Thus we see
$$
\cR_{\gA,\cC}\,=\,
\ZZ\oplus\ZZ\eps\oplus\ZZ\eps^2\oplus\ZZ\eps^3\oplus\ZZ\eps^4
\,,\qquad \eps^5=0\,,
$$
and, with Pochhammer symbol notation (\ref{eq:Pochhammer})
and $z=u_1^{-5}u_2u_3u_4u_5u_6$,
\begin{equation}\label{eq:quintic series deform}
\Phi_{\LL,\ugamma,\ueps}(\vu)\,=\,-5\eps u_1^{-1}
\sum_{n\geq 0} (-1)^n
\frac{(1+5\eps)_{5n}}{((1+\eps)_n)^5} z^{n+\eps}\,.
\end{equation}
So, the function $\Phi_{\LL,\ugamma,\ueps}(\vu)$ takes values in the 
vector space $\eps\cR_{\gA,\cC}\otimes\CC$. 
Now let $\mathrm{Ann}(\eps)=\{x\in\cR_{\gA,\cC}\;|\;x\eps\,=\,0\}$
denote the annihilator ideal of $\eps$ and let 
$\oR_{\gA,\cC}:=\modquot{\cR_{\gA,\cC}}{\mathrm{Ann}(\eps)}$.
Then as a vector space $\eps\cR_{\gA,\cC}\otimes\CC$
is isomorphic to $\oR_{\gA,\cC}\otimes\CC$.
The latter has, however, the advantage of being a ring.
Let $\oeps$ denote the class of $\eps$ in 
$\oR_{\gA,\cC}$.
Then
$$
\oR_{\gA,\cC}\,=\,
\ZZ\oplus\ZZ\oeps\oplus\ZZ\oeps^2\oplus\ZZ\oeps^3
\,,\qquad \oeps^4=0\,.
$$
Moreover we can write
\begin{equation}\label{eq:quintic series deform reduce}
\Phi_{\LL,\ugamma,\ueps}(\vu)\,=\,-5u_1^{-1}z^{\oeps}
\sum_{n\geq 0} (-1)^n
\frac{(1+5\oeps)_{5n}}{((1+\oeps)_n)^5} z^n\,,
\end{equation}
and view it as a function with values in the ring
$\oR_{\gA,\cC}\otimes\CC$.
We expand this function with respect to
the basis $\{1,\oeps,\oeps^2,\oeps^3\}$:
\begin{equation}\label{eq:quintic series expand}
\Phi_{\LL,\ugamma,\ueps}(\vu)\,=\,
\Phi_0(\vu)\:+\:\Phi_1(\vu)\oeps\:+\:\Phi_2(\vu)\oeps^2
\:+\:\Phi_3(\vu)\oeps^3\,.
\end{equation}
By-passing all motivations, justifications and interpretations from 
string theory and Hodge theory (see, however, Section 
\ref{subsection:Schwarz} and \cite{CK} p. 263) we define the \emph{canonical coordinate}
\begin{equation}\label{eq:quintic q}
q\,:=\,-\exp\left(\frac{\Phi_1(\vu)}{\Phi_0(\vu)}\right)
\end{equation}
and the \emph{prepotential}
\begin{equation}\label{eq:quintic prepotential}
\cF (q)\,=\,\frac{5}{2}\left(
\frac{\Phi_1(\vu)}{\Phi_0(\vu)}\frac{\Phi_2(\vu)}{\Phi_0(\vu)}
\:-\:\frac{\Phi_3(\vu)}{\Phi_0(\vu)}
\right)\,.
\end{equation}
Note that these are functions of $z$, because 
they are constructed from quotients of solutions to the same GKZ 
system. In fact $q\,=\, -z+O(z^2)$
and we can invert this relation so as to get $z$ as a function $z(q)$ of $q$. We then want to view the prepotential as a function of $q$.
The recipe for extracting
results about the enumerative geometry of the general quintic 
threefold is to take
\begin{equation}\label{eq:quintic expansion}
\cF (q)\,=\,\textstyle{\frac{5}{6}}\log^3 q \:+\:
\displaystyle{\sum_{j\geq 1}}\: N_j\,\mathrm{Li}_3 (q^j)\,,
\end{equation}
where $\mathrm{Li}_3$ is the \emph{trilogarithm function}
$\mathrm{Li}_3 (x)\,:=\,\sum_{n\geq 1}\frac{x^n}{n^3}$.
Then one of the miracles of mirror symmetry is that all numbers $N_j$ are positive integers and that in fact $N_j$ equals the number of rational 
curves of degree $j$ on a general quintic threefold \cite{COGP,CK}.
The first few of these $N_j$ are shown in Table \ref{table quintic}.
\begin{table}[t]
\begin{center}
$$
\begin{array}{lll}
N_1&=& 2875\\
N_2&=& 609250\\
N_3&=& 317206375\\
N_4&=& 242467530000\\
N_5&=& 229305888887625\\
N_6&=& 248249742118022000\\
N_7&=& 295091050570845659250\\
N_8&=& 375632160937476603550000\\
N_9&=& 503840510416985243645106250
\end{array}
$$
\caption{The numbers $N_j$ for the quintic threefold.}
\label{table quintic}
\end{center}
\end{table}

Actual computations proceed as follows: compute $F_0,\ldots,f_3$ from
$$
\sum_{n\geq 0} (-1)^n
\frac{(1+5\oeps)_{5n}}{((1+\oeps)_n)^5} z^n\,=\,
F_0(z)\:+\:F_1(z)\oeps\:+\:F_2(z)\oeps^2
\:+\:F_3(z)\oeps^3
$$
and $f_i(z):=\frac{F_i(z)}{F_0(z)}$. 
That implies
\begin{equation}\label{eq:efficient quintic expand}
\Phi_{\LL,\ugamma,\ueps}(\vu)\,=\,
\Phi_0(\vu)z^{\oeps}\left(1\:+\:f_1(z)\oeps\:+\:f_2(z)\oeps^2
\:+\:f_3(z)\oeps^3\right)
\end{equation}
Comparing the expansions of $\log \Phi_{\LL,\ugamma,\ueps}(\vu)$
which result from (\ref{eq:quintic series expand}) and
(\ref{eq:efficient quintic expand}), i.e. 
\begin{eqnarray*}
\log \Phi_{\LL,\ugamma,\ueps}(\vu)&=&\log(\Phi_0(z))\:+\:
\frac{\Phi_1(\vu)}{\Phi_0(\vu)}\oeps\:+\:
\left(-\frac{1}{2}\left[\frac{\Phi_1(\vu)}{\Phi_0(\vu)}\right]^2
\,+\,\frac{\Phi_2(\vu)}{\Phi_0(\vu)}\right)\oeps^2
\\
&& \:+\:
\left(\frac{1}{3}\left[\frac{\Phi_1(\vu)}{\Phi_0(\vu)}\right]^3\,-\,
\frac{\Phi_1(\vu)}{\Phi_0(\vu)}\frac{\Phi_2(\vu)}{\Phi_0(\vu)}
\:+\:\frac{\Phi_3(\vu)}{\Phi_0(\vu)}
\right)\oeps^3
\\
&=&
\log(\Phi_0(z))\:+\:(\log z\,+\,f_1(z))\,\oeps\:+\:
(-\textstyle{\frac{1}{2}}f_1(z)^2\,+\,f_2(z))\,\oeps^2
\\
&& +\:
(\textstyle{\frac{1}{3}}f_1(z)^3\,-\,f_1(z)f_2(z)\,+\,f_3(z))\,\oeps^3\,,
\end{eqnarray*}
we see that $q$ and $\cF (q)$ can easily be computed from the already 
known $f_1,\,f_2,\,f_3$:
\begin{eqnarray*}
q&=&-z\,\exp(f_1(z))\,,
\\
\cF (q)&=&\textstyle{\frac{5}{2}}\left(
\textstyle{\frac{1}{3}}\log^3 (-q) \:-\:
\left(\textstyle{\frac{1}{3}}f_1(z(q))^3\,-\,f_1(z(q))f_2(z(q))
\,+\,f_3(z(q))\right)\right)\!.
\end{eqnarray*}

\subsection{The intersection of two cubics in $\PP^5$}
\label{subsection:cubics intersect}
This is one of the examples discussed in \cite{LT}.
Here we treat it as a highly resonant GKZ system.
As in the case of the quintic
the matrix $\mf{B}$ in Table \ref{table CY} is of the form
$\mf{B}\,=\,(\widetilde{\mf{B}}\;\II_d)$ and we apply row operations to
the matrix $(\II_{N-d}\;-\widetilde{\mf{B}}^t)$  
so that the Laurent polynomial $P_\gA(\vx)$ in 
(\ref{eq:mirror polynomial}) assumes a pleasant form:
$$
\left(\begin{array}{rrrrrrrr}
1&0&0&0&0&0&0&3\\
0&1&0&0&0&0&0&3\\
0&0&1&0&0&0&0&-1\\
0&0&0&1&0&0&0&-1\\
0&0&0&0&1&0&0&-1\\
0&0&0&0&0&1&0&-1\\
0&0&0&0&0&0&1&-1
\end{array}\right)
\quad\leadsto\quad
\left(\begin{array}{rrrrrrrr}
1&0&1&1&0&0&0&1\\
0&1&0&0&1&1&1&0\\
0&0&1&0&0&0&0&-1\\
0&0&0&1&0&0&0&-1\\
0&0&0&0&1&0&0&-1\\
0&0&0&0&-1&1&0&0\\
0&0&0&0&-1&0&1&0
\end{array}\right)
$$
We let $\gA\,=\,\{\va_1,\ldots,\va_N\}$ denote the columns of the
right-hand matrix. Then
$$
P_\gA(\vx)\:=\:x_1\, P_{\gA,1}(x_3,x_4,x_5)\:+\:
x_2\, P_{\gA,2}(x_5,x_6,x_7)
$$
with
\begin{eqnarray*}
P_{\gA,1}(x_3,x_4,x_5)&=&
u_1+u_3x_3+u_4x_4+u_8(x_3x_4x_5)^{-1}
\\
P_{\gA,2}(x_5,x_6,x_7)&=&u_2+u_5x_5(x_6x_7)^{-1}+u_6x_6+u_7x_7\,.
\end{eqnarray*}
This way of combining the two Laurent polynomials in five variables,
$P_{\gA,1}$ and $P_{\gA,2}$, to one Laurent polynomial in seven
variables $P_\gA$ is known as \emph{Cayley's trick} (see \cite{gkz2, gkz3,gkz4}).
The two polynomials $P_{\gA,1}$ and $P_{\gA,2}$, suitably homogenized,
define a family of Calabi-Yau complete intersection threefolds in
$\PP^1\times\PP^2\times\PP^2$.
The corresponding period integrals are 
(cf. Section \ref{subsection:Euler} and \cite{gkz2,gkz4})
\begin{equation}\label{eq:cubics period}
I^-_\gs(\bu)
\,:=\,
\frac{1}{(2\pi \ci)^5}\int_\gs 
P_{\gA,1}(x_3,x_4,x_5)^{-1}\,
P_{\gA,2}(x_5,x_6,x_7)^{-1}\,
\frac{dx_3}{x_3}\frac{dx_4}{x_4} \frac{dx_5}{x_5}
\frac{dx_6}{x_6}\frac{dx_7}{x_7}\,.
\end{equation}
One can show as in Sections \ref{subsection:integrand k}
and \ref{subsection:Euler}
that these integrals
viewed as functions of $u_1,\ldots,u_8$, satisfy the GKZ system of differential equations (\ref{eq:GKZ1})-(\ref{eq:GKZ2}) with
$\gA$ as above and $\vc=-\va_1-\va_2$.

If the numbers $\left|u_3u_1^{-1}\right|$, $\left|u_4u_1^{-1}\right|$,
$\left|u_8u_1^{-1}\right|$ and $\left|u_5u_2^{-1}\right|$,
$\left|u_6u_2^{-1}\right|$, $\left|u_7u_2^{-1}\right|$ are 
sufficiently small, the domain of integration for one of the above period integrals $I^-_\gs(\bu)$ can be taken to be 
$\gs=\{|x_3|=|x_4|=|x_5|=|x_6|=|x_7|=1\}$. 
This period integral admits the series expansion, with 
$z=u_1^{-3}u_2^{-3}u_3u_4u_5u_6u_7u_8$,
\begin{equation}\label{eq:2cubics series}
I^-_\gs(\bu)\,=\,u_1^{-1}u_2^{-1}\sum_{n\geq 0} 
\left(\frac{(3n)!}{(n!)^3}\right)^2 z^n\,.
\end{equation}
The series $\Phi_{\LL,\ugamma,\ueps}(\vu)$
for $\LL=\ZZ(-3,-3,1,1,1,1,1,1)$,
$\ugamma=(-1,-1,0,0,0,0,0,0)$
and
$\ueps=(\eps_1,\ldots,\eps_8)=
(-3\eps,-3\eps,\eps,\eps,\eps,\eps,\eps,\eps)$ reads
\begin{equation}\label{eq:2cubics series deform}
\Phi_{\LL,\ugamma,\ueps}(\vu)\,=\,9\eps^2 u_1^{-1}u_2^{-1}
\sum_{n\geq 0} 
\left(\frac{(1+3\eps)_{3n}}{((1+\eps)_n)^3}\right)^2 z^{n+\eps}\,,
\end{equation}
and is evaluated in
$
\cR_{\gA,\cC}\,=\,
\ZZ\oplus\ZZ\eps\oplus\ZZ\eps^2\oplus\ZZ\eps^3\oplus\ZZ\eps^4
\oplus\ZZ\eps^5\,,\qquad \eps^6=0\,.
$
\\
The function $\Phi_{\LL,\ugamma,\ueps}(\vu)$ actually takes 
values in the vector space $\eps^2\cR_{\gA,\cC}\otimes\CC$. 
As in the case of the quintic, we replace this space 
by the isomorphic space $\oR_{\gA,\cC}\otimes\CC$, where $\mathrm{Ann}(\eps^2):=\{x\in\cR_{\gA,\cC}\;|\;x\eps^2\,=\,0\}$
and $\oR_{\gA,\cC}:=\modquot{\cR_{\gA,\cC}}{\mathrm{Ann}(\eps^2)}$.
Let $\oeps$ denote the class of $\eps$ in 
$\oR_{\gA,\cC}$.
Then
$$
\oR_{\gA,\cC}\,=\,
\ZZ\oplus\ZZ\oeps\oplus\ZZ\oeps^2\oplus\ZZ\oeps^3
\,,\qquad \oeps^4=0\,.
$$
Proceeding as in the case of the quintic we write
\begin{eqnarray*}
\Phi_{\LL,\ugamma,\ueps}(\vu)&=&
9u_1^{-1}u_2^{-1}
\sum_{n\geq 0} 
\left(\frac{(1+3\oeps)_{3n}}{((1+\oeps)_n)^3}\right)^2 z^{n+\oeps}
\\
&=&
\Phi_0(\vu)\:+\:\Phi_1(\vu)\oeps\:+\:\Phi_2(\vu)\oeps^2
\:+\:\Phi_3(\vu)\oeps^3\\
&=&\Phi_0(\vu)z^{\oeps}\left(1\:+\:f_1(z)\oeps\:+\:f_2(z)\oeps^2
\:+\:f_3(z)\oeps^3\right)\,,
\end{eqnarray*}
and almost exactly as for the quintic we extract from
$\log \Phi_{\LL,\ugamma,\ueps}(\vu)$ a \emph{canonical coordinate} and a \emph{prepotential}:
\begin{eqnarray}
\label{eq:cubics q}
q&:=&z\,\exp(f_1(z))\,=\, z+O(z^2)\,,
\\
\label{eq:cubics prepotential}
\cF (q)&:=&\textstyle{\frac{9}{2}}\left(
\textstyle{\frac{1}{3}}\log^3 q \:-\:
\left(\textstyle{\frac{1}{3}}f_1(z(q))^3\,-\,f_1(z(q))f_2(z(q))
\,+\,f_3(z(q))\right)\right)\!.
\end{eqnarray}
Finally we compute the numbers $N_j$ from the expansion
\begin{equation}\label{eq:cubics expansion}
\cF (q)\,=\,\textstyle{\frac{3}{2}}\log^3 q \:+\:
\displaystyle{\sum_{j\geq 1}}\: N_j\,\mathrm{Li}_3 (q^j)\,.
\end{equation}
The first few of the numbers $N_j$ are shown in Table 
\ref{table cubics} and agree with those in \cite{LT}.
\begin{table}[t]
\begin{center}
$$
\begin{array}{lll}
N_1&=& 1053\\
N_2&=& 52812\\
N_3&=& 6424326\\
N_4&=& 1139448384\\
N_5&=& 249787892583\\
N_6&=& 62660964509532\\
N_7&=& 17256453900822009\\
N_8&=& 5088842568426162960\\
N_9&=& 1581250717976557887945
\end{array}
$$
\caption{The numbers $N_j$ for the intersection of two cubics in $\PP^5$.}
\label{table cubics}
\end{center}
\end{table}
\subsection{The hypersurface of degree $(3,3)$ in $\PP^2\times\PP^2$}
\label{subsection:(3,3)}
Again the matrix $\mf{B}$ in Table \ref{table CY} is of the form
$\mf{B}\,=\,(\widetilde{\mf{B}}\;\II_d)\,,$ and we apply row operations
to the matrix $(\II_{N-d}\;-\widetilde{\mf{B}}^t)$:
$$
\left(\begin{array}{rrrrrrr}
1&0&0&0&0&3&3\\
0&1&0&0&0&-1&0\\
0&0&1&0&0&0&-1\\
0&0&0&1&0&-1&0\\
0&0&0&0&1&0&-1
\end{array}\right)
\quad\leadsto\quad
\left(\begin{array}{rrrrrrr}
1&1&1&1&1&1&1\\
0&1&0&0&0&-1&0\\
0&0&1&0&0&0&-1\\
0&0&0&1&0&-1&0\\
0&0&0&0&1&0&-1
\end{array}\right)
$$
We let $\gA\,=\,\{\va_1,\ldots,\va_N\}$ denote the columns of the
right-hand matrix. Then
$$
P_\gA(\vx)\,=\, x_1\left(u_1+u_2x_2+u_3x_3+u_4x_4+u_5x_5
+u_6(x_2x_4)^{-1}+u_7(x_3x_5)^{-1}\right)\,.
$$
The periods of the mirror Calabi-Yau hypersurfaces are given by integrals
\begin{equation}\label{eq:33period}
I^-_\gs(\bu)
\,:=\,\frac{1}{(2\pi \ci)^4}\int_\gs P_\gA(1,x_2,x_3,x_4,x_5)^{-1}\,
\frac{dx_2}{x_2}\frac{dx_3}{x_3}\frac{dx_4}{x_4} \frac{dx_5}{x_5}\,;
\end{equation}
As shown in Section \ref{subsection:integrand k}
these integrals
viewed as functions of $u_1,\ldots,u_7$, satisfy the GKZ system of differential equations (\ref{eq:GKZ1})-(\ref{eq:GKZ2}) with
$\gA$ as above and $\vc=-\va_1$.

If the numbers $\left|u_ju_1^{-1}\right|$ for $2\leq j\leq 7$ are 
sufficiently small, the domain of integration for one of the above period integrals $I^-_\gs(\bu)$ can be taken to be 
$\gs=\{|x_2|=|x_3|=|x_4|=|x_5|=1\}$.
This integral admits the expansion:
\begin{equation}\label{eq:(3,3) series}
I^-_\gs(\bu)\,=\,u_1^{-1}\sum_{n_1,\,n_2\geq 0} (-1)^{n_1+n_2}
\frac{(3n_1+3n_2)!}{(n_1!)^3(n_2!)^3} z_1^{n_1}z_2^{n_2}
\end{equation}
with $z_1=u_1^{-3}u_2u_4u_6$ and $z_2=u_1^{-3}u_3u_5u_7$.
Now look at the series $\Phi_{\LL,\ugamma,\ueps}(\vu)$
for 
$\LL=\ZZ(-3,1,0,1,0,1,0)\oplus\ZZ(-3,0,1,0,1,0,1)$,
$\ugamma=(-1,0,0,0,0,0,0)$ 
and\\
$\ueps=(\eps_1,\eps_2,\eps_3,\eps_4,\eps_5,\eps_6,\eps_7)=
\delta_1(-3,1,0,1,0,1,0)+\delta_2(-3,0,1,0,1,0,1)$. 
Using Corollary \ref{cor:list} and the matrix $\mf{B}$ from Table
\ref{table CY} for this example one easily checks that the triangulation
is given by the list $T_\cC$ consisting of the nine sets one gets
by deleting from $\{1,2,3,4,5,6,7\}$ one even and one odd number $>1$.
The minimal sets not contained in a set on the list $T_\cC$ are
$\{2,4,6\}$ and $\{3,5,7\}$. Thus we see
\begin{eqnarray*}
\cR_{\gA,\cC}&=&\ZZ\oplus\ZZ\delta_1\oplus\ZZ\delta_2\oplus
\ZZ\delta_1^2\oplus\ZZ\delta_1\delta_2\oplus\ZZ\delta_2^2\oplus
\ZZ\delta_1^2\delta_2\oplus\ZZ\delta_1\delta_2^2
\oplus\ZZ\delta_1^2\delta_2^2\,,
\\
&& \delta_1^3\,=\,\delta_2^3\,=\,0\,.
\end{eqnarray*}
Thus, with $z_1=u_1^{-3}u_2u_4u_6$ and $z_2=u_1^{-3}u_3u_5u_7$,
$$
\Phi_{\LL,\ugamma,\ueps}(\vu)\,=\,-3(\delta_1+\delta_2)
u_1^{-1}\hspace{-.5em}\sum_{n_1,\,n_2\geq 0} (-1)^{n_1+n_2}
\frac{(1+3\delta_1+3\delta_2)_{3n_1+3n_2}}
{((1+\delta_1)_{n_1}(1+\delta_2)_{n_2})^3} 
z_1^{n_1+\delta_1}z_2^{n_2+\delta_2}\,.
$$
The function $\Phi_{\LL,\ugamma,\ueps}(\vu)$ takes values in the 
vector space $(\delta_1+\delta_2)\cR_{\gA,\cC}\otimes\CC$. 
As in the previous cases, we replace this space by the isomorphic one $\oR_{\gA,\cC}\otimes\CC$, 
where $\mathrm{Ann}(\delta_1+\delta_2):=
\{x\in\cR_{\gA,\cC}\;|\;x(\delta_1+\delta_2)\,=\,0\}$
and $\oR_{\gA,\cC}:=
\modquot{\cR_{\gA,\cC}}{\mathrm{Ann}(\delta_1+\delta_2)}$.
Let $\ode_1$ and $\ode_2$ denote the classes of $\delta_1$ and
$\delta_2$, respectively, in $\oR_{\gA,\cC}$.
Then
\begin{eqnarray*}
\oR_{\gA,\cC}
&=&
\ZZ\oplus\ZZ\ode_1\oplus\ZZ\ode_2
\oplus\ZZ\ode_2^2\oplus\ZZ\ode_1^2\oplus\ZZ\ode_1^2\ode_2\,,
\\
\ode_1\ode_2=\ode_1^2+\ode_2^2\,,&&\hspace{-1em}
\ode_1^2\ode_2=\ode_1\ode_2^2\,,\qquad
\ode_1^3=\ode_2^3=\ode_1^2\ode_2^2=0\,.
\end{eqnarray*}
Proceeding as in the previous examples we write
\begin{eqnarray*}
\Phi_{\LL,\ugamma,\ueps}(\vu)&=&-3u_1^{-1}\sum_{n_1,\,n_2\geq 0} (-1)^{n_1+n_2}
\frac{(1+3\ode_1+3\ode_2)_{3n_1+3n_2}}
{((1+\ode_1)_{n_1}(1+\ode_2)_{n_2})^3} 
z_1^{n_1+\ode_1}z_2^{n_2+\ode_2}
\\
&&\hspace{-5em}=\,\Phi_0(\vu)\,+\,\Phi_{1,1}(\vu)\ode_1\,+\,
\Phi_{1,2}(\vu)\ode_2
\,+\,\Phi_{2,1}(\vu)\ode_2^2\,+\,\Phi_{2,2}(\vu)\ode_1^2\,+\,
\Phi_3(\vu)\ode_1^2\ode_2
\\
&&\hspace{-7em}=\,
\Phi_0(\vu)z_1^{\ode_1}z_2^{\ode_2}
\left(1\,+\,f_{1,1}(\vz)\ode_1\,+\,f_{1,2}(\vz)\ode_2
\,+\,f_{2,1}(\vz)\ode_2^2\,+\,f_{2,2}(\vz)\ode_1^2\,+\,
f_3(\vz)\ode_1^2\ode_2\right).
\end{eqnarray*}
Here $\vz=(z_1,z_2)$.
From the $\ode_1$ and $\ode_2$ components 
we construct two \emph{canonical coordinates}
(cf. (\ref{eq:can coord}))
\begin{equation}\label{eq:33cancoord}
q_1\,:=\,-z_1\exp(f_{1,1}(\vz))\,,\qquad 
q_2\,:=\,-z_2\exp(f_{1,2}(\vz))\,.
\end{equation}
We view $z_1,z_2$ as functions of $q_1,q_2$ via the inverse of 
relation (\ref{eq:33cancoord}).\\
The \emph{prepotential} in this case is 
(cf. (\ref{eq:prepotential}))
\begin{equation}\label{eq:33 prepotential}
\cF (q)\,=\,\frac{3}{2}\left(
\frac{\Phi_{1,1}(\vu)}{\Phi_0(\vu)}\frac{\Phi_{2,1}(\vu)}{\Phi_0(\vu)}\:+\:
\frac{\Phi_{1,2}(\vu)}{\Phi_0(\vu)}\frac{\Phi_{2,2}(\vu)}{\Phi_0(\vu)}
\:-\:\frac{\Phi_3(\vu)}{\Phi_0(\vu)}
\right)\,.
\end{equation}
The $-$-signs in (\ref{eq:33cancoord}) and the factor $3$ 
in (\ref{eq:33 prepotential}) are needed to
match the calculations below with the
results in \cite{HKTY} Appendix B2.

We expand $\log\Phi_{\LL,\ugamma,\ueps}(\vu)$ on
the basis 
$\{1,\ode_1,\,\ode_2,\,\ode_2^2,\,\ode_1^2,\,\ode_1^2\ode_2\}$
of $\oR_{\gA,\cC}$.
The $\ode_1^2\ode_2$-coordinate is on the one hand
$$
\log(-q_1)\log(-q_2)\log(q_1q_2)-\textstyle{\frac{2}{3}}\cF (q)
$$
and on the other hand it is
$$
f_{1,1}(\vz)^2f_{1,2}(\vz)+f_{1,1}(\vz)f_{1,2}(\vz)^2
\,-f_{1,1}(\vz)f_{2,1}(\vz)-f_{1,2}(\vz)f_{2,2}(\vz)\,+\,f_3(\vz).
$$

Computing the coefficients $N_{j_1,j_2}$ in the expansion
$$
\cF(q_1,q_2)\,=\,\textstyle{\frac{3}{2}}
\log(-q_1)\log(-q_2)\log(q_1q_2)\,+\,
\sum_{j_1,j_2\geq 0,\,j_1+j_2>0} N_{j_1,j_2}
\mathrm{Li}_3 (q_1^{j_1}q_2^{j_2})
$$
is now somewhat more involved than in the previous examples.
We leave it as an\\ 
\textbf{exercise in Mathematica, Maple or PARI programming}.\\
A table of the numbers $N_{j_1,j_2}$ for this example appears in \cite{HKTY} Appendix B2
under the name $X_{(3|3)}(1,1,1|1,1,1)$. In \cite{HKTY} one
finds many more $2$-parameter models.

\subsection{The Schwarz map for some extended GKZ systems}
\label{subsection:Schwarz}
In this section we briefly discuss how the 
$\oR_{\gA,\cC}\otimes\CC$-valued function 
$\Phi_{\LL,\ugamma,\ueps}(\vu)$ which we met in the preceding examples,
can be viewed as a Schwarz map and what are some special features of the image.

First note that, since $\dim\oR_{\gA,\cC}\otimes\CC\,<\,
\dim\cR_{\gA,\cC}\otimes\CC\,=\,\mathrm{volume}\: \Delta_\gA$,
the components of $\Phi_{\LL,\ugamma,\ueps}(\vu)$ with respect to some basis of $\oR_{\gA,\cC}\otimes\CC$ can not suffice as a basis for the solution space of the GKZ system. They do however constitute a basis for the solution space of some extension of the GKZ system (see \cite{HKTY}).
So, strictly speaking we are not talking about the Schwarz map for the
GKZ system, but for an extension thereof.
Since we do explicitly have all these basis solutions for the extended 
system, we need not care about this system itself.

In the examples, coming from (families of) Calabi-Yau threefolds,
the ring $\oR_{\gA,\cC}$ is graded and splits in homogeneous pieces,
$$
\oR_{\gA,\cC}\,=\,\oR_{\gA,\cC}^{(0)}\oplus\oR_{\gA,\cC}^{(1)}
\oplus\oR_{\gA,\cC}^{(2)}\oplus\oR_{\gA,\cC}^{(3)}\,,
$$
with degrees $0,\,1,\,2,\,3$ and ranks $1,\,d,\,d,\,1$, respectively;
recall $d\,=\,\rank\LL$.
We fix a basis for $\oR_{\gA,\cC}$ by fixing bases for the homogeneous pieces
$$
\oR_{\gA,\cC}^{(0)}:\,
e_0=1\,,\quad
\oR_{\gA,\cC}^{(1)}:\,
e_{1,1},\ldots,e_{1,d}\,,\quad
\oR_{\gA,\cC}^{(2)}:\,
e_{2,1},\ldots,e_{2,d}\,,\quad
\oR_{\gA,\cC}^{(3)}:\,
e_3\,,
$$
and expand $\Phi_{\LL,\ugamma,\ueps}(\vu)$ with respect to this basis
\begin{equation}\label{eq:Phi expand}
\Phi_{\LL,\ugamma,\ueps}(\vu)\,=\,\Phi_{0}(\vu)e_0\,+\,
\sum_{i=1}^d\Phi_{1,i}(\vu)e_{1,i}\,+\,
\sum_{i=1}^d\Phi_{2,i}(\vu)e_{2,i}\,+\,\Phi_{3}(\vu)e_3\,.
\end{equation}
The Schwarz map lands in the projective space
$$
\PP\left(\oR_{\gA,\cC}\otimes\CC\right)
$$
and $\Phi_{0}(\vu),\ldots,\Phi_{3}(\vu)$ are homogeneous coordinates
for the image points. Since these functions are solutions of the same GKZ system their quotients, and hence also the Schwarz map, are defined on
some open set in $\cV_\gA$ near the special point $\vp_\cC$ corresponding to the maximal cone $\cC$ in the secondary fan.
The map is multi-valued and we do fully control the local monodromy.

The image of the Schwarz map has dimension equal to
$\dim \cV_\gA\,=\,\rank \LL\,=\,d$, whereas the projective space
$\PP\left(\oR_{\gA,\cC}\otimes\CC\right)$ has dimension $1+2d$.
For the description of the image of the Schwarz map we want to
profit from the description of the moduli of Calabi-Yau 
threefolds by Bryant and Griffiths \cite{BG}.
In the theory of moduli of Calabi-Yau threefolds one writes the 
holomorphic $3$-form in coordinates with respect to a basis of the 
third cohomology space given by topological $3$-cycles. These coordinates
are the period integrals of the $3$-form. 
We know $\Phi_{0}(\vu)$ explicitly as a period integral 
(see (\ref{eq:quintic period}), (\ref{eq:cubics period}), 
(\ref{eq:33period})), but we still need an argument for the other 
coordinates in (\ref{eq:Phi expand}) to be periods.
Such an argument may be that inspection of the local monodromy shows
that the extended GKZ system of differential equations satisfied by
the known period integral $\Phi_{0}(\vu)$ is irreducible, for then all periods must be linear combinations of 
$\Phi_{0}(\vu),\ldots,\Phi_{3}(\vu)$.

Having matched (\ref{eq:Phi expand}) with the coordinates (= periods)
of the holomorphic $3$-form with respect to a basis of topological $3$-cycles, we must check that the basis $e_0,\ldots,e_3$
satisfies the requirements for application of the Bryant-Griffiths theory,
i.e. we need to know that with respect to the alternating bilinear form
$\langle,\rangle$ on the third cohomology space of the Calabi-Yau threefold
\begin{equation}\label{eq:intersection products}
\langle e_0,e_3\rangle\,=\,-\langle e_3,e_0\rangle\,=\,
-\langle e_{1,i},e_{2,i}\rangle\,=\,
\langle e_{2,i},e_{1,i}\rangle\,=\,1\quad\textrm{for}\;i=1,\ldots,d,
\end{equation}
and all other $\langle e_r,e_s\rangle\,=\,0$.
In an example at the end of this section we show how to derive
(\ref{eq:intersection products}) from the explicitly known local
monodromy and logarithmic pieces of $\Phi_{\LL,\ugamma,\ueps}(\vu)$.

We are now all set for applying \cite{BG}. First define the
\emph{canonical coordinates}
\begin{equation}\label{eq:can coord}
q_i\,:=\,\exp\left(\frac{\Phi_{1,i}(\vu)}{\Phi_0(\vu)}\right)
\qquad\textrm{for}\;i=1,\ldots d\,.
\end{equation}
The derivations $q_i\frac{\partial}{\partial q_i}$ act on the  cohomology spaces of the Calabi-Yau threefolds in the family.
Griffiths' transversality and the Riemann bilinear relations imply
\begin{equation}\label{eq:transversality}
\left\langle\frac{\Phi_{\LL,\ugamma,\ueps}(\vu)}{\Phi_0(\vu)},\,
q_i\frac{\partial}{\partial q_i}\left(
\frac{\Phi_{\LL,\ugamma,\ueps}(\vu)}{\Phi_0(\vu)}\right)
\right\rangle\,=\,0\,.
\end{equation}
Write $\gf_3\,:=\,\frac{\Phi_3(\vu)}{\Phi_0(\vu)}$ and
$\gf_{a,j}\,:=\,\frac{\Phi_{a,j}(\vu)}{\Phi_0(\vu)}$
for $a=1,\,2$, $j=1,\ldots,d$. These are (multivalued) functions
of $q_1,\ldots,q_d$, and in fact $\gf_{1,j}\,=\,\log q_j$. Then the left hand side of (\ref{eq:transversality}) evaluates to
$$
q_i\frac{\partial\gf_3}{\partial q_i}\,-\,\sum_{j=1}^d\left(
\gf_{1,j} q_i\frac{\partial\gf_{2,j}}{\partial q_i}\right)\,+\,
\gf_{2,i}\;=\;
q_i\frac{\partial}{\partial q_i}\left(\gf_3\,-\,\sum_{j=1}^d
\gf_{1,j}\gf_{2,j}\right)\,+\,2\,\gf_{2,i}\,.
$$
According to (\ref{eq:transversality}) this equals $0$ and thus
\begin{equation}\label{eq:special geometry}
\gf_{2,i}\,=\,q_i\frac{\partial\cF}{\partial q_i}
\end{equation}
where 
\begin{equation}\label{eq:prepotential}
\cF\,:=\,\frac{1}{2}\left(-\gf_3\,+\,\sum_{j=1}^d\gf_{1,j}\gf_{2,j}
\right)
\end{equation}
is the so-called \emph{prepotential}.

\

\noindent
\textbf{Example.}
Thus we recover the canonical coordinate (\ref{eq:quintic q})
for the quintic in $\PP^4$ up to a $-$-sign and the prepotential
(\ref{eq:quintic prepotential}) up to a factor $5$ (which is the
degree of the quintic). And similarly for the intersection of two 
cubics in $\PP^3$ and the hypersurface of degree $(3,3)$ in
$\PP^2\times\PP^2$.
The factors `sign' and  `degree' are needed to match the results of
our calculations with
the tables of enumerative data in the literature.
Moreover, if the wrong sign is used, the numbers $N_{j_1,\ldots,j_d}$
are often even not integers.

\

\noindent
\textbf{Conclusion.} \emph{The above discussion shows that the canonical
coordinates and the prepotential act like a parametrization for the
image of the Schwarz map: the image points have coordinates
$(1,t_1,\ldots,t_{2d+1})$ with}
$$
\begin{array}{rclll}
t_j&=&\log q_j&\textrm{for}& j=1,\ldots,d\,,\\[1ex]
t_{d+j}&=&q_j\frac{\partial\cF}{\partial q_j}
&\textrm{for}& j=1,\ldots,d\,,\\[1ex]
t_{2d+1}&=&-2\,\cF\,+\,
\sum_{j=1}^d t_jt_{d+j}\,.
\end{array}
$$

\

\noindent
\textbf{Remark.}
On the graded ring $\oR_{\gA,\cC}$ there is an involution ${.}^*$
given for homogeneous elements by
$x^*\,=\,(-1)^{\mathrm{deg}\,x}x$.
We fix the linear map
$\tau: \oR_{\gA,\cC}\stackrel{\mathrm{project}}{\longrightarrow}
\oR_{\gA,\cC}^{(3)}\stackrel{\simeq}{\longrightarrow}\ZZ$.
Then, in the examples of Sections \ref{subsection:quintic}, \ref{subsection:cubics intersect}, \ref{subsection:(3,3)}
the alternating bilinear form defined by the ordered basis and
the relations (\ref{eq:intersection products}), is in fact
$$
\langle x,y\rangle\,=\, \tau (x^*y)\,.
$$
Moreover, in those examples the trick of expanding 
$\log\Phi_{\LL,\ugamma,\ueps}(\vu)$ in two ways showed that the
prepotential is a polynomial of degree $3$ in 
$\log q_1,\ldots,\log q_d$ plus a power series in $q_1,\ldots,q_d$.

\

\noindent
\textbf{Example.} Let us check that (\ref{eq:intersection products})
holds for the ordered basis 
$$
\{e_0,e_{1,1},e_{1,2},e_{2,1},e_{2,2},e_3\}\,:=\,
\{1,\ode_1,\,\ode_2,\,\ode_2^2,\,\ode_1^2,\,\ode_1^2\ode_2\}
$$
in the example of Section \ref{subsection:(3,3)}.
The alternating bilinear form (on the third cohomology space in a family of Calabi-Yau threefolds) is invariant under the local monodromy
operators, which in this case are given by multiplication by
$\exp(2\pi \ci \ode_1)$ and $\exp(2\pi \ci \ode_2)$.
This means that the matrices for multiplication with $\ode_1$ and 
$\ode_2$ and the Gramm matrix of the alternating bilinear form,
everything with respect to the above basis, must satisfy
$\mathrm{Gramm}^t\,=\,-\mathrm{Gramm}$ and
\begin{equation}\label{eq:gramm}
\mathrm{mat}(\ode_a)\cdot\mathrm{Gramm}\,=\,
-\mathrm{Gramm}\cdot\mathrm{mat}(\ode_a)^t\,,\qquad a=1,2\,.
\end{equation}
One easily checks
$$
\mathrm{mat}(\ode_1)\,=\,
\left(\begin{array}{rrrrrr}
0&0&0&0&0&0\\
1&0&0&0&0&0\\
0&0&0&0&0&0\\
0&0&1&0&0&0\\
0&1&1&0&0&0\\
0&0&0&1&0&0
\end{array}\right)\,,\qquad
\mathrm{mat}(\ode_2)\,=\,
\left(\begin{array}{rrrrrr}
0&0&0&0&0&0\\
0&0&0&0&0&0\\
1&0&0&0&0&0\\
0&1&1&0&0&0\\
0&1&0&0&0&0\\
0&0&0&0&1&0
\end{array}\right)\,.
$$
The general anti-symmetric matrix solution to (\ref{eq:gramm})
has, up to a non-zero scalar factor, the form
$$
\left(\begin{array}{rrrrrr}
0&0&0&0&0&1\\
0&0&0&-1&0&0\\
0&0&0&0&-1&0\\
0&1&0&0&0&x\\
0&0&1&0&0&y\\
-1&0&0&-x&-y&0
\end{array}\right)\,.
$$
When we evaluate (\ref{eq:transversality}) using this Gramm matrix,
we find, for $i=1,2$,
\begin{equation}\label{eq:transversality2}
\begin{array}{l}
\displaystyle{
q_i\frac{\partial\gf_3}{\partial q_i}\,-\,
\gf_{1,1} q_i\frac{\partial\gf_{2,1}}{\partial q_i}
\,-\,
\gf_{1,2} q_i\frac{\partial\gf_{2,2}}{\partial q_i}\,+\,
\gf_{2,i}\,+\,}
\\
\displaystyle{
+\,
x\left(\gf_{2,1}q_i\frac{\partial\gf_3}{\partial q_i}
-\gf_3q_i\frac{\partial\gf_{2,1}}{\partial q_i}\right)
\,+\,
y\left(\gf_{2,2}q_i\frac{\partial\gf_3}{\partial q_i}
-\gf_3q_i\frac{\partial\gf_{2,2}}{\partial q_i}\right)
\,=\,0\,.}
\end{array}
\end{equation}
We want to estimate the growth of the various terms in this
expression by looking at the logarithmic pieces. 
So, recall that in this example
$$
\Phi_{\LL,\ugamma,\ueps}(\vu)\,=\,
\Phi_0(\vu)z_1^{\ode_1}z_2^{\ode_2}\times
(\textrm{power series in }z_1,z_2)
$$
and
\begin{eqnarray*}
z_1^{\ode_1}z_2^{\ode_2}&=&
(1+(\log z_1)\ode_1+\textstyle{\frac{1}{2}}(\log z_1)^2\ode_1^2)
(1+(\log z_2)\ode_2+\textstyle{\frac{1}{2}}(\log z_2)^2\ode_2^2)
\\
&=& 1\,+\,(\log z_1)\ode_1\,+\,(\log z_2)\ode_2\,+\,
\left(\textstyle{\frac{1}{2}}(\log z_2)^2 +(\log z_1)(\log z_2)\right)
\ode_2^2
\\
&&+\,
\left(\textstyle{\frac{1}{2}}(\log z_1)^2 +(\log z_1)(\log z_2)\right)
\ode_1^2\,+\,\textstyle{\frac{1}{2}}
(\log z_1)(\log z_2)(\log z_1z_2)\ode_1^2\ode_2\,.
\end{eqnarray*}
Moreover, up to addition of power series, $\log q_1\,\asymp\,\log z_1$
and $\log q_2\,\asymp\,\log z_2$.
Thus we see that the highest order logarithmic contributions are
\begin{eqnarray*}
-\gf_3q_1\frac{\partial\gf_{2,1}}{\partial q_1}\,+\,
\gf_{2,1}q_1\frac{\partial\gf_3}{\partial q_1}&\asymp&
-\,
\textstyle{\frac{1}{2}}((\log q_1)^2(\log q_2)+(\log q_1)(\log q_2)^2)
(\log q_2)
\\
&&\hspace{-6em}
+\,\left(\textstyle{\frac{1}{2}}(\log q_2)^2 +(\log q_1)(\log q_2)
\right)
\left((\log q_1)(\log q_2)+\textstyle{\frac{1}{2}}(\log q_2)^2\right)
\\
&&\hspace{-7em}=\:
\textstyle{\frac{1}{2}}(\log q_1)^2(\log q_2)^2\,+\,
\textstyle{\frac{1}{2}}(\log q_1)(\log q_2)^3\,+\,
\textstyle{\frac{1}{4}}(\log q_2)^4
\end{eqnarray*}
and
\begin{eqnarray*}
-\gf_3q_1\frac{\partial\gf_{2,2}}{\partial q_1}\,+\,
\gf_{2,2}q_1\frac{\partial\gf_3}{\partial q_1}&&
\\
&&\hspace{-7em} \asymp\,
-\,
\textstyle{\frac{1}{2}}((\log q_1)^2(\log q_2)+(\log q_1)(\log q_2)^2)
(\log q_1+\log q_2)
\\
&&\hspace{-6em}
+\,\left(\textstyle{\frac{1}{2}}(\log q_1)^2 +(\log q_1)(\log q_2)
\right)
\left((\log q_1)(\log q_2)+\textstyle{\frac{1}{2}}(\log q_2)^2\right)
\\
&&\hspace{-7em} =\,\textstyle{\frac{1}{4}}(\log q_1)^2(\log q_2)^2
\end{eqnarray*}
So if we consider (\ref{eq:transversality2}) for $i=1$ and 
$q_2\downarrow 0$ the dominant $\textstyle{\frac{1}{4}}(\log q_2)^4$
term forces $x=0$. Having $x=0$ we consider (\ref{eq:transversality2}) 
for $i=1$ and $q_1=q_2\downarrow 0$. Once again there is a dominant $\textstyle{\frac{1}{4}}(\log q_2)^4$, forcing $y=0$.
\\
\emph{This finishes the proof for the fact that in the example
of Section \ref{subsection:(3,3)} the ordered basis
$\{1,\ode_1,\,\ode_2,\,\ode_2^2,\,\ode_1^2,\,\ode_1^2\ode_2\}$
satisfies (\ref{eq:intersection products}).}

\

\noindent
\textbf{Exercise.} Apply the techniques of the preceding example
and show that the basis $\{1,\oeps,\oeps^2,\oeps^3\}$
in Sections \ref{subsection:quintic} 
and \ref{subsection:cubics intersect}
satisfies (\ref{eq:intersection products}).


%

\end{document}